\documentclass[runningheads]{llncs}
\usepackage[T1]{fontenc}
\usepackage{graphicx}
\usepackage{amsmath}
\usepackage{caption}
\usepackage[pass]{geometry}
\usepackage{latexsym}
\usepackage{placeins}
\usepackage[figuresright]{rotating}
\usepackage{setspace}

\numberwithin{equation}{section}
\numberwithin{figure}{section}
\numberwithin{table}{section}
\numberwithin{theorem}{section}
\numberwithin{lemma}{section}
\numberwithin{corollary}{section}

\DeclareSymbolFont{AMSb}{U}{msb}{m}{n}
\DeclareMathSymbol{\R}{\mathalpha}{AMSb}{"52}
\DeclareMathSymbol{\Z}{\mathalpha}{AMSb}{"5A}
\DeclareMathSymbol{\MetricM}{\mathalpha}{AMSb}{"4D}

\newcommand\STHGP[1]{{\mathrm{STHGP(#1)}}}

\newcommand\BothMode{{\ttfamily -{}-{}both}}
\newcommand\LargeMode{{\ttfamily -{}-{}large}}
\newcommand\SmallMode{{\ttfamily -{}-{}small}}

\spnewtheorem{MyLemma}[theorem]{Lemma}{\bfseries}{\itshape}
\spnewtheorem*{MyProof}{Proof}{\bfseries}{\upshape}
\spnewtheorem{MyCorollary}[theorem]{Corollary}{\bfseries}{\itshape}

\newlength{\ALindent}
\newcommand\AL[3][]{\makebox[\ALindent][l]{#1}\hskip #2\ALindent plus 0pt minus 0pt#3\hfill\break\vspace{-2pt}}

\newcommand{\MyBeginFig}{\begin{figure}[htbp]\centering}
\newcommand{\MyBeginFigure}[1]{\begin{figure}[#1]\centering}
\newcommand{\MyEndFig}{\end{figure}}
\newcommand{\MyBeginTable}{\begin{table}[htbp]\centering}
\newcommand{\MyEndTable}{\end{table}}

\newcommand{\MyBeginMeanTable}{\begin{table}[htbp]\centering}
\newcommand{\MyEndMeanTable}{\end{table}}

\newcommand{\MyBeginSummaryTable}{\begin{sidewaystable}[htbp]\centering}
\newcommand{\MyEndSummaryTable}{\end{sidewaystable}}

\newcommand{\MyBeginPlot}{\begin{figure}[h]\centering}
\newcommand{\MyEndPlot}{\end{figure}}

\newcommand{\MyBeginFSTTable}{\begin{table}[htbp]\centering}
\newcommand{\MyEndFSTTable}{\end{table}}

\newcommand{\MyBeginConcatTable}{\begin{sidewaystable}[htbp]\centering}
\newcommand{\MyEndConcatTable}{\end{sidewaystable}}

\newcommand{\SDEV}[1]{{\tiny $\pm #1$}}
\newcommand{\NEWLINE}{\vspace*{-0.2cm}\\}

\newcommand{\RunningFlag}{0}
\newcommand{\ClrRunningFlag}{\global\def\RunningFlag{0}}

\newcommand{\CrashedFlag}{0}
\newcommand{\ClrCrashedFlag}{\global\def\CrashedFlag{0}}

\newcommand{\UnprunedFlag}{0}
\newcommand{\ClrUnprunedFlag}{\global\def\UnprunedFlag{0}}
\newcommand{\SetUnprunedFlag}{\global\def\UnprunedFlag{1}}

\newcommand{\UNPRUNED}{$\strut^{\ddagger}$\SetUnprunedFlag}

\newcommand{\ExplainRunning}{\if 1\RunningFlag{} $*~=$ still running. \fi}
\newcommand{\ExplainCrashed}{\if 1\CrashedFlag{} $\dagger~=$ crashed. \fi}
\newcommand{\ExplainUnpruned}{%
	\if 1\UnprunedFlag{} $\ddagger~=$ used unpruned FSTs. \fi}

\newcommand{\InitializeTableFootnoteFlags}{%
	\ClrRunningFlag%
	\ClrCrashedFlag%
	\ClrUnprunedFlag%
}
\newcommand{\DisplayTableFootnoteFlags}{%
	\strut%
	\ExplainRunning%
	\ExplainCrashed%
	\ExplainUnpruned%
}

\begin{document}
\title{Quantitative Indicators for Strength of Inequalities with
  Respect to a Polyhedron}%
\subtitle{Part II: Applications and Computational Evidence}
\titlerunning{Strength of Inequalities with Respect to a Polyhedron, Part II: Applications}
% If the paper title is too long for the running head, you can set
% an abbreviated paper title here
%
\author{David M. Warme%
%\inst{1}
\orcidID{0009-0004-2307-9812}}
\authorrunning{D. M. Warme}
% First names are abbreviated in the running head.
% If there are more than two authors, 'et al.' is used.
%
\institute{Group W, Vienna VA 22180, USA
\email{dwarme@groupw.com}%
%\url{http://www.springer.com/gp/computer-science/lncs}
}%
\maketitle

\begin{abstract}

``Strength'' is an important property of inequalities used in integer
and mixed-integer optimization, both in theory and practice.
Unfortunately, no good formal characterization for strength exists,
nor is it well-understood.

The first paper explored two quantitative strength indicators (extreme
point ratio (EPR) and centroid distance (CD)), applying them
to the subtour inequalities of the spanning tree in hypergraph
polytope $\STHGP{n}$.
Although it was known that subtour inequalities of small cardinality
were strong, EPR and CD agree that subtour inequalities of {\em large}
cardinality are significantly {\em stronger}.

In this second paper, we exploit this previously unknown property
algorithmically, presenting strong computational
evidence that the EPR and CD indicators are highly predictive of
actual computational strength.
Previous branch-and-cut implementations for optimizing over
$\STHGP{n}$ find violated subtour inequalities of only relatively
small cardinality, strengthening only by reducing the cardinality of
violated subtours.
We present new methods that strengthen violated subtour inequalities
by {\em augmentation} (instead of reduction).
Combining strengthening via reduction and augmentation yields violated
subtour inequalities of both small and large cardinality, covering
both classes deemed ``strong'' by EPR and CD.
Across all instance classes studied, the
computational results are remarkable --- culminating with an optimal
solution of a 1,000,000 terminal random Euclidean Steiner tree
instance.

The conclusion is that the EPR and CD strength indicators presented in
the first paper have strong predictive power regarding actual
computational strength (at least regarding $\STHGP{n}$ subtour
inequalieis).
The ability to accurately measure the strength of inequalities has
numerous applications of great importance, both in theory and
practice.

\keywords{Integer programming	\and%
  Combinatorial optimization	\and%
  Steiner tree			\and%
  Enumerative Combinatorics	\and%
  Extreme point ratio		\and%
  Centroid distance		\and%
  GeoSteiner}%
{%
\renewcommand{\keywordname}{{\bf Mathematics Subject Classification (2020):}}
\keywords{90C10 \and 90C11 \and 90C27 \and 05A15}
}%
\end{abstract}

\section{Introduction}
\label{sec:intro}

In the first paper \cite{WarmeIndicators1}, we reviewed several
indicators from the literature that we claim to be strongly correlated
with computational ``strength'' of inequalities with respect to a
given polyhedron.
These indicators are functions $m : (H,P) \mapsto \R^+$, where $H$ is
the hyperplane corresponding to a given linear inequality, and $P$ is
a polyhedron.
Let $X$ be the set of all extreme points of polyhedron $P$.
The extreme point ratio (EPR) indicator is $|H \cap X| \, / \, |X|$,
for which larger ratios indicate stronger inequalities.
The centroid distance (CD) indicator measures the distance between the
centroid $C = (1 \, / \, |X|) \, \sum_{x \in X} x$ and $H$, for which
shorter distances indicate stronger inequalities.
(There are two versions of CD indicator: the ``weak'' CD measures
distance along a line residing entirely within the affine hull of $P$.
The ``strong'' CD measures distance along a line residing entirely
within $P$.)

\begin{figure}[!ht]
\begin{center}
\begin{minipage}[t]{2.25in}
\includegraphics[width=2.25in,clip=]{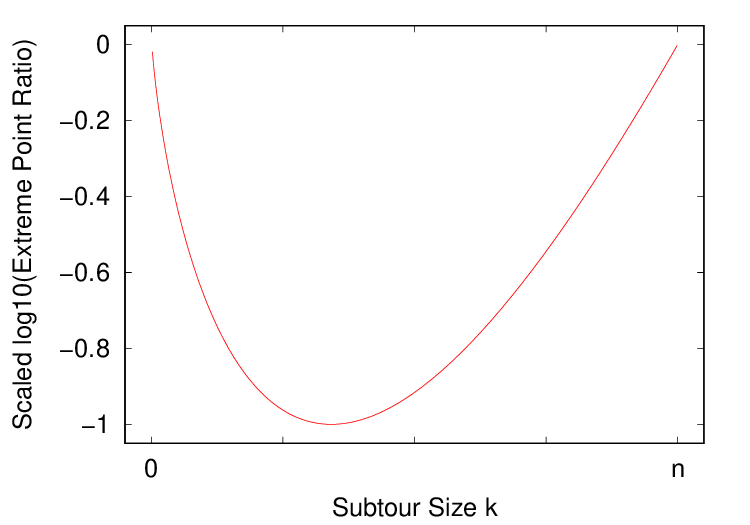}
 \captionof{figure}{EPR for STHGP subtours of size
  $k=|S|$, $n=1000$ (larger = stronger)}
 \label{fig:sthgp-sec-epr-indicators-1000}
\end{minipage}
\quad
\begin{minipage}[t]{2.25in}
\includegraphics[width=2.25in,clip=]{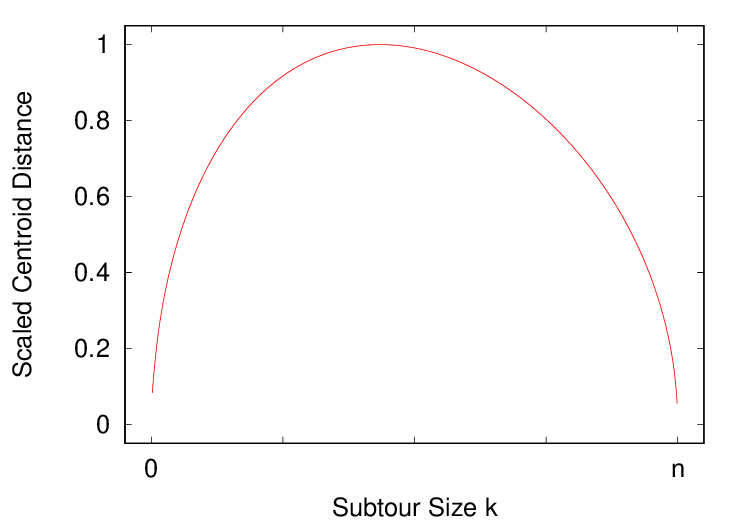}
 \captionof{figure}{CD for STHGP subtours of size
  $k=|S|$, $n=1000$ (smaller = stronger)}
 \label{fig:sthgp-sec-cd-indicators-1000}
\end{minipage}
\end{center}
\vspace*{\floatsep}
\begin{center}
\begin{minipage}[t]{2.25in}
\includegraphics[width=2.25in,clip=]{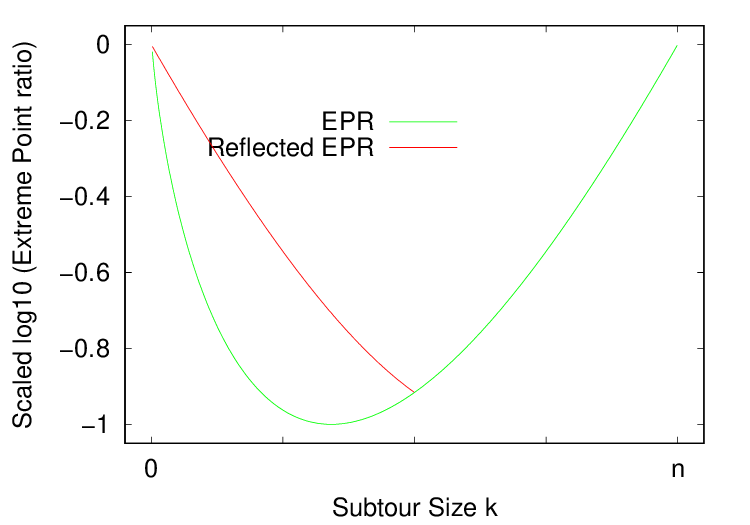}
 \captionof{figure}{EPR for STHGP subtours showing
   large cardinality subtours are stronger than those of small
   cardinality (larger = stronger)}
 \label{fig:reflected-epr-indicators-1000}
\end{minipage}
\quad
\begin{minipage}[t]{2.25in}
\includegraphics[width=2.25in,clip=]{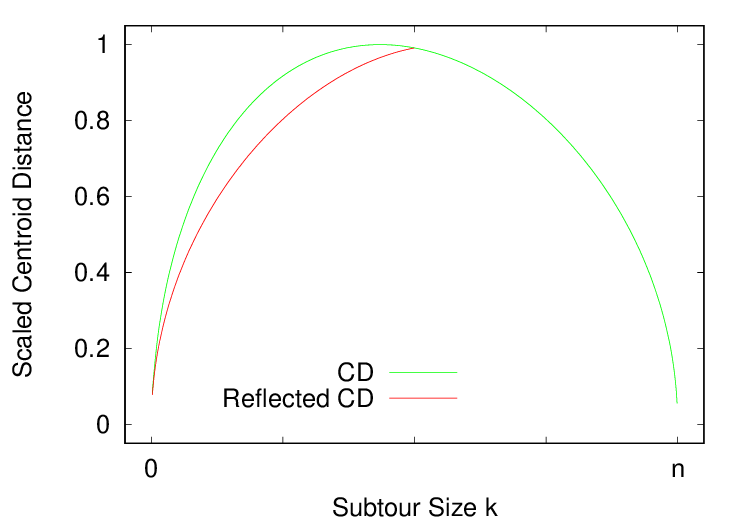}
 \captionof{figure}{CD for STHGP subtours showing
   large cardinality subtours are stronger than those of small
   cardinality (smaller = stronger)}
 \label{fig:reflected-cd-indicators-1000}
\end{minipage}
\end{center}
\end{figure}

Closed-forms were obtained in \cite{WarmeIndicators1} for both the EPR
and CD indicators applied to the subtour inequalities of the Spanning
Tree in Hypergraph Polytope (STHGP).
Figures~\ref{fig:sthgp-sec-epr-indicators-1000}
and~\ref{fig:sthgp-sec-cd-indicators-1000}
present plots of the EPR and CD indicators for STHGP subtours
$S \subset V$
as a function of $k=|S|$ with $n=1000$.
It was known that subtours of small cardinality were computationally
very strong, and both indicators confirm this.
Both indicators, however, also agree that subtours of large
cardinality (i.e., containing a large fraction of all vertices) are
even {\em stronger} than subtours of small cardinality.
This is illustrated in
Figures~\ref{fig:reflected-epr-indicators-1000}
and~\ref{fig:reflected-cd-indicators-1000}
in which the right half of
each curve has been reflected onto the left half.
The region between the two curves indicates how much stronger the
large subtours are predicted to be than the corresponding small
subtour.
This property was neither known, nor suspected until revealed by the
indicators, and suggests that adding violated subtours $S$ to the LP
where $|S|$ is large could significantly speed up computations by
reducing the number of optimize/separate iterations needed to satisfy
all of the subtour inequalities.

In this paper we test these ideas computationally within GeoSteiner,
which has been the fastest software in existence for solving geometric
Steiner trees in the plane since 1998
\cite{Warme98},
\cite{JuhlWarmeWinterZachariasen2018}.
Given a finite set $N \subseteq {\R^{2}}$ (called ``terminals''),
the geometric Steiner tree problem seeks a set of line segments
connecting $N$ having minimal total length.
The solution usually includes additional points (called ``Steiner
points'') introduced as necessary, serving as ``junctions'' that
reduce the total length.
A distance metric must be specified, giving rise to different Steiner
tree variants (Euclidean, rectilinear, etc.).
The Euclidean Steiner tree problem is one of the oldest optimization
problems in all of
mathematics~\cite{BrazilGrahamThomasZachariasen2014}.

A closely related problem is the Steiner problem in graphs (SPG).
Given an edge weighted graph $G=(V,E)$ with vertices partitioned into
terminals $T$ and Steiner vertices $S$, find edges $E' \subseteq E$
forming a tree of minimum total weight connecting $T$.
Geometric Steiner tree problems can be reduced to instances of SPG.
Polzin and Vahdati Daneshmand~\cite{Polzin2003},
\cite{VahdatiDaneshmand2004},
\cite{PolzinVahdatiDaneshmandDIMACS2014}
long had the best methods for solving SPG.
Recently, Rehfeldt and Koch~\cite{RehfeldtKoch2023} reported important
progress on SPG.

Significantly condensed versions of this paper and its
companion~\cite{WarmeIndicators1} have been submitted for journal
publication.

We define our notation in Section~\ref{sec:notation}.
In Section~\ref{sec:geosteiner-approach} we outline the major
algorithmic phases of GeoSteiner.
In Section~\ref{sec:geosteiner-forces}
we explore the algorithmic forces operating within GeoSteiner that
cause it to effectively find only violated subtours $S$ of relatively
small cardinality.
In Section~\ref{sec:geosteiner-mods}
we present all of the modifications to the experimental version of
GeoSteiner presented here, including the subtour separation and
constraint strengthening algorithms that intentionally target subtours
of large cardinality to exploit their superior strength, as predicted
by both EPR and CD indicators.
% No motivation for such modifications existed prior to the insight
% these indicators provided.
%
In Section~\ref{sec:computational-results}
we present computational results for the improved version of
GeoSteiner.  These results provide truly remarkable improvements over
the previous algorithm~\cite{JuhlWarmeWinterZachariasen2018},
including an optimal solution for a random 1,000,000 terminal
Euclidean instance.
We also present computations specifically designed to directly compare
the strength of subtours having large cardinality versus those of
small cardinality (and combining both together).
Finally, we describe computational experiments for the Steiner tree in
hypergraph problem (having nearly identical EPR and CD curves):
targeting subtours of both small and large cardinality produces
similar speedups on this problem.
This provides compelling computational evidence that the EPR and CD
indicators are highly predictive of actual computational strength, at
least regarding the subtours of STHGP.
In Section~\ref{sec:discussion} we discuss various ramifications and
applications of these results.
We conclude in Section~\ref{sec:conclusion}.
Detailed computational results appear in Appendix~\ref{sec:appendix}.

\section{Notation}
\label{sec:notation}

We define $\R$ to be the real numbers, and $\Z$ to be the integers.
The corresponding non-negative sets are denoted
$\Z^+$ and $\R^+$, respectively.
A hypergraph $H = (V,E)$ consists of a finite set $V$ of vertices,
and a set $E$ of edges.
Each edge $e \in E$ satisfies $e \subseteq V$ and $|e| \ge 2$.
For any $S \subseteq V$, let:

\begin{displaymath}
	\delta(S) ~=~
	  \lbrace e \in E : \emptyset \subset (e \cap S) \subset e \rbrace.
\end{displaymath}
We extend this definition to also apply to single vertices so that for
any $v \in V$, $\delta(v) = \delta(\lbrace v \rbrace)$.
Let $S,T \subseteq V$ such that $S \cap T = \emptyset$.
We define
\begin{displaymath}
	(S:T) = \lbrace e \in E :
		(e \cap S \ne \emptyset) \wedge
		(e \cap T \ne \emptyset) \rbrace.
\end{displaymath}
Let $x \in \R^{|E|}$.
For any $F \subseteq E$, $x(F)$ denotes $\sum_{e \in F} x_e$.
For $S,T \subseteq V$ such that $S \cap T = \emptyset$ we define
\begin{eqnarray*}
	x(S) &\equiv& \sum_{e \in E} \max(|e \cap S| \,-\, 1,\, 0) \, x_e, \\
	x(S:T) &\equiv& \sum_{e \in (S:T)} x_e, \\
	x(\delta(S)) &\equiv& x(S:V-S).
\end{eqnarray*}
The first of these gives the total weight of edges within $S$
(where edge $e$ having $k = |e \cap S|$ vertices within $S$ is taken
to be equivalent to $\max(k-1,0)$ conventional edges having the same
weight $x_e$).
The second of these gives the weight of edges crossing cut $(S:T)$.
% Note that $x(S:T) \,=\, x(T:S).$
For $v \in V$, we also define
$x(\delta(v)) \,\equiv\, x(\delta(\lbrace v \rbrace))$.

Given a hypergraph $H=(V,E)$ and an $x \in R^{|E|}$, we define the
{\em support hypergraph} of $H$ to be the hypergraph
$\bar{H} = (V,\bar{E})$,
where $\bar{E} \,=\, \lbrace e \in E : x_e > 0 \rbrace$.

\section{GeoSteiner Algorithmic Overview}
\label{sec:geosteiner-approach}

GeoSteiner computes a geometric Steiner tree in two main phases.
Given a finite set $N$ of terminals in $\R^2$, the ``FST Generation''
phase computes a ``sufficient'' set $F$ of Full Steiner Trees (FSTs).
An FST $t$ is a geometrically embedded tree such that:
(1) $t$ connects some subset $U \subseteq N$;
(2) terminals appearing in $t$ are all leaves of $t$ (degree 1);
(3) all non-leaf vertices in $t$ (Steiner vertices) are interior
nodes of $t$ (degree $>$ 1); and
(4) $t$ is locally optimal.
(This is an approximation: FSTs satisfy additional conditions.)
The set $F$ is ``sufficient'' in that it guarantees that an optimal
Steiner tree exists as the concatenation of some subset
$F' \subseteq F$ of FSTs.
A completely different algorithm is used for each supported distance
metric $\MetricM$ (based upon a corresponding structure theorem
for $\MetricM$),
and implicitly considers all possible subsets $U$ of $N$, and all
possible tree topologies over $U$, discarding only combinations that
can be proven to be sub-optimal or redundant.
Although this can give rise to an exponentially large set of FSTs,
most point sets ``in general position'' yield only $O(n)$ FSTs for $n$
given terminals.
The FST generation algorithms are described in:
Euclidean~\cite{WinterZachariasen},
rectilinear~\cite{ZachariasenFST},
and uniform directions (hexagonal and
octilinear)~\cite{BrazilThomasWengZachariasen2006}.
Despite the large space that is implicitly enumerated, GeoSteiner's
FST generation algorithms are surprisingly efficient.
For ``well behaved'' point sets yielding a linear number of FSTs,
these algorithms empirically run in approximately $O(n^{2.5})$ time.
This success notwithstanding, there are point sets that do provoke bad
behavior.
For the Euclidean distance metric, points lying on regular integer
grids give rise to an exponential number of
FSTs~\cite{CockayneHewgillImproved},%
~\cite{WinterZachariasen}.
For the rectlinear metric there is a known ``fractal'' pattern of
terminals that yields an exponential number of
FSTs~\cite{FossmeierKaufmann}.

During the second ``FST concatenation'' phase, we are given the
terminals $N$ and the set $F$ of FSTs, and the goal is to identify a
subset $F' \subseteq F$ such that the FSTs $F'$ form a tree that spans
all terminals $N$ having minimum total length.

Between these two phases is an optional {\em FST pruning} phase that
uses global information about the entire set of FSTs to mark some FSTs
as {\em pruned} or {\em required}.
A {\em pruned} FST cannot appear in any optimal Steiner tree for $N$
whereas a {\em required} FST {\em must} appear in every optimal
Steiner tree for $N$.
FST pruning can greatly reduce and simplify the subsequent FST
concatenation problem~\cite{JuhlWarmeWinterZachariasen2018}.

\begin{figure}[ht]
 \centering
 \includegraphics[width=4.0in]{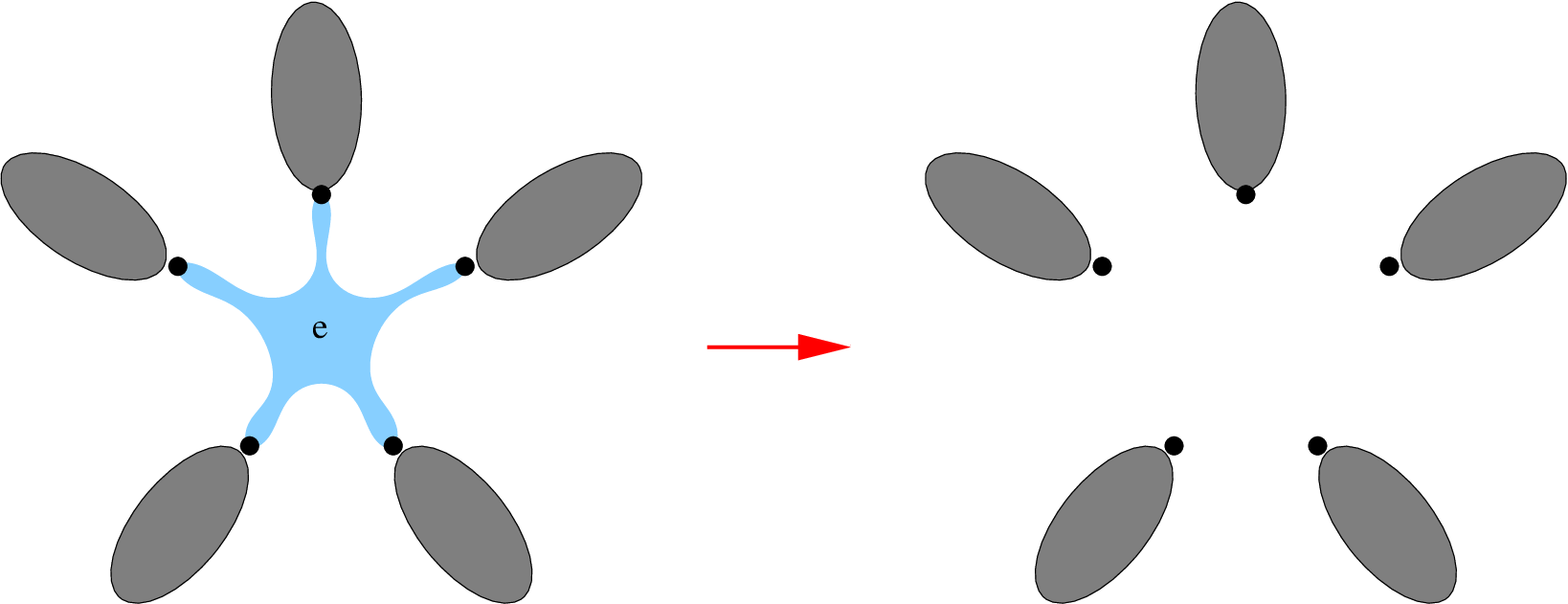}
 \caption{Removing edge $e$ with $k=|e|$ increases the number of
   connected components by $k-1$ ($k=5$ in this example)}
 \label{fig:reverse-kruskal}
\end{figure}

\begin{figure}[ht]
 \centering
 \includegraphics[width=2.5in]{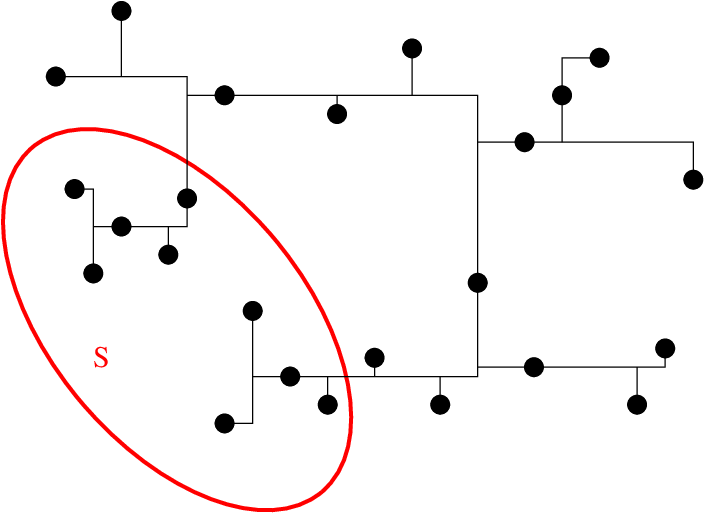}
 \caption{Subhypergraph of a tree induced by $S \subseteq V$, where
   $|S|=k$, $k \ge 2$ yields a forest of $j$ trees, $1 \le j \le k$
   ($j=2$ and $k=9$ in this example)}
 \label{fig:fig-forest}
\end{figure}

GeoSteiner reduces the FST concatenation phase to the problem of
finding a minimum spanning tree in a hypergraph $H = (V,E)$ with given
vector $c$ of edge weights.
The hypergraph is constructed by trivially setting $V = N$, and for
each FST $t \in F$, there is a corresponding edge $e \in E$ whose
vertices consist of the terminals appearing in $t$, and whose edge
weight $c_e$ is the geometric length of $t$, as measured in the
appropriate distance metric.
(The edges $E$ contain only terminals --- the interior Steiner
vertices of the FST are abstracted away in the resulting hypergraph.)
GeoSteiner then uses branch-and-cut to find a Minimum Spanning Tree in
Hypergraph $H$ using the following integer program:
\begin{align}
\nonumber
&\mathrm{Minimize:} \\
\nonumber
&	&&	c \cdot x \\
\nonumber
&\mathrm{Subject~To:} \\
\label{eq:total-degree-equation}
&	&& x(V) ~=~ |V| \, - \, 1, \\
\label{eq:subtour-inequalities}
&	&& x(S)	~\le~ |S| \,-\, 1
	&& \text{for all $S \subset V$ such that $|S| \ge 2$}, \\
\nonumber
&	&& x_e \in \lbrace 0,\, 1 \rbrace
	&& \text{for all $e \in E$}.
\end{align}
Let $T \subseteq E$ be such that $H_T = (V,T)$ is a spanning tree of
$H$.
($H_T$ consists of a single connected component with no cycles.)
Equation~(\ref{eq:total-degree-equation})
satisfies all incidence vectors $x$ corresponding to such spanning
trees $T$ as we now show.
Consider a process in which all edges $e \in T$ are removed one at a
time.
Removing an edge $e$ of cardinality $k=|e|$ replaces the connected
component containing $e$ with $k$ connected components (as shown in
Figure~\ref{fig:reverse-kruskal}), increasing the number of connected
components by $k-1$.
When all edges have been removed, the number of connected components
has increased from 1 to $|V|$, for a total increase of $|V|-1$
connected components.
Thus equation~(\ref{eq:total-degree-equation}) holds for every
incidence vector $x$ corresponding to a spanning tree of $H$.
\cite{Warme98} (page 43) shows this to be the only such equation, so
that~(\ref{eq:total-degree-equation}) is the affine hull of
$\STHGP{n}$.
Consider the subhypergraph induced by any $S \subseteq V$, where
$|S|=k$, $k \ge 2$.
The result is a forest of $j$ trees, where $1 \le j \le k$, as shown
in Figure~\ref{fig:fig-forest}.
Applying the same counting argument to this induced subhypergraph
yields $x(S) \,=\, |S| \,-\, j$, whose right-hand-side is maximized
when $j=1$.
This implies~(\ref{eq:subtour-inequalities}) is satisfied by all
indicence vectors $x$ that are spanning trees of $H$.
It is proven in~\cite{Warme98} that the integer feasible solutions of
this integer program are identically the incidence vectors of all
spanning trees of $H$.
This correspondence is independent of $c$, which can therefore be
any vector in $\R^{|E|}$.

\section{GeoSteiner Finds Only Small Cardinality Subtours}
\label{sec:geosteiner-forces}

The EPR and CD metrics predict that subtours of large cardinality are
stronger than the corresponding subtours of small cardinality.
This section explains why GeoSteiner versions 5.3 and earlier
generally use only subtours of relatively small cardinality, thereby
missing out on potential performance gains that might be obtained by
also using subtours of large cardinality.

We first note that there is actually one special case in which
GeoSteiner versions 5.3 and earlier do find violated subtours having
large cardinality.
Let $\bar{H} = (V,\bar{E})$ be the support hypergraph, and
$H_1 = (V_1,E_1),\, H_2 = (V_2,E_2),\,
\ldots,\,
H_k = (V_k,E_k)$ be the connected components of $\bar{H}$, with $k \ge
2$.
For each $1 \le j \le k$ GeoSteiner generates subtour inequalities for
both $V_j$ and $V-V_j$ (both likely violated).
When $|V_j|$ is small, $|V-V_j|$ is large, and vice versa.
In most cases the support hypergraph has $k \ge 2$ connected
components for only the first several optimize / separate iterations,
limiting the effectiveness of this technique.

GeoSteiner's main separation algorithm for subtour inequalities is
described in \cite{Warme98}, pages 70--74.
Given a hyperaph $H = (V,E)$ with edge weights $x$
(and corresponding support hypergraph $\bar{H} = (V,\bar{E})$),
it finds an $S \subset V$ such that subtour $S$ is maximally
violated by $x$
by running at most $|V| - 1$ max-flow/min-cut computations on an
acyclic flow network having $|V| + |\bar{E}| + 2$ vertices and
$|V| \, + \, |\bar{E}| \, +\, \sum_{e \in \Bar{E}} |e|$ directed
arcs.

For all of the distance metrics supported by GeoSteiner, there is a
small, fixed integer $k$ (e.g., $k = 20$ usually suffices) such that
the number of terminals in a single FST is at most $k$ (with high
probability) over all ``well behaved'' point sets
(i.e., ``in general position'').
Furthermore, as mentioned previously, we usually have
$|E| \,=\, O(|V|)$ (which implies $|\bar{E}| \, = \, O(|V|)$).
This means that the flow network is usually of size $O(n)$, where
$n = |V|$.
Under these conditions the separation algorithm runs in $O(|V|^4)$ time
(using the strongly-polynomial $O(|V| \, |E|^2)$ Dinitz
algorithm~\cite{Dinitz1970}).
When using an algorithm this costly, it is advantageous to perform
various reductions on the support hypergraph $\bar{H}$ to assure that
the expensive flow separator is invoked only on instances that are as
small as possible.

For all $v \in V$, we define
\begin{displaymath}
	b_v ~=~ x(\delta(v))
\end{displaymath}
For all $S \subset V$ with $|S| \ge 2$ we define
$f(S) \,=\, |S| \,-\, x(S)$,
so that $f(S) < 1$ implies that the corresponding subtour
inequality~(\ref{eq:subtour-inequalities}) is violated.

The following reduction lemmas were proved in~\cite{Warme98}:
\begin{MyLemma}
\label{lem:reduction}
Let $x \in \R^{|E|}$,
$U \subset V$ such that $|U| \ge 2$, and
$v \in V-U$.
If $b_v \le 1$ and
$f(U \cup \lbrace v \rbrace) < 1$
then $f(U) \le f(U \cup \lbrace v \rbrace) < 1$.
\end{MyLemma}
Thus vertices $v$ such that $b_v \le 1$ can be iteratively eliminated
from the support hypergraph until no such vertices remain.
If the original hypergraph contains a violated subtour, then the
reduced hypergraph must also contain a violated subtour.

\begin{MyLemma}
\label{lem:reduce-cc}
Let $H = (V, E)$ be a hypergraph with weights $x_e$ for all $e \in E$.
Let $\bar H = (V, \bar E)$ be the support hypergraph of $H$.
Let the connected components of $\bar H$ be
$H_1 = (V_1,E_1),\, H_2 = (V_2,E_2),\, \ldots,\, H_k = (V_k,E_k)$.
Let $S \subseteq V$ and
$S_j = S \cap V_j$ for all $1 \le j \le k$.
If $f(S) < 1$ then there is some $j$ such that $f(S_j) < 1$.
\end{MyLemma}
It therefore suffices to independently seperate the connected
components of the support hypergraph.

\begin{MyLemma}
\label{lem:reduce-bcc}
Let $H = (V, E)$ be a hypergraph with weights $x_e$ for all $e \in E$.
Let $\bar H = (V, \bar E)$ be the support hypergraph of $H$.
Let $A,B,C$ be a partition of $V$
and $E_A,E_B$ be a partition of $\bar E$
such that $|C|=1$,
$E_A = \lbrace e \in \bar E : e \subseteq (A \cup C) \rbrace$ and
$E_B = \lbrace e \in \bar E : e \subseteq (B \cup C) \rbrace$.
If $S \subset V$ such that $f(S) < 1$
then $f(S \cap (A \cup C)) < 1$ or $f(S \cap (B \cup C)) < 1$.
\end{MyLemma}
So it is sufficient to independently separate the biconnected
components.
These three reductions can be combined recursively to yield separation
sub-problems that are significantly smaller than the original support
hypergraph.
For typical geometric instances with $n = |V|$, one usually obtains
reduced separation problems with at most $0.35 \, n$ vertices,
often {\em much} smaller.
Since the subsets $S \subset V$ identified as violations are subsets
of these individual reduced separation sub-problems, it is not
surprising that the violated subtours $S$ discovered tend to have
cardinality at most $0.35 n$, and usually much smaller.

Once the flow-based seperator has discovered a violated subtour
$S$, GeoSteiner performs ``constraint strengthening'' by recursively
applying the reductions of
Lemmas~\ref{lem:reduction}, \ref{lem:reduce-cc} and
\ref{lem:reduce-bcc} to the sub-hypergraph induced by $S$ to further
reduce the cardinality of the violation, possibly splitting $S$ (via
connected and biconnected components) into several much smaller (and
vastly stronger) violated subtours.

If the original violation $S$ has cardinality exceeding approximately
$0.35 \, n$, then the EPR and CD indicators
from~\cite{WarmeIndicators1}
indicate that these cardinality reductions actually {\em weaken} the
constraint --- at least until the ``weakest subtour'' threshold
has been passed.
When starting with a violation $S$ above this cardinality threshold,
the indicators suggest that the cardinality of $S$ must be
{\em increased} in order to strengthen the resulting subtour
inequality.
Previous versions of GeoSteiner strengthen subtour violations only via
reduction, not via augmentation.

This is why GeoSteiner (versions 5.3 and earlier) is virtually
incapable of identifying violated subtours $S$ having large
cardinality.
In order for these versions of GeoSteiner to generate a violated
subtour of large cardinality, the following events must happen
simultaneously:
\begin{itemize}
 \item The initial reductions must be relatively ineffective, yielding
   a ``reduced'' separation sub-problem containing most of $V$.
 \item The flow-separator must identify a violation $S$ of large
   cardinality within this sub-problem.
 \item The resulting $S$ must mostly survive all reductions
   subsequently applied to ``strengthen'' subtours emerging from the
   flow-separator.
\end{itemize}
This confluence of events seems to be impossible for the geometric
instances that typically arise within GeoSteiner, but could
potentially happen within certain more general classes of weighted
hypergraphs.

\section{Improvements to GeoSteiner}
\label{sec:geosteiner-mods}

We now present this paper's algorithmic improvements over GeoSteiner
version 5.3.
(The results reported
in~\cite{JuhlWarmeWinterZachariasen2018}
were obtained using GeoSteiner 4.0.
Version 5.3 contains little that would affect computational results
over those produced by version 4.0.)

The most important improvement is strengthening via augmentation.
Choosing branch variables via pseudo costs is the next most
significant improvement.
All of the remaining improvements do not broadly affect solution
times.
Some of these address specific algorithmic misbehaviors that manifest
on small subsets of the instances.
We present computational results only for the old vs new code.
We do not provide results independently assessing
the effectiveness of these changes individually, or in smaller
combinations.

\subsection{Constraint Strengthening via Augmentation}
\label{sec:strengthening-via-augmentation}

We now present the core algorithmic change of this paper, yielding one
of the most significant performance improvements in GeoSteiner
since~\cite{Warme98}.
This new algorithm provides an additional method for strengthening
subtour inequalities.

Let $S \subset V$ be a violated subtour.
GeoSteiner's existing method for strengthening $S$ applies the various
``reductions'' to the support hypergraph induced by $S$, resulting in
violations $S'$ such that $|S'| \le |S|$.
The new strengthening algorithm explicitly targets subtours of large
cardinality by ``augmenting'' $S$ with additional vertices using
methods that work by conceptually running the existing reductions ``in
reverse.''
We first prove some useful lemmas.
\begin{MyLemma}
Let $S,T \subseteq U \subseteq V$ such that
$S \cup T \,=\, U$ and $S \cap T \,=\, \emptyset$.
Then
\begin{equation}
\label{eq:lemma-U=S-cut-T-equation}
	x(U) ~=~ x(S) \,+\, x(S:T) \,+\, x(T).
\end{equation}
\end{MyLemma}
\begin{MyProof}
We prove the equivalence separately for each term $x_e$.
Let $e \in E$, $s \,=\, |e \cap S|$ and $t \,=\, |e \cap T|$.
The $x_e$ terms of~(\ref{eq:lemma-U=S-cut-T-equation})
can then be written
\begin{equation}
\label{eq:lemma-U=S-cut-T-proof-1}
	\max(s + t \,-\, 1,\, 0)
	~=~
	\max(s \,-\, 1,\, 0)
	\,+\,
	\lbrack (s \ge 1) \wedge (t \ge 1) \rbrack
	\,+\,
	\max(t \,-\, 1,\, 0)
\end{equation}
(where we use Knuth $\lbrack boolean \rbrack$ notation to denote the
value 1 when $boolean$ is true, and 0 when $boolean$ is false).
We must have either $s = 0$, $s = 1$ or $s \ge 2$.
Likewise, we must have either $t = 0$, $t = 1$ or $t \ge 2$.
It is straightforward case analysis to verify
that~(\ref{eq:lemma-U=S-cut-T-proof-1}) holds for all 9 of these
cases.
This proves equivalence for the $x_e$ term.
Since we did not specify any particular $e \in E$, the equivalence
holds for every $e \in E$.
\hfill\qed
\end{MyProof}
\begin{MyCorollary}
\label{cor:x(V)=x(S)+x(S:V-S)+x(V-S)}
For all $S \subseteq V$,
\begin{displaymath}
	x(V) ~=~ x(S) \,+\, x(S:V-S) \,+\, x(V-S).
\end{displaymath}
\end{MyCorollary}
We note the following alternate forms for subtour inequalities.
\begin{MyLemma}
\label{lem:alternate-subtour-forms}
Let $H \,=\, (V,E)$ be a hypergraph,
$S \subset V$,
$T \,=\, V \,-\, S$ and
$x \in \R^{|E|}$ such that
$|S| \ge 2$ and
\begin{displaymath}
	x(V) ~=~ |V| \,-\, 1.
\end{displaymath}
Then
\begin{eqnarray}
\nonumber
&&	x(S) \,\le\, |S| \,-\, 1 \\
\label{eq:alternate-subtour-form}
\Longleftrightarrow
&&	x(T) \,+\, x(T:V-T) \,\ge\, |T| \\
\nonumber
\Longleftrightarrow
&&	x(V-S) \,+\, x(S:V-S) \,\ge\, |V-S|.
\end{eqnarray}
\end{MyLemma}
\begin{MyProof}
\begin{align*}
&	(x(V) \,=\, |V| \,-\, 1) \wedge
	(x(S) \,\le\, |S| \,-\, 1) \\
\Longleftrightarrow
&	(x(S) \,+\, x(S:T) \,+\, x(T) \,=\, |V| \,-\, 1) \wedge
	(x(S) \,\le\, |S| \,-\, 1) \\
\Longleftrightarrow
&	x(S) \,+\, x(S:T) \,+\, x(T) \,-\, x(S)
	\,\ge\,
	(|V| \,-\, 1) \,-\, (|S| \,-\, 1) \\
\Longleftrightarrow
&	x(S:T) \,+\, x(T)
	\,\ge\,
	|V| \,-\, |S| \\
\Longleftrightarrow
&	x(T) \,+\, x(T:V-T) \,\ge\, |T| \\
\Longleftrightarrow
&	x(V-S) \,+\, x(S:V-S) \,\ge\, |V-S|.
\tag*{\qed}
\end{align*}
\end{MyProof}

\noindent
We call~(\ref{eq:alternate-subtour-form}) the
``anti-subtour inequality for $T$.''
It places a lower bound on the total edge weight within $T$ and
crossing out of $T$.
In contrast, the subtour inequality for $S$ places an upper bound on
the total edge weight within $S$.
Note that computations are free to use either the subtour inequality
for $S$, or the equivalent anti-subtour inequality for $T=V-S$,
whichever is more sparse.
(GeoSteiner has done this since version 4.0.)

\subsubsection{Single-Vertex Augmentation}
\label{sec:single-vertex-augmentation}
The first augmentation method reverses the
``iterative vertex removal''
reduction of Lemma~\ref{lem:reduction}.
\begin{theorem}[Single-vertex augmentation]
\label{th:1-vert-aug}
Let $S \subset V$,
$t \in V-S$ and
$S' \,=\, S \,\cup\, \lbrace t \rbrace$.
If $x(S) \,>\, |S| \,-\, 1$ and
$x(S:\lbrace t \rbrace) \,\ge\, 1$,
then $x(S') \,>\, |S'| \,-\, 1$.
\end{theorem}
\begin{MyProof}
\begin{align*}
x(S')
&=	x(S \,\cup\, \lbrace t \rbrace) \\
&=	x(S) \,+\, x(S:\lbrace t \rbrace) \\
&\ge	x(S) \,+\, 1 \\
&>	(|S| \,-\, 1) \,+\, 1 \\
&=	|S| \\
&=	|S'| \,-\, 1
\tag*{\qed}
\end{align*}
\end{MyProof}
This provides the conditions under which a violated subtour $S$ can be
augmented with vertex $t$ to obtain a violated inequality $S'$ of
larger cardinality.
This augmentation step is iterated greedily until no further valid
augmentations exist.

This can be implemented in
$O(|V| \, \log(|V|) \,+\, \sum_{e \in E} |e|)$
time using a heap containing the set $F$
of ``frontier'' vertices:
\begin{displaymath}
	F \,=\, \lbrace v \in V \,-\, S \,:\,
x(S:\lbrace v \rbrace) > 0 \rbrace
\end{displaymath}
Each vertex $v$ in the heap has corresponding key value
$x(S:\lbrace v \rbrace)$.
The heap is ordered to place vertex $v$ having the largest key at the
top of the heap.
The algorithm always terminates when the heap is empty.
Otherwise,
let $v$ (having key $z$) be the top element of the heap, and
$y > 0$ be the amount by which the current subtour $S$ is violated by
$x$.
If $y + z < 1$ then stop (augmenting by $v$ would result in an $S$ no
longer violated by $x$).
Otherwise perform the augmention steps:
(1) remove top element of heap;
(2) $S = S \cup \lbrace v \rbrace$;
(3) $y = y + z - 1$; and
(4) Each edge $e$ incident to $v$ with $x_e > 0$ having no other
vertices in $S$ now crosses cut $(S:V-S)$.
Update frontier vertices and keys:
For each vertex $u \in e$, $u \ne v$, (a) if $u$ is not in the heap,
insert it with key of zero; and (b) increase $u$'s key by $x_e$.

This augmentation step maintains the invariant $y > 0$ so that the
violation does not disappear.
Note, however, that the algorithm can augment with $v \in V-S$
such that
$x(S:\lbrace v \rbrace) \,<\, 1$, which {\em reduces} the amount of
violation.
If we assume that the algorithm stops instead of employing any such
``violation reducing'' augmentation, we can argue that the algorithm
always yields the same result $S$, regardless of the order in which
the augmentations are executed.
First note that heap keys only increase in value.
For all $(v,z)$ in the heap once $z \ge 1$, this condition never
becomes false and the augmentation by $v$ will eventually happen.
There is no way to ``miss'' the opportunity for any augmentation,
regardless of the order in which the augmentations are performed.
Thus the purely greedy algorithm always yields the same resulting $S$
(at least until any ``violation reducing'' augmentations are made).
For each valid sequence of ``violation reducing'' augmentations, the
resulting $S$ is the same regardless of the order in which the other
augmentations are performed.
Our implementation is the strictly greedy algorithm presented here ---
it makes no attempt to optimize the sequence of such
``violation reducing'' augmentations.

\subsubsection{Augmentation via Complementary Connected Components}
\label{sec:complementary-connected-components-augmentation}

The next augmentation method is analogous to reversing the ``connected
components'' reduction of Lemma~\ref{lem:reduce-cc}.

\begin{theorem}[Augmentation via connected components]
\label{th:augment-via-cc}
Let $H = (V,E)$ be a hypergraph,
$S \subset V$, and
$x \in \R^{|E|}$ such that
$|S| \,\ge\, 2$,
\begin{equation}
\label{eq:cc-theorem-total-degree-equation}
	x(V) \,=\, |V| \,-\, 1,
\end{equation}
and
\begin{displaymath}
	x(S) \,>\, |S| \,-\, 1.
\end{displaymath}
Let $\bar{K} \,=\, (\bar{V},\bar{E})$ be the support subhypergraph of
$H$ induced by $V-S$.
Let $(V_1,E_1),\, (V_2,E_2),\, \ldots,\, (V_k,E_k)$ be the connected
components of $\bar{K}$.
Let $S_i \,=\, V \,-\, V_i$ for $1 \le i \le k$.
Then there is at least one $i$, $1 \le i \le k$ for which
$x(S_i) \,>\, |S_i| \,-\, 1.$
\end{theorem}
\begin{figure}[ht]
 \centering
 \includegraphics[width=4.5in]{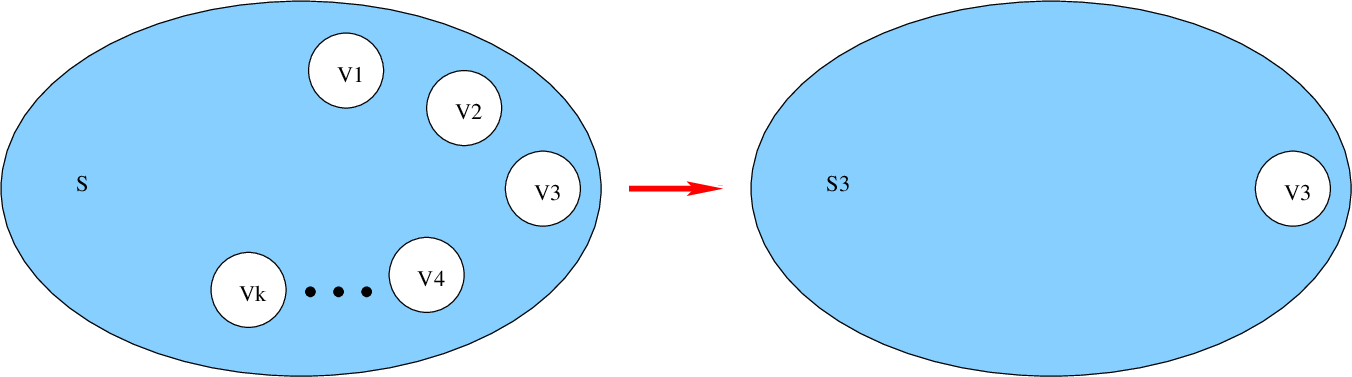}
 \caption{Illustration of Theorem~\ref{th:augment-via-cc}
	(shaded area is $S/S3$, outer ellipse is $V$)}
 \label{fig:cc1}
\end{figure}
\begin{MyProof}
The situation is as shown in the left side of Figure~\ref{fig:cc1}.
Assume that $x(S_i) \,\le\, |S_i| \,-\, 1$ for all $1 \le i \le k$.
For each $i$ such that $1 \le i \le k$, we have (by the alternate form
of Lemma~\ref{lem:alternate-subtour-forms}):
\begin{eqnarray}
\nonumber
&&	x(S_i) \,\le\, |S_i| \,-\, 1 \\
\label{eq:cc-holes}
\Longleftrightarrow
&&	x(V_i) \,+\, x(V_i:S_i) \,\ge\, |V_i|.
\end{eqnarray}
Because the $(V_i,E_i)$ are connected components of $\bar{K}$, there
is no edge $e \in \bar{E}$ crossing the cut $(V_i:V_j)$ for any
$1 \le i < j \le k$.
All edges crossing the cut $(V_i:S_i)$ must therefore cross only cut
$(V_i:S)$.
Subtracting the~(\ref{eq:cc-holes})
for all $1 \le i \le k$
from
Equation~(\ref{eq:cc-theorem-total-degree-equation})
yields:
\begin{eqnarray*}
&&	x(V) \,-\, \left ( \sum_{i=1}^k x(V_i) \,+\,
			   \sum_{i=1}^k x(V_i:S_i)
		   \right )
		\,\le\, |V| \,-\, 1 \,-\, \sum_{i=1}^k |V_i| \\
\Longrightarrow
&&	x(V) \,-\, \left ( \sum_{i=1}^k x(V_i) \,+\,
			   \sum_{i=1}^k x(V_i:S)
		   \right )
		\,\le\, |V| \,-\, 1 \,-\, \sum_{i=1}^k |V_i| \\
\Longrightarrow
&&	x(V) \,-\, (x(V-S) \,+\, x(V-S:S))
		\,\le\, |V| \,-\, 1 \,-\, |V-S| \\
\Longrightarrow
&&	x(V) \,-\, x(V-S) \,-\, x(V-S:S)
		\,\le\, |S| \,-\, 1 \\
\Longrightarrow
&&	\left ( x(S) + x(S:V-S) + x(V-S) \right )
		- x(V-S) - x(V-S:S))
		\le |S| - 1 \\
\Longrightarrow
&&	x(S) \,\le\, |S| \,-\, 1,
\end{eqnarray*}
which contradicts the given condition $x(S) \,>\, |S| \,-\, 1$.
The assumption that $x(S_i) \le |S_i| \,-\, 1$ for all $1 \le i \le k$
is therefore false, and there must be at least one $i$, such that
$1 \le i \le k$ for which $x(S_i) \,>\, |S_i| \,-\, 1$.
\hfill{} \qed
\end{MyProof}

Theorem~\ref{th:augment-via-cc} permits further step-increase in
subtour cardinality (from $|S|$ to $|S_i|$)
beyond that provided by Theorem~\ref{th:1-vert-aug}, while still
maintaining a violated subtour, as shown in the right side of
Figure~\ref{fig:cc1}.

Let $H$ and $x$ be as in Theorem~\ref{th:augment-via-cc}, let $I$ be
the support hypergraph of $H$, and let
$(V_1,E_1)$, $(V_2,E_2)$, $\ldots$, $(V_k,E_k)$
be the connected components of $I$.
As noted in Lemma~\ref{lem:reduce-cc},
for each $1 \le i \le k$ it suffices to
independently apply the separation oracle within $(V_i,E_i)$.
Let $S_i \subseteq V_i$ be such a violated subtour, and consider
strengthening of $S_i$ via Theorem~\ref{th:augment-via-cc}.
A direct application of Theorem~\ref{th:augment-via-cc} considers the
sub-hypergraph $J$ of $I$ induced by $V-S_i$, and the corresponding set
$K$ of all connected components of $J$.
For all $1 \le j \le k$, $j \ne i$ one finds that $(V_j,E_j)$ is a
member of $K$.  Thus, each of these $V_j$ get reconsidered as
potentially violated anti-subtour inequalities every time a new
violated subtour in $V_i$ is discovered and strengthened via
Theorem~\ref{th:augment-via-cc}.
This is a wasteful repetition of computational effort.
It suffices instead to consider the sub-hypergraph $L$ of ($V_i,E_i)$
induced by $V_i-S_i$, and the corresponding set $M$ of all connected
components of $L$.
Although one loses the guarantee of Theorem~\ref{th:augment-via-cc}
(that at least one of these connected components represents a violated
anti-subtour inequality), this guarantee is recovered over the union of
independently separating and strengthening each $(V_i,E_i)$, $1 \le i
\le k$.
This eliminates all such duplication of effort during strengthening
via Theorem~\ref{th:augment-via-cc}.

\subsubsection{Augmentation via Complementary Biconnected Components}
\label{sec:complementary-biconnected-components-augmentation}

The next method is analogous to reversing the ``biconnected
components'' reduction of Lemma~\ref{lem:reduce-bcc}.

\noindent
{\bf Definition:} A hypergraph is {\em biconnected} if it is
connected, and remains so after the removal of any single vertex.

\noindent
{\bf Definition:} An {\em articulation vertex} $t$ is a vertex of a
connected graph (hypergraph) whose removal splits the graph
(hypergraph) into two or more connected components.
In a connected graph (hypergraph), it is the articulation vertices
that form the boundaries between its biconnected components.

\begin{theorem}[Augmentation via biconnected components]
\label{th:augment-via-bcc}
Let $H = (V,E)$ be a hypergraph,
$S \subset V$ and
$x \in \R^{|E|}$ such that
$|S| \ge 2$,
\begin{displaymath}
	x(V) \,=\, |V| \,-\, 1,
\end{displaymath}
and
\begin{displaymath}
	x(S) \,>\, |S| \,-\, 1.
\end{displaymath}
Suppose
\begin{displaymath}
  x(S:\lbrace v \rbrace) \,<\, 1
  ~~~~~~~~~\hbox{for all $v \in V-S$}.
\end{displaymath}
Let $(V_i,E_i)$ and $S_i$ be from Theorem~\ref{th:augment-via-cc} such
that $x(S_i) \,>\, |S_i| \,-\, 1$.
Let
$(\overline{U}_1,\overline{E}_1)$,
$(\overline{U}_2,\overline{E}_2)$,
$\ldots$,
$(\overline{U}_m,\overline{E}_m)$
be the biconnected components of $(V_i,E_i)$.
Let $\overline{S}_j \,=\, V \,-\, \overline{U}_j$
for $1 \le j \le m$.
Then there is at least one $j$ for which
$x(\overline{S}_j) \,>\, |\overline{S}_j| \,-\, 1$.
\end{theorem}
\begin{MyProof}
If $(V_i,E_i)$ has no articulation vertices, then it is biconnected
and $(\overline{U}_1,\overline{E}_1)$ $= (V_i,E_i)$ is the only
biconnected component ($m \,=\, 1$).
We can choose to let $j \,=\, 1$.

Otherwise, $(V_i,E_i)$ has at least one articulation vertex.
Let $t$ be any articulation vertex of $(V_i,E_i)$.
Let $r \ge 2$ be the number of connected components that $(V_i,E_i)$
splits into upon removal of $t$, and let
$P_k$ be the vertex set of the $k$-th connected component,
for $1 \le k \le r$.
The situation is as shown in Figure~\ref{fig:bcc-proof}.
\begin{figure}[ht]
 \centering
 \includegraphics[width=4.5in,clip=]{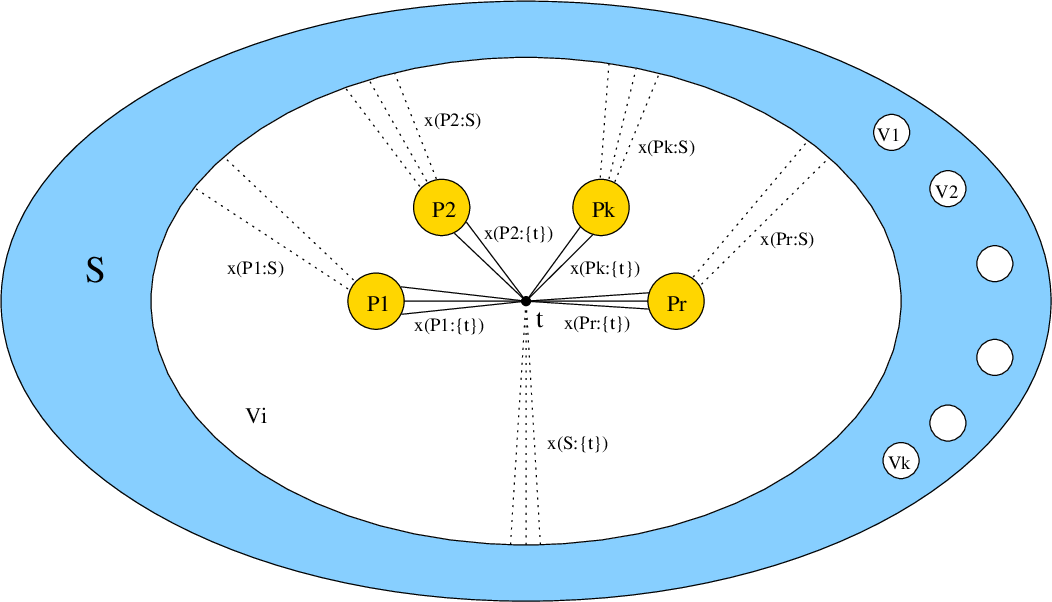}
 \caption{Proof of Theorem~\ref{th:augment-via-bcc}}
 \label{fig:bcc-proof}
\end{figure}

Let $Q_k \,=\, P_k \cup \lbrace t \rbrace$ for $1 \le k \le r$.
Assume that all of the (alternate form) subtour inequalities $Q_k$ are
satisfied, e.g.
\begin{equation}
\label{eq:Qk-subtours}
	x(Q_k) \,+\, x(Q_k:V-Q_k) \,\ge\,|Q_k|
	~~~~~~\hbox{for all $1 \le k \le r$}.
\end{equation}
We obtain a contradiction as follows.
Note that no edge crosses from $P_i$ to $P_j$ for any
$1 \le i < j \le r$.
Any edge crossing cut
$(P_k:V-P_k)$ can only cross
$(P_k:\lbrace t \rbrace)$,
$(P_k:S)$, or both,
as shown in Figure~\ref{fig:bcc-proof}.
We are given
\begin{eqnarray}
\nonumber
&&	x(S_i) \,>\, |S_i| \,-\, 1 \\
\nonumber
\Longrightarrow
&&	x(V_i) \,+\, x(V_i:V-V_i) \,<\, |V_i|
	~~~~~~~~\hbox{(alternate form)} \\
\nonumber
\Longrightarrow
&&	x(V_i) \,+\, x(V_i:S) \,<\, |V_i| \\
\label{eq:theta-source}
\Longrightarrow
&&	x(V_i) \,+\, x(V_i:S) \,-\, |V_i| \,<\, 0
\end{eqnarray}
We note that
\begin{displaymath}
	|V_i| ~=~ 1 \,+\, \sum_{k-1}^r |P_k|
\end{displaymath}
because of the need to count vertex $t$.
We define $\theta \in \R$ to be the left-hand side
of~(\ref{eq:theta-source}), expressed in terms of $t$ and the
$P_k$,
$1 \le k \le r$ (refer to Figure~\ref{fig:bcc-proof}):
\begin{eqnarray*}
\theta &=&
	\left ( \sum_{k=1}^r
		x(P_k) \,+\, x(P_k:\lbrace t \rbrace) \,+\, x(P_k:S)
		\,-\, |P_k|
	\right )
	\,-\, 1 \,+\, x(S:\lbrace t \rbrace) \\
	&<& 0
\end{eqnarray*}
We now sum all of the satisfied (alternate form)
$Q_k$ subtours~(\ref{eq:Qk-subtours}):
{\allowdisplaybreaks
\begin{eqnarray*}
&&
	\sum_{k=1}^r	x(Q_k) \,+\, x(Q_k:V-Q_k)
	\,\ge\,
	\sum_{k=1}^r	|Q_k| \\
&\Longrightarrow&
	\sum_{k=1}^r	x(Q_k) \,+\, x(Q_k:S)
	\,\ge\,
	\sum_{k=1}^r	|Q_k| \\
&\Longrightarrow&
	\sum_{k=1}^r
	    \left \lbrack
			\vphantom{\sum_{k=1}}
		x(P_k)
		\,+\, x(P_k:\lbrace t \rbrace)
		\,+\, x(P_k:S)
		\,+\, x(S:\lbrace t \rbrace)
	    \right \rbrack
	\,\ge\,
	r \,+\, \sum_{k=1}^r |P_k| \\
&\Longrightarrow&
	\sum_{k=1}^r
	    \left \lbrack
			\vphantom{\sum_{k=1}}
		x(P_k)
		\,+\, x(P_k:\lbrace t \rbrace)
		\,+\, x(P_k:S)
		\,+\, x(S:\lbrace t \rbrace)
		\,-\, |P_k|
	    \right \rbrack
	\,\ge\, r \\
&\Longrightarrow&
	r \, x(S:\lbrace t \rbrace)
	\,+\,
	\left \lbrack
	\sum_{k=1}^r
		x(P_k)
		\,+\, x(P_k:\lbrace t \rbrace)
		\,+\, x(P_k:S)
		\,-\, |P_k|
	\right \rbrack
	\,\ge\, r \\
&\Longrightarrow&
r \,\, x(S:\lbrace t \rbrace) \,+\,
	\left \lbrack
		\theta \,+\, 1 \,-\, x(S:\lbrace t \rbrace)
	\right \rbrack
	\,\ge\, r \\
&\Longrightarrow&
(r \,-\, 1) \, x(S:\lbrace t \rbrace) \,+\, \theta \,+\, 1 \,\ge\, r \\
&\Longrightarrow&
\theta \,\ge\, r \,-\, 1 \,-\, (r \,-\, 1) \, x(S:\lbrace t \rbrace) \\
&\Longrightarrow&
\theta \,\ge\, (r \,-\, 1) \, (1 \,-\, x(S:\lbrace t \rbrace)) \\
\end{eqnarray*}
}
Because $\theta < 0$ and $r \ge 2$, this implies that
$x(S:\lbrace t \rbrace) > 1$.
We are given that $x(S:\lbrace v \rbrace) \,<\, 1$
for all $v \in V-S$.
This is contradicted for $v=t$.
There must therefore be at least one $Q_k$ ($1 \le k \le r$)
for which the (alternate form) subtour~(\ref{eq:Qk-subtours}) is
violated.
Let $\bar{S}_k \,=\, V \,-\, Q_k$.
This implies that subtour $\bar{S}_k$ is violated:
$x(\bar{S}_k) \,>\, |\bar{S}_k| \,-\, 1$.

We did not specify which articulation point $t$ to choose, so the
above argument can be repeated for every articulation vertex $t$ of
$(H_i,E_i)$.
These articulation points partition $(H_i,E_i)$ into biconnected
components
\begin{displaymath}
(\overline{U}_1,\overline{E}_1),\,
(\overline{U}_2,\overline{E}_2),\,
\ldots,\,
(\overline{U}_m,\overline{E}_m).
\end{displaymath}
At least one of the $(\overline{U}_j,\overline{E}_j)$ must therefore
be such that (alternate form) subtour $\overline{U}_j$ is violated:
$x(\overline{U}_j) \,+\, x(\overline{U}_j:V-\overline{U}_j)
\,<\, |\overline{U}_j|$.
Let $\overline{S}_j \,=\, V \,-\, \overline{U}_j$.
Then subtour $\overline{S}_j$ is a violated subtour:
$x(\overline{S}_j) \,>\, |\overline{S}_j| \,-\, 1$.
\hfill{} \qed
\end{MyProof}

Theorem~\ref{th:augment-via-bcc} permits additional growth of $S$ beyond that
of Theorem~\ref{th:augment-via-cc}, while still maintaining a violated
subtour.

Our improved constraint strengthening implementation within GeoSteiner
processes every violated subtour $S$ generated by the deterministic
flow separation algorithm described in~\cite{Warme98}.
It first applies the existing ``strengthening by reduction'' methods
to $S$.
It then applies the new ``strengthing by augmentation'' methods to $S$
as follows.
Apply the greedy single-vertex augmentation algorithm of
Section~\ref{sec:single-vertex-augmentation} to $S$, yielding an
augmented violation $S'$.
Then compute the connected components
$(V_1,E_1),\, \ldots,\, (V_k,E_k)$
of the support hypergraph induced by
$V-S'$.
For each violated anti-subtour $V_i$, $1 \le i \le k$:
(1) generate anti-subtour inequality $V_i$; and
(2) compute the biconnected components
$(\bar{U}_j, \bar{E}_j)$ of $V_i$ ---
generating each anti-subtour $\bar{U}_j$ that is violated.
Subtours $V-V_i$ and $V-\overline{U_i}$ are usually of quite large
cardinality --- which the EPR and CD predict to be very strong.
Because they are violated by the {\em current} LP solution $x$, they
are also very pertinent to the problem being solved
(unlike the complementary subtours of
Section~\ref{sec:complementary-subtours} below, which are almost never
violated by the current LP solution $x$ --- only ``incidentally'' by
future LP solutions $x$).

\subsection{Choosing Branch Variables via Pseudo Costs}
\label{sec:pseudo-costs}

The branch variable selection algorithm used (by default) in
GeoSteiner versions 5.3 and earlier is useful on problems that require
at most a few dozen branch-and-bound nodes, but becomes hopelessly
inefficient when much more branching is required.
Let $F$ be the list of fractional variables, which the algorithm
places into a heuristically sorted order.
It then performs ``strong branching''
(\cite{ApplegateBixbyChvatalCook2011} Section 14.3,
actually solving the LP for both branches)
on the members of $F$ until $k$ consecutive members tested
fail to find an improved branch, where
$k \,=\, 2 \, \lfloor \log_2(|F|) \rfloor$.
It chooses the variable whose minimum ``down branch'' versus ``up
branch'' objectives is maximized.
When solving the ``down branch'' and ``up branch'' LPs, the code also
invokes the primal upper bound heuristic on each resulting LP solution.
(This primal heuristic turns out to be significantly more costly than
solving the LP for each branch, especially on very large problems.)
The default algorithm stops immediately on any variable for which one
or more of the branches yields an infeasible LP or a cutoff.
Although this causes the variable to be fixed without branching, the
same variable selection process usually starts over from scratch,
testing the same variables in the same order as before.
This algorithm wastes significant amounts of CPU time, e.g., on random
rectilinear instances above 4000 points, becoming completely
intolerable above 8000 points.

We implemented a new branch variable selection algorithm using
pseudo-costs to choose branch
variables~\cite{BenichouGauthierGirodetHentgesBirbiereVincent1971},
\cite{LinderothSavelsbergh1999}.
Let $x$ be an LP solution, $z = c \cdot x$ be its optimal objective
value, and $e \in E$ such that $x_e$ is fractional ($0 < x_e < 1$).
We define
\begin{eqnarray*}
c^0_e &=& {{z^0_e \,-\, z} \over {x_e}},	\hbox{ and} \\
c^1_e &=& {{z^1_e \,-\, z} \over {1 \,-\, x_e}},
\end{eqnarray*}
where $z^0_e$ is the LP objective value with branch condition $x_e =
0$ imposed, and $z^1_e$ is the LP objective value with branch
condition $x_e = 1$ imposed.
We call $c^0_e$ the ``downward pseudo cost of $x_e$,''
and $c^1_e$ the ``upward pseudo cost of $x_e$.''
Using pseudo costs, a good branch variable $x_e$ can be selected using
\begin{displaymath}
	e = \mathrm{argmax}_{e \in \tilde{E}}
		(\min (c^0_e \, x_e,\, c^1_e \, (1 \,-\, x_e))),
\end{displaymath}
where
$\tilde{E} = \lbrace e \in E : \hbox{$x_e$ is fractional} \rbrace$.
If the pseudo costs are always accurate, this is equivalent to
``full'' strong branching at each node.

Calculating pseudo costs for a single fractional variable is mildly
expensive, requiring the solution of two LPs.
Calculating them for all fractional variables assures accuracy, but is
{\em very} expensive.

It has been observed that pseudo costs generally remain fairly stable
during the branch-and-bound process.
This property can be exploited to choose very good branch variables at
reasonably low amortized cost.

We tried a simple method having a single global array of pseudo costs
(indexed by $e \in E$, initially undefined), where pseudo costs were
calculated as needed and cached in this array.
At the root node, this means that pseudo costs are calculated for all
fractional variables.
This method performs very well on problems requiring a few thousand
nodes or less, but gets progressively worse with larger
branch-and-bound trees, because cached pseudo costs diverge from
their true values, which are node-specific.

In our improved pseudo cost method:
\begin{itemize}
 \item Pseudo costs are stored for each branch-and-bound node,
 \item Child nodes inherit pseudo costs from their parents, and
 \item Pseudo costs are occasionally re-calculated so that they remain
	more accurate during branch-and-bound.
\end{itemize}
Unfortunately, we don't know which pseudo costs change from parent to
child node (unless we re-calculate them all).
We instead use a heuristic to re-calculate ``some'' pseudo costs at
each node, seeking to balance pseudo cost accuracy (and branch
variable quality) versus computational cost.

Our heuristic is presented in Figure~\ref{fig:pseudo-cost-heuristic}.
The set $L \subseteq E$ of edges for which pseudo costs are calculated
consists of
\begin{itemize}
 \item $U$: edges having undefined pseudo costs;
 \item $O$: the $k$ oldest pseudo costs; and
 \item $S$: a ``stratified random sample'' of remaing fractional
   variables.
   The starting offset for this sample is not actually random, but
   calculated from a hash of all branch conditions imposed upon this node
   up to the root node.
\end{itemize}
The primal heuristic is invoked only on the children of the branch
variable finally selected.

\MyBeginFig
{
\ttfamily
\scriptsize
\AL{0}{update\_pseudo\_costs (node)}
\AL{0}{$\lbrace$}
\AL{1}{	F = $\lbrace e \in E : \hbox{$x_e$ is fractional} \rbrace$;}
\AL{1}{	U = $\lbrace e \in F : \hbox{pseudo costs for $x_e$ undefined.} \rbrace$;}
\AL{1}{	R = F - U;}
\AL{1}{	Sort R by pseudo cost age, newest to oldest;}
\AL{1}{	k = $\lfloor \sqrt{|\mathrm{R}|} / 4 \rfloor$;}
\AL{1}{	O = oldest k elements of R;}
\AL{1}{	R = R - O;}
\AL{1}{	stride = $\lfloor |\mathrm{R}| / k \rfloor$;}
\AL{1}{	offset = hash (node) mod stride;}
\AL{1}{	S = $\emptyset$;}
\AL{1}{ i = offset;}
\AL{1}{	While i < |R| Do}
\AL{2}{		S = S $\cup$ R$\lbrack$i$\rbrack$;}
\AL{2}{		i = i + stride;}
\AL{1}{	End}
\AL{1}{	L = U $\cup$ O $\cup$ S;}
\AL{1}{	For each $e \in$ L Do}
\AL{2}{		Re-calculate pseudo costs for $x_e$;}
\AL{1}{	End}
\AL{0}{$\rbrace$}%
}
\begin{center}
\captionof{figure}{Pseudo cost re-calculation heuristic}%
\label{fig:pseudo-cost-heuristic}
\end{center}
\MyEndFig

We do not claim this method of choosing branch variables anywhere
close to state-of-the-art.
It is a simple method that generally performs better than GeoSteiner's
current default branch variable selection method.
Because it evaluates pseudo costs for {\em every} fractional variable
at the root node, this new method is generally a bit slower than
GeoSteiner's default method on problems having small branch-and-bound
trees.
It does, however, perform substantially better on problems having
medium to large branch-and-bound trees.

\subsection{Improved Local Cuts}
\label{sec:local-cuts}

We made several improvements to GeoSteiner's local cut
implementation originally described
in~\cite{JuhlWarmeWinterZachariasen2018}.
As usual, the first step is to choose a projection
$\hat{x} = P.x$, where $x$ is a vector in the original solution space,
$\hat{x}$ is a vector in the projected space, and
$\dim(\hat{x}) \ll \dim(x)$.
The chosen projection now contracts long chains of $k \ge 2$ edges
having LP weight of 1 down to a single 2-vertex edge.
Because there are $2^k$ subsets of such a chain, each such subset
could appear in max-weight forest subproblems (testing validity of
the current projected cut), resulting in large numbers of generated
forests and local cut iterations.
Contracting such chains to a single edge eliminates this misbehavior.

Let $\hat{E}$ be the set of projected edges, and
$\hat{C} \subseteq \hat{E}$ be the contracted edges.
Let $\hat{e} \in \hat{C}$, and
$e_1,\, e_2,\, \ldots,\, e_k$ be its original chain of edges.
The $\hat{e}$ component of the projection is
\begin{displaymath}
	\hat{x}_{\hat{e}} ~=~
		x_{e_1} \,+\, x_{e_2} \,+\, \cdots \,+\, x_{e_k},
\end{displaymath}
for which only two cases are relevant:
$\hat{x}_{\hat{e}} = k$ (contracted edge is present), and
$\hat{x}_{\hat{e}} \le k-1$ (contracted edge not present).
For each contracted edge $\hat{e}_i \in \hat{C}$,
let $k_i \ge 2$ be the number of edges in the original chain.
The present implementation utilizes projected vectors for which
$\hat{x}_{\hat{e}_i} \in [0,k_i]$, which are no longer 0-1 vectors.

The Farkas Lemma works as-is on such vectors, producing candidate
hyperplanes that separate these vectors from the projected fractional
solution $x$.
The validity test, however, must change.
Let $a.x \le b$ be the current candidate projected inequality
(constructed via the Farkas Lemma).
The validity step seeks a projected integer feasible solution (IFS)
$\hat{x}$ that maximally violates $a.\hat{x} \le b$.
In projected IFS $\hat{x}$,
$\hat{x}_{\hat{e}_i}$ can be any integer in $[0,k_i]$.
The max-weight forest sub-problems are computed just as before
(and still yield 0-1 values for each edge, whether contracted or not),
but with slightly different edge weights $c$.
For all $\hat{e} \in \hat{E} \,-\, \hat{C}$, we use
$c_{\hat{e}} = a_{\hat{e}}$, as before.
For all $\hat{e}_i \in \hat{C}$, however, we use
\begin{displaymath}
c_{\hat{e}_i} ~=~
	\begin{cases}
		a_{\hat{e}_i}	&
		\text{if $a_{\hat{e}_i} \ge 0$;} \\
		k_i \, a_{\hat{e}_i} &
		\text{if $a_{\hat{e}_i} < 0$.}
	\end{cases}
\end{displaymath}
In the case where $a_{\hat{e}_i} \ge 0$, the max-weight forest
sub-problem objective has an additional constant term of
$(k_i \,-\, 1) \, a_{\hat{e}_i}$.
Let $\bar{x}$ be the 0-1 incidence vector of the maximum (or
sufficiently large) weight forest.
We construct the corresponding projected integer feasible solution
as follows.
For all $\hat{e} \in \hat{E} \,-\, \hat{C}$,
let $\hat{x}_{\hat{e}} \,=\, \bar{x}_{\hat{e}}$.
For all $\hat{e}_i \in \hat{C}$, let
\begin{displaymath}
\hat{x}_{\hat{e}_i} ~=~
	\begin{cases}
		\bar{x}_{\hat{e}_i} \,+\, k_i \,-\, 1	&
		\text{if $a_{\hat{e}_i} \ge 0$;} \\
		k_i \, \bar{x}_{\hat{e}_i}		&
		\text{if $a_{\hat{e}_i} \,<\, 0$.} \\
	\end{cases}
\end{displaymath}
This reduction greatly reduces local cut iterations when long chains
of integral edges are present in the projected portion of the support
hypergraph.

A second mode of misbehavior is that some of the max-weight forest
sub-problems require large amounts of CPU time to solve.
As explained in~\cite{JuhlWarmeWinterZachariasen2018},
the max-weight forest sub-problems reduce to MST in hypergraph
instances that are then solved by recursively calling GeoSteiner's
branch-and-cut.
This misbehavior was due to local cuts being enabled on these
recursive calls ---
essentially all of which were solved without any branching.
On these smaller projected instances, local cuts were turning
GeoSteiner's branch-and-cut algorithm into a pure cutting-plane
algorithm.
(Branching commences only when no further violated local cuts are
found.)
It is well-known that pure cutting plane methods can stall at states
of high dual degeneracy that persist through huge numbers of generated
cuts.
Because the objective does not improve during such stalls, slack rows
cannot be deleted, greatly expanding the LP and CPU time.
Disabling local cuts during these recursive calls completely cured
this misbehavior.

A third mode of local cut misbehavior happened when this same stalling
behavior occurs during the main MST in hypergraph instance (degenerate
local cuts --- many consecutive iterations during which only local
cuts are generated, with no improvement in the objective function).
This was corrected by: incrementing a (per node) counter after
separation runs yielding only local cuts; zeroing this counter when
the node objective improves; and skipping local cuts when this counter
reaches a fixed limit (10 by default).

In GeoSteiner 5.3, local cuts are limited (by default) to projected
components having at most 80 vertices and 256 edges.
In the present experimental version, these defaults were reduced to
64 vertices and 96 edges
(which actually reduces the effectiveness of local cuts because it
``gives up'' on local cuts sooner).
The computational results reported here use these differing limits.
With virtually no local cut misbehavior remaining, however, these
defaults should probably be restored to their previous values, making
local cuts more aggressive overall than before (with chain contraction
reducing the size of components, allowing more of them to satisfy
these limits, thereby yielding additional local cuts).

\subsection{Integrality Delta}
\label{sec:integrality-delta}

When solving rectilinear instances, the objective coefficients (FST
lengths) are effectively integers.
This permits a more aggressive policy for deciding to ``cut off'' a
branch-and-bound node.
Let $U$ be the current upper bound, and $z$ be the objective value of
a branch-and-bound node.
Under normal circumstances a cutoff is only valid when $z \ge U$.
When the objective coefficients are all integral, however, the node
can be cut off when $z > U - \Delta$, where $\Delta = 1$ is the
minimum positive difference in objective values across all integer
feasible solutions.
To account for numerical error, the present implementation (for
integral objectives only) is more conservative.
Let $Limit_1 = U - 1$,
$Limit_2 = Limit_1 + 1/32$ and
$Limit_3 = Limit_1 + 1/4$.
If $(double) Limit_2 > (double) Limit_1$,
then cut off when $z \ge Limit_2$.
Otherwise, if $(double) Limit_3 > (double) Limit_1$,
then cut off when $z \ge Limit_3$.
Otherwise, cut off when $z > Limit_3$.
There are a few rectilinear instances for which this significantly
reduces the branch-and-cut ``end game.''

\subsection{New Method for Solving over all Constraints in the Pool}
\label{sec:SOS-heuristic}

Given an LP solution $x$, GeoSteiner versions 5.3 and earlier
scan the pool for {\em all} inequalities violated by $x$ and add
{\em all} of these inequalities to the LP at once.
This scan is repeated (after re-solving the updated LP) until an LP
solution $x$ is obtained that satisfies all inequalities in the
constraint pool.
(The code also deletes slack constraints inside this loop if the
LP objective strictly improves.)
This algorithm has the undesirable side-effect of filling up the LP
tableaux with many very dense rows, including presumeably many rows
that are redundant.
This effect is particularly acute on problems that require many
optimize/separate iterations in order to satisfy all subtour
inequalities at the root node, especially large random rectilinear
instances.
To address this problem we devised the
``Sparse, Orthogonal, Stratified'' (SOS) heuristic illustrated
in Figure~\ref{fig:SOS-heuristic}.
\MyBeginFig
{
\ttfamily
\scriptsize
\AL{0}{SOS\_heuristic ($x$)}
\AL{0}{$\lbrace$}
\AL{1}{	$L ~=~$ list of all constraints in pool violated by $x$;}
\AL{1}{	$B ~=~$ list of all constraints in pool binding for at least 1
  active node;}
\AL{1}{	$L ~=~ L \,-\, B$;}
\AL{1}{	Sort $L$ from least number of non-zeros to most non-zeros;}
\AL{1}{	$C ~=~ \emptyset$;}
\AL{1}{	Repeat at most $k$ times while $L$ is non-empty:}
\AL{2}{		$D ~=~ \emptyset$;}
\AL{2}{		For each constraint $i$ in $L$ Do}
\AL{3}{			If $\mathrm{nonZeroColumns}($i$) \cap
  			    \mathrm{nonZeroColumns}($D$) ~=~ \emptyset$}
\AL{3}{			Then}
\AL{4}{				Move constraint $i$ from $L$ into set $D$;}
\AL{3}{			EndIf}
\AL{2}{		End}
\AL{2}{		$C ~=~ C \,\cup\, D$;}
\AL{1}{	End}
\AL{1}{	return ($C$);}
\AL{0}{$\rbrace$}%
}
\begin{center}
\caption{Algorithm SOS\_heuristic}%
\label{fig:SOS-heuristic}%
\end{center}
\MyEndFig

The sorting step causes the heuristic to favor sparse rows over dense
rows.
The ``If'' statement chooses rows that are orthogonal to each other
(no non-zero columns in common).
Repeating this scan $k$ times causes several ``strata'' of such
constraints to be chosen during each such iteration over the pool.
Thus the name Sparse, Orthogonal, Stratified.
($k=8$ by default.)

This heuristic is reasonably effective at reducing the size and
density of the LPs, especially on problems where many
optimize/separate iterations are required.
It is not totally effective, however.
On very large, long-running problems we find that the LP can still
be substantially reduced in size by starting over with the initial set
of constraints (together with the current node's branching
constraints) and re-running the SOS heuristic until all constraints in
the pool are satisfied.
Manual testing indicates that this ``restart'' operation is
computationally quite costly.
Exploiting restarts algorithmically would require some policy to
decide when to invoke the ``restart'' operation, and do so only when
the benefit usually outweighs the cost.
(No restarts were used in any of the computational results reported
herein.)

Although reasonably effective at keeping the LP smaller on very large
problems, smaller problems generally experience a slight increase in
solution time because of having to execute more ``SOS'' passes over
the constraint pool before achieving an $x$ that satisfies all
constraints in the pool.
We accept this small penalty in exchange for smaller LPs and better
performance on large problems.

\subsection{Zero Weight Cutsets}

GeoSteiner uses a special subtour separation algorithm on support
hypergraphs having $k \ge 2$ connected components.
When $k \le Threshold$ (5 by default), geosteiner generates a subtour
(if violated) for each partition $(S,V-S)$ of the $k$ connected
components into two non-empty sets.
Otherwise, for each connected component $S \subset V$ it generates two
subtours (if violated), one for $S$ and one for $V-S$.
The first of these methods was eliminated in the experimental version
presented here.
Only the latter type of ``simple'' subtours are now used.

\subsection{Complementary Subtours}
\label{sec:complementary-subtours}

To exploit the greater ``predicted'' strength of large cardinality
subtours, we tried the following simple idea. 
Whenever some subtour $S$ gets added to the constraint pool, also add
the complementary subtour $V-S$ to the constraint pool.
When $|S|$ is small, $|V-S|$ will be large.
Although subtour $V-S$ might not be violated by the present LP
solution $x$, it might become violated by future LP solutions $x'$, at
which time the next scan of the constraint pool will add it to the LP.

Such ``complementary'' subtours do in fact get pulled into the problem
occasionally (5--25\% of all SOS pool scans contain one or more
such complementary subtours).
These provide no statistically significant impact on solution
times, however.
Although these are subtours of large cardinality, they do not seem to
be the ``correct'' large-cardinality subtours for the instance.

\subsection{New Contraint Representation}
\label{sec:constraint-representation}

Violated subtours $S$ were previously represented as a dense bitmask
of vertices whose $v$-th bit is one, for each $v \in S$.
At $n=1,000,000$, each such bitmask occupies 125,000 bytes of memory.
We therefore modified the code so that
(1) subtours are now stored as variable-length lists of vertices in
(or not in) $S$; and
(2) a hash table is now used to detect duplicate subtours the moment
they are created rather than storing multiple duplicate copies that
must be deduplicated at a later time.

\subsection{New Constraint Pool Size Management}
\label{sec:constraint-pool}

In GeoSteiner versions 5.3 and earlier, the constraint pool can grow
unacceptably large.
We institute a (user-adjustable) hard limit on the total number of
non-zero coefficients for constraints in the constraint pool.
This limit is now exceeded only when the set of binding constraints
among all active branch-and-bound nodes must exceed this limit.

The methods used here are quite standard.
(1) Compute the number of rows and non-zeros needed by the set of
constraints being added to the pool.
(2) Compute minimum number of rows and non-zeros by which the current
constraint pool must be reduced in order to make room for the new
constraints.
(3) Sort the existing constraints in the pool from largest number of
non-zeros down to least number of non-zeros (ignoring constraints that
are binding for any active branch-and-bound node).
(4) Greedily remove dense constraints until the size limits have been
reached (or all discardable constraints have been removed).
(5) Add the new constraints to the pool.

These changes have almost no measurable impact upon the overall run
time of GeoSteiner's branch-and-cut algorithm, but provide a much
more reliable upper bound on the size of the constraint pool.

\section{Computational Results}
\label{sec:computational-results}

All computations were performed on a computer with:
\begin{itemize}
 \item CPU: AMD Ryzen Threadripper 2950X (16 cores, base clock 3.5GHz
	max boost clock 4.4GHz, 40MB combined cache)
 \item Memory: 128GB DDR4 3200MHz ECC
 \item Motherboard: Asus ROG Strix X399-E Gaming
\end{itemize}
The operating system is the 64-bit version of Rocky Linux 8.6.
The source code is compiled using GCC 8.5.0-10 with the
optimization flag -O3.
The LP solver is CPLEX version 12.5.1.
Because GeoSteiner implements its own native branch-and-cut framework,
it uses CPLEX only as an LP solver, and generally uses default
settings for all CPLEX LP solver algorithmic parameters.

The FST generation and pruning code has not changed, so there is no
need for separate old / new runs of these phases.
Since only the FST concatenation code has changed, it is the change in
FST concatenation time that is of most interest in these results.
All CPU times are reported in seconds.
The test suite consists of the following instances:
\begin{itemize}
  \item OR-Library (31 instances x 4 metrics)
  \item Random instances (10 sizes x 15 instances/size x 4 metrics)
  \item Large random instances (sizes 15k, 20k, 25k, 30k, 40k, 50k and
    100k points, 15 instances/size --- Euclidean for all
    of these, and rectilinear for instances up to 30k points).
  \item Very large random instances --- 15 instances each at 250k and
    500k points, plus one instance at 1,000,000 points (Euclidean
    metric only).
  \item TSPLIB instances (46 instances x 4 metrics)
\end{itemize}
This totals 1104 instances, all ideally solved with both the old and
new code.
In reality, the old code cannot solve the very large random
Euclidean instances, or rectilinear instances above 8000
terminals, so we did not attempt these.
We use a time limit of 24 CPU hours for the FST concatenation phase of
each instance (except for the 1,000,000 point Euclidean instance).
(Note that much more lenient time limits were used for the FST
generation and FST pruning phases.
Some of these instances required several days or weeks of CPU time to
obtain the FSTs.
This was done to maximize the number of test instances, and because
the algorithmic improvements being studied do not reside within the
FST generation or FST pruning code.)

We performed FST generation and FST pruning on all instances before
moving on to FST concatenation (which was performed twice on most
instances: once with the old code and once with the new code).
In most cases, 8 instances were executed simultaneously (using 8 out
of 16 available CPU cores), although this was not tightly controlled.
(In certain cases, e.g., the end of a sequence of queued jobs, the
number of active cores dropped below 8.)
Because of cache interference between cores and limited memory
bandwidth, one should expect the reported execution times to be within
some reasonable factor (e.g., 15--30 percent?) of what might be
expected if these computations were repeated (in some not carefully
controlled order).
There were a few FST concatenation instances that exceeded the 128 Gb
physical memory limit of the machine.
These occasions forced the Linux kernel to perform swapping,
but the access patterns were always ``reasonable'' (not thrashing)
and the impact upon other processes was unlikely to be acute.
The main exceptions are the very large random Euclidean FST
concatenation runs:
(1) instance \#1 at 250k, 500k and 1,000,000 points were run all alone
(on an otherwise idle CPU) in order to get accurate CPU times;
(2) instances \#2--15 at 250k and 500k points were run 5 instances at
a time.

Of the original 1104 instances, no FSTs are available for the
12 instances shown in Table~\ref{tab:no-fsts-instances}:
\begin{table}
\setlength{\tabcolsep}{8.0pt}
\begin{center}
  \captionof{table}{Instances for which no FSTs are available}
  \begin{tabular}{|lcccc|}\hline
	Instance	& Euclidean & Rectilinear & Hexagonal & Octilinear \\
	\hline
	p654		&	&	& X	& 	\\
	d2103		& X	&	&	&	\\
	u2319		& X	&	& X	&	\\
	fl3795		& X	&	&	&	\\
	pla7397		& X	&	& X	&	\\
	pla33810	& X	&	& X	&	\\
	pla85900	& X	&	& X	& X	\\
	\hline
  \end{tabular}
  \label{tab:no-fsts-instances}
\end{center}
\end{table}
This leaves 1092 instances for which FST concatenation was attempted.

\subsection{Overall Impact of Improvements}
\label{sec:impact-of-improvements}

As shown in Table~\ref{tab:subtour-timeout-instances},
there are just 4 of these 1092 instances for which the new code
fails to obtain the optimal subtour relaxation within the 24 hour CPU
time limit:
\begin{table}[!ht]
\setlength{\tabcolsep}{8.0pt}
\begin{center}
  \captionof{table}{Instances for which new code does not satisfy all
    subtours within 24 CPU hours}
  \begin{tabular}{|llr|}\hline
	Instance~~ & Metric	& ~~Iterations \\
	\hline
	u2319	 & Octilinear	& $> 67$ \\
	d15112	 & Rectilinear	& $> 172$ \\
	d18512	 & Rectilinear	& $> 168$ \\
	pla33810 & Octilinear	& $> 120$ \\
	\hline
  \end{tabular}
  \label{tab:subtour-timeout-instances}
\end{center}
\end{table}
In all 4 cases, the problem is not an excessive number of iterations,
but rather LPs that take a long time to solve.
(The author believes that the new code would likely satisfy all
subtours for these instances in a reasonable number of iterations if
sufficient additional CPU time were provided for solving these
difficult LPs.)
Table~\ref{tab:most-subtour-iteration-instances}
presents the instances requiring the largest number of iterations.
\begin{table}[!ht]
\setlength{\tabcolsep}{4.0pt}
\begin{center}
  \captionof{table}{Instances requiring the most iterations to
    converge all subtours}
  \begin{tabular}{|llr| |llr|}\hline
	Instance & Metric & Iterations~ & ~Instance & Metric & Iterations \\
	\hline
	brd14051 & Octilinear	&	156~ & ~usa13509 & Octilinear  & 207 \\
	d18512   & Octilinear	&	157~ & ~usa13509 & Rectilinear & 208 \\
	1,000,000 point & Euclidean &	163~ & ~usa13509 & Hexagonal   & 211 \\
	brd14051 & Octilinear	&	164~ & ~d18512 & Hexagonal     & 217 \\
	brd14051 & Rectilinear	&	165~ & ~d15112 & Octilinear    & 294 \\
	d15112 & Euclidean	&	187~ & ~d15112 & Hexagonal     & 334 \\
	\hline
  \end{tabular}
  \label{tab:most-subtour-iteration-instances}
\end{center}
\end{table}
All other instances require fewer than 150 iterations to satisfy all
subtour inequalities.

The enhancement described in
Section~\ref{sec:strengthening-via-augmentation}
is so computationally successful that the new code
(with just the four exceptions listed above)
now rapidly satisfies the root node subtour contstraints for
{\em every}
instance in our test set.
At most 334 optimize/separate iterations are required for each of
these problems.
Because so few iterations are required, there is little opportunity
for useless or redundant constraints to accumulate within the LP.
Since the LP contains generally stronger constraints (which are more
sparse),
the LP A matrix remains much smaller,
the LPs solve much more quickly,
and the solution vector $x$ is much cleaner numerically (i.e.,
rational numbers with small denominators).
Cleaner $x$ cause the reductions to remain much more effective,
decreasing the size of instances presented to the expensive flow
separator.
These effects greatly decrease the overall cost of the separation
algorithm.
It is significant that simply adding some new, stronger constraints to
the mix of violated subtours generated produces so many virtuous
feedback cycles within the branch-and-cut algorithm as a whole.

The improved stregthening has increased the range of practically
solvable instances by a factor of 10, which has exposed new
bottlenecks within the existing GeoSteiner code.
(The changes described in Section~\ref{sec:constraint-representation}
address the most significant of these new bottlenecks.)

The net result is that failure to converge the subtour relaxation has
essentially vanished from our large suite of test problems.
(Convergence of local cuts could potentially be a lingering issue, but
the fix for degenerate local cuts described in
Section~\ref{sec:local-cuts} eliminated such misbehavior from the
entire suite of test instances.)
Once FST generation and pruning succeed, branching is now the major
remaining issue causing computations to take so long they are
abandoned.
The improved branching method described in
Section~\ref{sec:pseudo-costs} significantly improves performance
during the branching phase of the computation for most instances.

\subsection{Detailed Computational Results}
\label{sec:detailed-computational-results}

Figure~\ref{fig:cdf-new-iterations} plots the cumulative distribution
function (CDF) of the number of optimize/separate iterations needed
for the new code to obtain the optimal subtour relaxation --- across
all 1092 instances whose solution was attempted.

Figure~\ref{fig:cdf-new-time} plots the CDF of the CPU time needed for
the new code to obtain the optimal subtour relaxation --- across all
1092 instances whose solution was attempted.
The more diffuse range of $x$ values here (as compared with
Figure~\ref{fig:cdf-new-iterations}) demands the use of a logarithmic
$x$ axis, and results from the fact that all instance sizes are
included in the CDF.
These plots demonstrate that the number of iterations required by the
new code remains low despite increasing problem sizes and solution
times.

Figure~\ref{fig:cdf-both-iterations} plots the cumulative distribution
function (CDF) of the number of optimize/separate iterations needed
for both the old and new code to obtain the optimal subtour relaxation
--- across the 971 instances attempted by both the old and new code.
The new code converged all subtours for 967 of these instances,
whereas the old code only converged 904 of these instances within the
24 CPU hour time limit.
This difference is visible as maximum $y$ value (probability)
achieved.
The huge numbers of iterations required on many instances by the old
code demands a logarithmic $x$ axis when comparing CDFs of both the
old and new code.
This plot illustrates the large reductions in the number of iterations
achieved by the new ``strengthening via augmentation'' methods, and
the correspondingly dramatic improvement in convergence rate.

Figure~\ref{fig:cdf-both-time} plots the CDFs of CPU time needed to
converge subtours for both the old and the new code (across the 971
instances attempted by both the old and new code).
Figure~\ref{fig:cdfs-random-all-metrics} plots the CDFs of the number
of iterations needed to converge subtours for random instances with
the new code.
This plot has one CDF for each of the four metrics, and for each
instance size $1000 \le n \le 10000$.
At most 65 iterations were required (6000 point hexagonal instance
15), with all other instances needing at most 11 to 51 iterations.
Note that each CDF displays a small variance (curves are tall and
narrow).

Tables~\ref{tab:mean-random-e} through~\ref{tab:mean-random-u4}
present mean statistics (and standard deviations) for the
random instances as a function of instance size.
There is one table for each metric.
The standard deviations appear as a $\pm$ entry below each
corresponding mean.
The ``Iter'' fields give the mean number of LP/separate iterations
needed to satisfy all subtour inequalities at the root node.
The ``Relax Time'' fields give the mean CPU time needed to satisfy the
all subtour inequalities at the root node.
The ``Nodes'' fields gives the mean number of branch-and-bound nodes
generated.
The ``Done'' fields give C/N, where C is the number of instances
completed within the time limit, and N is the total number of
instances attempted.
(Instances that did not complete are omitted from the means and
standard deviations.)
The ``Concat Time'' fields give the mean CPU time needed for FST
Concatenation.

Table~\ref{tab:summary-tsplib-relax} summarizes the subtour relaxation
statistics (iterations and time) for the TSPLIB instances, all metrics.
Table~\ref{tab:summary-tsplib-concat} summarizes the FST Concatenation
statistics (nodes and time) for the TSPLIB instances, all metrics.

Tables~\ref{tab:concat-orlib-e} through~\ref{tab:concat-tsplib-2-u4}
present detailed computational results for the FST Concatenation phase
of each instance.
These tables normally contain two rows per instance%
\footnote{The exceptions are for Tables~\ref{tab:concat-rp-5-r}
  and~\ref{tab:concat-huge-1-r}
  which contain only results for the new code.
  These instances are so far beyond the capabilities of the old code
  that computations were not even attempted.}:
the first row are the results from GeoSteiner version 5.3;
the second row are the results for the new experimental version of
GeoSteiner as described in Section~\ref{sec:geosteiner-mods}.
To ease identifying pairs of corresponding old/new rows, the first row
of each pair (old code) is shifted slightly to the left, while the
second row of each pair (new code) is shifted slightly to the right.

Note that GeoSteiner dynamically generates two types of constraints:
subtours and local cuts.
Local cuts are only attempted when the current LP solution satisfies
all of the subtour inequalities.
In light of this, we partition FST concatenation into three distinct,
consecutive computational phases:
\begin{itemize}
  \item The ``Subtour Relaxation phase,'' during which only violated
    subtour inequalities are generated.  This phase ends when all
    subtour inequalities are satisfied.
  \item The ``Root Node phase,'' during which local cuts and
    subtours are generated in an interleaved fashion (local cuts are
    attempted for every LP solution that satisfies all subtours).
    This phase ends when all subtours are satisfied and no more
    violated local cuts are identified.
    Note that the gap is reduced substantially (in most cases) during
    this phase because: (1) local cuts substantially increase the LP
    lower bound; and (2) the improved LP solutions can in turn elicit
    improved upper bounds from the primal heuristic.
  \item The ``Branching phase,'' during which branch-and-cut is used
    to obtain integer solutions from the LP solutions obtained at each
    branch-and-bound node.
\end{itemize}

Tables \ref{tab:concat-orlib-e} through \ref{tab:concat-tsplib-2-u4}
report the following fields from the Subtour Relaxation phase:
{\em Iter} is the number of LP/separate iterations needed to
    satisfy all of the subtour inequalities;
{\em NZ/Col} is the average number of non-zero coefficients
    per column in the LP A-matrix;
{\em CC} is the size (in vertices) of the largest ``congested
    component,'' (maximum size of sub-problems processed by the
    $O(n^4)$ separation oracle for subtour inequalities);
{\em Gap} is $10^6 (ub - lb) / ub$, where $ub$ is either:
	(1) the integer optimal objective value (if known), or
	(2) the best upper bound known when the computation was
	    terminated, and $lb$ is the LP objective value of the
            Subtour Relaxation (this is in ``parts-per-million'');
{\em Frac} is the number of fractional variables in the final
    Subtour Relaxation LP solution;
{\em Time} is the CPU time needed to initially satisfy all
    subtour inequalities.

The following fields are reported from the Root Node phase:
{\em Iter} is the number of {\em additional} LP/separate
    iterations needed (after Subtour Relaxation) to further improve
    the LP lower bound using both local cuts and subtours;
{\em Rnds} is the number of ``rounds'' of local cut generation
    (i.e., number of invocations of the local cut separation oracle
    that returned one or more violated local cuts);
{\em Cuts} is the total number of local cut inequalities
    generated;
{\em NZPC} is the average number of non-zero coefficients per
    local cut generated;
{\em Gap} is $10^6 (ub - lb) / ub$, where $ub$ is either:
	(1) the integer optimal objective value (if known), or
	(2) the best upper bound known when the computation was
	    terminated,  and $lb$ is the final LP objective value of
            the root node (this is in ``parts-per-million'');
{\em Frac} is the number of fractional variables in the final
    LP solution of the Root Node;
{\em Time} is the {\em additional} CPU time (after Subtour
    Relaxation) needed to finish processing of the Root Node.

The following fields are reported from the Branching phase:
{\em Nodes} is the number of branch-and-bound nodes created;
{\em Proc} is the number of branch-and-bound nodes processed
    (excludes nodes that are cut off before being processed);
{\em Time} is the CPU time expended during the the Branching
    phase.

The {\em Total Time} field is the total CPU time expended during FST
Concatenation.

The breakthrough is clearly visible in these tables:
\begin{itemize}
  \item There are huge reductions in the number of iterations needed
    during the Subtour Relaxation phase.  Instances that previously
    required huge numbers of iterations / oracle calls now converge
    with relatively few iterations.
  \item On large instances, the NZ/Col field usually indicates an LP
    A-matrix that is much more sparse with the new code.  Such LPs
    solve much more quickly, and often result in solution vectors $x$
    that are much more numerically ``clean.''
  \item Instances often show greatly reduced CC values with the new
    code.
    Consider that this is the size of the largest instances given to
    the $O(n^4)$ subtour inequality separation oracle.  A 10x decrease
    in CC value can therefore represent up to a 10,000x reduction in
    CPU time for this separation oracle.
  \item The CPU times are dramatically smaller for most of the
    instances that the old code took the most time to solve.
\end{itemize}

One should also note the great utility of local cuts: (1) they often
provide a significant reduction of the gap; and (2) they are typically
very sparse inequalities.

Note that TSPLIB instance {\ttfamily p654} is the smallest instance
that has not yet been solved in all of the four metrics studied herein.

Tables \ref{tab:fsts-orlib-e} through \ref{tab:fsts-tsplib-u4} present
computational results for the FST generation and FST pruning phases.
Under FST Generation:
the FSTs field gives the total number of FSTs generated;
the Time field gives the FST generation CPU time.
(Euclidean FSTs report an additional Eq-Points field that gives the
total number of equilateral points generated.)
Under FST Pruning:
the Pruned field gives the total number of FSTs that were pruned
(eliminated);
the Required field gives the number of FSTs that were marked
``required,'' and {\em must} appear in every optimal solution;
the FSTs field gives the total number of FSTs that remain after
pruning (all 2-vertex MST edges are retained, even if pruned, as these
allow easy completion of a partial tree solution);
the Time field gives the FST pruning CPU time.
The Total Time field gives the total CPU time for both the FST
generation and FST pruning phases.

Note that the uniform directions FST generator (for hexagonal and
octilinear metrics) does not contain any of the numeric stability
improvements used in the Euclidean FST generator.
As a result, numeric errors during FST generation for the hexagonal
and octilinear metrics sometimes result in FSTs that are geometrically
inconsistent.
These geometric inconsistencies manifest as two logical
inconsistencies during the FST pruning phase:
(1) FSTs that get marked as both {\em pruned} and {\em required}; and
(2) final pruned FSTs that do not span all terminals.
(One or more terminals have no incident FSTs.)
The FST pruning code detects these situations, reports the error and
terminates unsuccessfully.
When this happens, we instead use the (inconsistent) FSTs without
pruning --- and appropriately annotate this occurrence for each such
instance in the tabulated results.
Since these numeric errors get worse with instance size, we do not
attempt hexagonal or octilinear solutions for random instances above
10,000 terminals.

\subsection{Large vs Small Subtours}
\label{sec:large-vs-small}

The above computational results establish the following:
\begin{eqnarray*}
	\hbox{Small subtours}\phantom{\hbox{+ Large subtours}}
&=&		\hbox{Good} \\
	\hbox{Small subtours} + \hbox{Large subtours}
&=&		\hbox{Very Good}
\end{eqnarray*}
To more completely confirm the accuracy of the EPR and CD indicators,
we now present computational evidence that:
(1) large subtours are stronger than small subtours; and
(2) large and small subtours together (as used in the main
computational experiments described herein) are even stronger.
Ideally we could just:
take a set $I$ of instances;
solve the $I$ using only large subtours;
solve the $I$ using only small subtours; and
compare the results.
Unfortunately it is not that simple.

GeoSteiner versions 5.3 and earlier demonstrate that
problems can be successfully solved using only relatively small
subtours.
However, it is quite difficult in practice to solve
problems using {\em only} relatively large subtours.
The number of iterations required explodes, and many small violated
components survive many consecutive iterations.
Some number of small subtours appear necessary to
achieve ``reasonable'' convergence.

\subsubsection{Code Changes to Explore Large vs Small Subtours}
\label{sec:code-changes-large-vs-small}
To address these issues, we prepared a special variant of the
GeoSteiner code (based upon the ``latest'' experimental code described
above) having the following properties:
(1) It stops upon converging all subtours;
(2) A \LargeMode{} operating mode that favors large subtours;
(3) A \SmallMode{} operating mode that favors small subtours;
(4) All previous ``bias'' of small subtours over large subtours (or
vice versa) is removed, to the extent possible;
(5) It reports the number of optimize/separate iterations needed to
converge all subtours;
(6) It reports a histogram that for each subtour size $k$ gives the
number of distinct subtour inequalities (on $k$ vertices) that enter
the LP during computation of the subtour relaxation.
(A subtour is counted only once, no matter how many times it enters or
leaves the LP.)
To compare these with the the main experimental code described in
Section~\ref{sec:geosteiner-mods}, we also added a \BothMode{} mode
that generates subtours normally with this same instrumentation.
Details of these code changes are as follows.

Code that seeds the initial constraint pool with:
(1) all distinct 2-vertex subtour inequalities, and
(2) inequalities that prohibit {\em incompatible} pairs of FSTs (as
determined during FST pruning) are disabled, as these favor small
subtours.

All separation algorithms were disabled {\em except} the main
$O(n^4)$ deterministic flow separator.
(Examples of removed separation algorithms include those for connected
components (cuts of weight zero), integral cycles, exhaustive
enumeration of small cardinality subsets, and a heuristic flow
formulation.)
The reductions normally used to decrease the size of instances going
into the main flow separator induce a heavy bias toward small
subtours.
We eliminate this bias as follows.
Let $C$ be the set of all reduced components
$H_i(x) = (V_i, E_i, x_i)$ for $1 \le i \le |C|$.
From $C$ we compute a single component
$H'(x) = (V',E',x')$
consisting of the {\em union} of all components
$H_i(x) = (V_i, E_i, x_i) \in C$.
We select a $v \in V'$ that minimizes $x(\delta(v))$ and construct
the flow network whose optimal solution yields the maximally violated
subtour $S$ such that $v \in S$ (as described in~\cite{Warme98}).
Let $S$ be the maximally violated subtour obtained.
If $S$ is violated, $S$ is added to the constraint pool and then
subjected to the constraint strengthening procedures described above.
Vertex $v$ is then removed from the problem as follows.
(1) Partition $C$ into sets $C_1$ (none of which contain vertex
$v$) and $C_2$ (all of which contain vertex $v$).
($C_2$ can contain more than one component, e.g., when $v$ is an
articulation vertex separating two or more biconnected components.)
(2) Let $C = C_1$.
(3) For each component $H_i \in C_2$, delete vertex $v$ from $H_i$ and
then perform all reductions on the modified $H_i$.
All remaining sub-components (if any) of $H_i$ are appended to $C$.
(4) If $C$ is empty, stop.
Otherwise, iterate back to the step wherein the union of $C$ is
computed and $v$ is selected.

One might ask whether use of the reductions (whether in union form or
not) retains any bias toward small subtours.
Computationally, it certainly finds large violated $S$.
For example, during early optimize/separate iterations this algorithm
identifies violated $S$ that are quite large (most of the vertices),
with violations of magnitude 70 or more.
(Violations of this magnitude are achieved by the combined effects of
numerous smaller violations.)
This is evidence that this new, ``unbiased'' separator is able to
truly maximize the magnitude of violation over the entire problem.
We argue this more formally by noting that the connected and
biconnected component reductions never discard any edges (nor any of
their indicent vertices), but merely partition the problem into
smaller components --- and these transformations are completely reversed
by taking the union of the reduced components.
This leaves single-vertex reduction, which iteratively deletes any
vertex $v$ for which $x(\delta(v)) \le 1$.
Let $S$ be a maximally-violated subtour,
and $v \in S$ such that $x(\delta(v)) < 1$.
Then removing $v$ from $S$ causes subtour $x(S) \le |S|-1$ to become
{\em more} violated because the right-hand-side decreases by 1,
whereas the left-hand-side decreases by an amount strictly less than 1.
Therfore, such $v$ cannot be in any maximally-violated subtour.
If on the other hand we have $x(\delta(v)) = 1$, then deleting $v$
imposes no bias because $v$ can be immediately restored during
single-vertex augmentation.

When adding subtour $S$ to the constraint pool, the code still also
adds subtour $V-S$ (whether violated or not).
This can be viewed as yet another measure to avoid bias between small
versus large subtours.

The final change concerns the new SOS heuristic
(Section~\ref{sec:SOS-heuristic}).
On each scan of the pool, and depending upon the mode of operation
(\LargeMode{} versus \SmallMode{}) this heuristic
was modified to selectively favor pulling violated constraints of
large or small cardinality from the constraint pool.
Instead of processing all violated constraints, it instead processes a
``filtered'' list of constraints.

According to the EPR indicator for $n=1000$, the weakest subtour is
at $k=342$.
The 20-th percentile strength subtours (80 percent of minimum
$\log_{10}(epr)$ value) occur at $k=146$ and $k=597$.
In \SmallMode{} mode, the filter takes all violated subtours
having at most $\lfloor 0.342|V| + 0.5\rfloor$ vertices.
If no such violations exist, then the filter takes all violated
subtours having at most $\lfloor 0.597|V| + 0.5 \rfloor$ vertices.
Otherwise it takes just one violated subtour (the first having the
fewest vertices).
In \LargeMode{} mode, the filter takes all violated subtours
having more than $\lfloor 0.342|V| + 0.5\rfloor$ vertices.
If no such violations exist, then the filter takes all violated
subtours having more than $\lfloor 0.146|V| + 0.5 \rfloor$ vertices.
Otherwise it takes just one violated subtour (the first having the
most vertices).

\subsubsection{Large vs Small Subtour Computational Results}
\label{sec:large-vs-small-computational-results}

We generated 1000 new random instances having 1000 terminals each.
FSTs were generated and pruned for these instances using all four
distance metrics (Euclidean, rectilinear, hexagonal and octilinear).
Due to the aforementioned numeric instability of the uniform
directions FST generator, FST pruning failed for hexagonal instances
736 and 800, so we used the unpruned (geometrically inconsistent) FSTs
for these two instances.
We solved each of these 4000 instances using \BothMode{}, \LargeMode{}
and \SmallMode{} modes, producing
iteration count and size histogram for each run.
It should be noted that the \SmallMode{} computation for
hexagonal instance 112 was highly degenerate.
After 17613 iterations and more than 66 CPU hours, we killed this one
instance.
The other 11999 runs finished successfully.

\begin{figure}[!ht]
\begin{center}
 \includegraphics[width=4.0in,clip=]{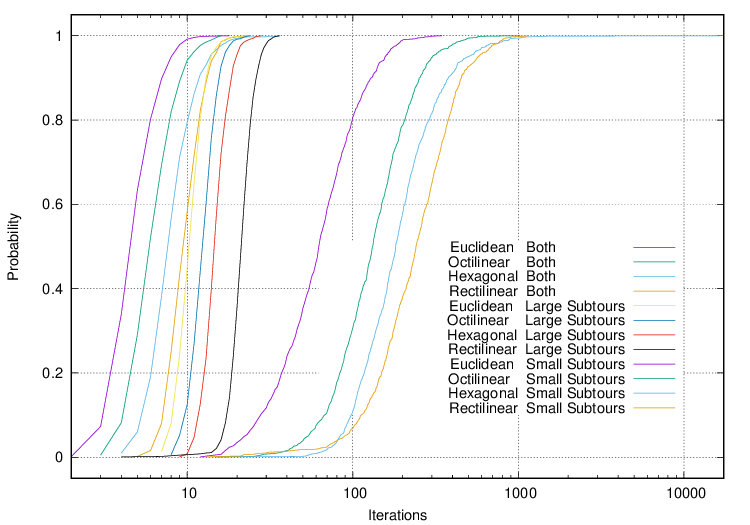}
 \captionof{figure}{CDF of optimize/separate iterations needed to
   converge all subtours, $n=1000$, all metrics, \BothMode{},
   \LargeMode{} and \SmallMode{} modes}
 \label{fig:cdf-both-large-small-subtours}
 \vspace*{\floatsep}
 \includegraphics[width=4.0in,clip=]{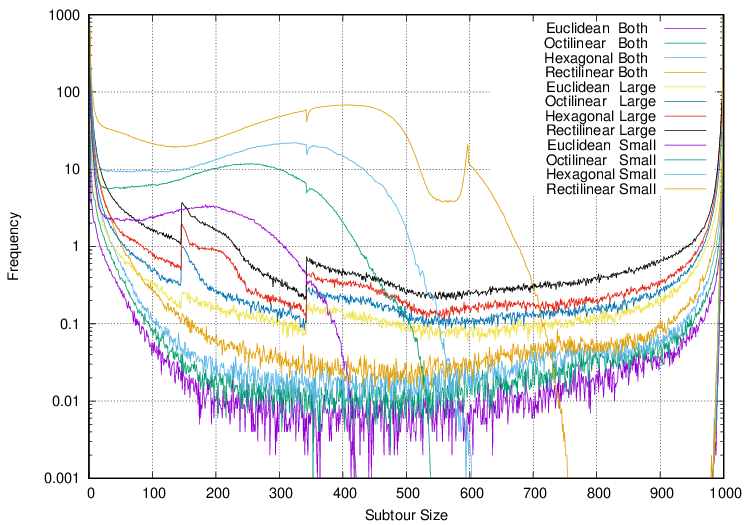}
 \captionof{figure}{Mean number of distinct $k$-vertex subtours ever
   present in LP during solution of subtour relaxation, $n=1000$,
   all metrics, \BothMode{}, \LargeMode{} and \SmallMode{} modes}
 \label{fig:histogram-both-large-small-subtours}
\end{center}
\end{figure}

Figure~\ref{fig:cdf-both-large-small-subtours}
plots the CDF of iterations needed to satisfy all subtours for each
metric for \BothMode{}, \LargeMode{} and \SmallMode{} runs.
The \BothMode{} runs produce the least number of iterations,
with \LargeMode{} producing slightly more.
The \SmallMode{} runs require substantially more iterations.
This plot also clearly shows a ranking of the metrics: Euclidean,
octilinear, hexagonal and rectilinear, with the Euclidean metric
converging most quickly, and rectilinear metric converging least
quickly.
Figure~\ref{fig:histogram-both-large-small-subtours}
plots the histogram data: for each subtour size $k$ the plot gives the
mean number of distinct $k$-vertex subtours ever present in the LP
during solution of the subtour relaxation.
(Artifacts in these plots for \LargeMode{} and \SmallMode{} modes are
clearly visible at the filter boundary subtour sizes
$k \in \lbrace 146,\, 342,\, 597 \rbrace$.

In Table~\ref{tab:histogram-both-large-small-subtours}
we tabulate the mean histogram data for the smallest and largest
several subtour sizes in order to:
(1) confirm the statistical difference between subtour size
distribution between \BothMode{}, \LargeMode{} and \SmallMode{} modes;
and
(2) highlight these data which are difficult to see at the left and
right edges of Figure~\ref{fig:histogram-both-large-small-subtours}.
\begin{sidewaystable}[htbp]
\begin{center}
 \caption{Mean number of distinct $k$-vertex subtours ever
   present in LP during solution of subtour relaxation, all metrics,
   \BothMode{}, \LargeMode{} and \SmallMode{} modes
   (small $\pm$ entries are standard deviations)}
 \begin{tabular}{|r|rrr|rrr|rrr|rrr|} \hline
 \multicolumn{1}{|c}{Subtour} &
 \multicolumn{3}{|c}{Euclidean} &
 \multicolumn{3}{|c}{Rectilinear} &
 \multicolumn{3}{|c}{Hexagonal} &
 \multicolumn{3}{|c|}{Octilinear} \\
 \multicolumn{1}{|c|}{Size} &
	{\small\BothMode} & {\small\LargeMode} & {\small\SmallMode} &
	{\small\BothMode} & {\small\LargeMode} & {\small\SmallMode} &
	{\small\BothMode} & {\small\LargeMode} & {\small\SmallMode} &
	{\small\BothMode} & {\small\LargeMode} & {\small\SmallMode} \\
	\hline
  2 &	185.39 & 44.43 & 153.24 &
	699.15 & 153.71 & 517.97 &
	651.93 & 148.76 & 474.42 &
	340.63 & 72.81 & 265.03 \NEWLINE
    &	\SDEV{31.269} & \SDEV{9.959} & \SDEV{25.224} &
	\SDEV{46.725} & \SDEV{16.530} & \SDEV{34.038} &
	\SDEV{44.529} & \SDEV{16.526} & \SDEV{28.956} &
	\SDEV{38.559} & \SDEV{12.927} & \SDEV{28.909} \\
  3 &	3.618 & 28.85 & 50.41 &
	52.33 & 128.13 & 208.02 &
	29.24 & 101.13 & 173.05 &
	10.62 & 48.60 & 99.26 \NEWLINE
    &	\SDEV{2.233} & \SDEV{7.380} & \SDEV{9.652} &
	\SDEV{8.661} & \SDEV{13.625} & \SDEV{17.923} &
	\SDEV{5.784} & \SDEV{11.798} & \SDEV{13.036} &
	\SDEV{3.643} & \SDEV{8.778} & \SDEV{11.938} \\
  4 &	4.161 & 23.24 & 37.23 &
	32.91 & 95.79 & 150.32 &
	15.00 & 72.16 & 111.40 &
	8.645 & 45.27 & 67.16 \NEWLINE
    &	\SDEV{2.166} & \SDEV{6.341} & \SDEV{7.457} &
	\SDEV{6.638} & \SDEV{11.863} & \SDEV{15.460} &
	\SDEV{4.117} & \SDEV{9.273} & \SDEV{11.799} &
	\SDEV{3.117} & \SDEV{8.306} & \SDEV{8.986} \\
  5 &	3.246 & 14.49 & 22.71 &
	22.528 & 60.62 & 105.02 &
	11.72 & 46.57 & 64.52 &
	6.156 & 27.47 & 38.65 \NEWLINE
    &	\SDEV{1.839} & \SDEV{4.535} & \SDEV{5.382} &
	\SDEV{5.138} & \SDEV{8.601} & \SDEV{14.117} &
	\SDEV{3.678} & \SDEV{6.928} & \SDEV{9.430} &
	\SDEV{2.544} & \SDEV{5.996} & \SDEV{6.768} \\
  6 &	2.605 & 11.12 & 18.53 &
	18.56 & 47.83 & 87.22 &
	10.92 & 35.71 & 47.48 &
	5.038 & 20.33 & 28.80 \NEWLINE
    &	\SDEV{1.664} & \SDEV{3.877} & \SDEV{4.420} &
	\SDEV{4.696} & \SDEV{7.642} & \SDEV{13.595} &
	\SDEV{3.571} & \SDEV{6.355} & \SDEV{7.933} &
	\SDEV{2.325} & \SDEV{5.726} & \SDEV{5.758} \\
  7 &	2.195 & 8.50 & 12.99 &
	16.36 & 36.05 & 72.28 &
	9.321 & 26.31 & 33.69 &
	4.380 & 15.03 & 20.59 \NEWLINE
    &	\SDEV{1.539} & \SDEV{3.265} & \SDEV{3.858} &
	\SDEV{4.222} & \SDEV{6.491} & \SDEV{13.161} &
	\SDEV{3.096} & \SDEV{5.191} & \SDEV{7.098} &
	\SDEV{2.185} & \SDEV{4.125} & \SDEV{5.169} \\
  8 &	2.05 & 7.36 & 11.09 &
	15.39 & 30.48 & 64.36 &
	8.078 & 21.81 & 27.65 &
	4.138 & 12.70 & 16.92 \NEWLINE
    &	\SDEV{1.505} & \SDEV{3.035} & \SDEV{3.443} &
	\SDEV{4.262} & \SDEV{6.007} & \SDEV{12.521} &
	\SDEV{3.003} & \SDEV{4.948} & \SDEV{6.671} &
	\SDEV{2.084} & \SDEV{3.768} & \SDEV{4.759} \\
  9 &	2.078 & 6.20 & 8.37 &
	14.11 & 25.28 & 56.95 &
	7.531 & 17.39 & 22.29 &
	4.019 & 10.44 & 13.57 \NEWLINE
    &	\SDEV{1.500} & \SDEV{2.742} & \SDEV{3.227} &
	\SDEV{4.114} & \SDEV{5.430} & \SDEV{12.353} &
	\SDEV{2.967} & \SDEV{4.344} & \SDEV{6.093} &
	\SDEV{2.057} & \SDEV{3.477} & \SDEV{4.356} \\
 10 &	1.954 & 5.49 & 7.31 &
	12.94 & 22.14 & 52.51 &
	6.876 & 14.98 & 19.37 &
	3.811 & 9.18 & 11.95 \NEWLINE
    &	\SDEV{1.486} & \SDEV{2.742} & \SDEV{3.079} &
	\SDEV{3.845} & \SDEV{5.085} & \SDEV{12.302} &
	\SDEV{2.856} & \SDEV{3.914} & \SDEV{5.910} &
	\SDEV{2.008} & \SDEV{3.204} & \SDEV{4.085} \\
	\hline
991 &	0.521 & 4.12 & 0.016 &
	2.043 & 11.78 & 0.100 &
	1.445 & 10.32 & 0.091 &
	0.8749 & 6.95 & 0.039 \NEWLINE
    &	\SDEV{0.7402} & \SDEV{2.190} & \SDEV{0.125} &
	\SDEV{1.485} & \SDEV{3.562} & \SDEV{0.307} &
	\SDEV{1.259} & \SDEV{3.146} & \SDEV{0.305} &
	\SDEV{0.9199} & \SDEV{2.686} & \SDEV{0.204} \\
992 &	0.6256 & 4.77 & 0.032 &
	2.390 & 14.39 & 0.181 &
	1.598 & 12.66 & 0.188 &
	1.104 & 8.20 & 0.065 \NEWLINE
    &	\SDEV{0.8734} & \SDEV{2.410} & \SDEV{0.176} &
	\SDEV{1.594} & \SDEV{3.883} & \SDEV{0.420} &
	\SDEV{1.315} & \SDEV{3.436} & \SDEV{0.448} &
	\SDEV{1.047} & \SDEV{3.010} & \SDEV{0.258} \\
993 &	0.6977 & 5.61 & 0.065 &
	2.809 & 18.17 & 0.309 &
	1.864 & 16.30 & 0.325 &
	1.117 & 9.98 & 0.131 \NEWLINE
    &	\SDEV{0.8795} & \SDEV{2.548} & \SDEV{0.277} &
	\SDEV{1.722} & \SDEV{4.467} & \SDEV{0.544} &
	\SDEV{1.459} & \SDEV{4.040} & \SDEV{0.574} &
	\SDEV{1.078} & \SDEV{3.312} & \SDEV{0.363} \\
994 &	0.7568 & 6.79 & 0.146 &
	3.424 & 23.64 & 0.723 &
	2.318 & 21.10 & 0.700 &
	1.351 & 12.34 & 0.274 \NEWLINE
    &	\SDEV{0.9114} & \SDEV{2.792} & \SDEV{0.388} &
	\SDEV{1.956} & \SDEV{5.264} & \SDEV{0.846} &
	\SDEV{1.571} & \SDEV{4.662} & \SDEV{0.849} &
	\SDEV{1.177} & \SDEV{3.895} & \SDEV{0.528} \\
995 &	0.9750 & 8.60 & 0.298 &
	4.213 & 31.66 & 1.45 &
	2.963 & 28.77 & 1.52 &
	1.668 & 16.28 & 0.607 \NEWLINE
    &	\SDEV{1.070} & \SDEV{3.442} & \SDEV{0.560} &
	\SDEV{2.062} & \SDEV{6.068} & \SDEV{1.241} &
	\SDEV{1.705} & \SDEV{5.726} & \SDEV{1.262} &
	\SDEV{1.307} & \SDEV{4.397} & \SDEV{0.758} \\
996 &	1.298 & 11.62 & 0.598 &
	5.664 & 45.60 & 3.50 &
	3.816 & 40.87 & 3.59 &
	2.153 & 22.92 & 1.52 \NEWLINE
    &	\SDEV{1.199} & \SDEV{4.117} & \SDEV{0.800} &
	\SDEV{2.338} & \SDEV{7.484} & \SDEV{1.854} &
	\SDEV{2.068} & \SDEV{6.655} & \SDEV{1.932} &
	\SDEV{1.494} & \SDEV{5.545} & \SDEV{1.292} \\
997 &	1.700 & 16.20 & 1.13 &
	8.559 & 69.45 & 5.19 &
	5.669 & 63.55 & 5.84 &
	2.900 & 30.81 & 2.24 \NEWLINE
    &	\SDEV{1.339} & \SDEV{5.047} & \SDEV{1.065} &
	\SDEV{3.077} & \SDEV{9.181} & \SDEV{2.387} &
	\SDEV{2.248} & \SDEV{8.977} & \SDEV{2.409} &
	\SDEV{1.719} & \SDEV{6.311} & \SDEV{1.516} \\
998 &	2.488 & 23.54 & 4.15 &
	11.32 & 80.25 & 19.12 &
	9.510 & 108.06 & 19.86 &
	4.374 & 45.64 & 7.81 \NEWLINE
    &	\SDEV{1.604} & \SDEV{6.179} & \SDEV{2.224} &
	\SDEV{3.358} & \SDEV{9.455} & \SDEV{4.596} &
	\SDEV{3.235} & \SDEV{12.130} & \SDEV{4.642} &
	\SDEV{2.085} & \SDEV{8.458} & \SDEV{2.938} \\
999 &	929.07 & 929.07 & 929.07 &
	972.25 & 972.25 & 972.25 &
	978.74 & 979.76 & 979.76 &
	951.69 & 951.69 & 951.69 \NEWLINE
    &	\SDEV{8.891} & \SDEV{8.891} & \SDEV{8.891} &
	\SDEV{6.213} & \SDEV{6.213} & \SDEV{6.213} &
	\SDEV{31.36} & \SDEV{4.965} & \SDEV{4.968} &
	\SDEV{7.416} & \SDEV{7.416} & \SDEV{7.416} \\
	\hline
 \end{tabular}
 \label{tab:histogram-both-large-small-subtours}
\end{center}
\end{sidewaystable}

The subtours of size $k=999$ are equivalent to $x(\delta(v)) \ge 1$
for each $v \in V$, and
all of these are included in the initial LP instance.
Although one might expect to see exactly 1000 such subtours for each
instance, fewer actually appear because distinct vertices $v$ can
yield the same inequality, and such duplicate rows are automatically
detected and omitted by GeoSteiner.

Although both small and large subtours are used in both
\LargeMode{} and \SmallMode{} modes,
Figure~\ref{fig:histogram-both-large-small-subtours} and
Table~\ref{tab:histogram-both-large-small-subtours} show the statistical
mix of subtours clearly leaning in favor of large and small
subtours, respectively.
This statistical difference in mix of subtour sizes is sufficient to
cause the substantial differences in iteration counts shown in
Figure~\ref{fig:cdf-both-large-small-subtours}.

\begin{figure}[!ht]
\begin{center}
 \includegraphics[width=4.0in,clip=]{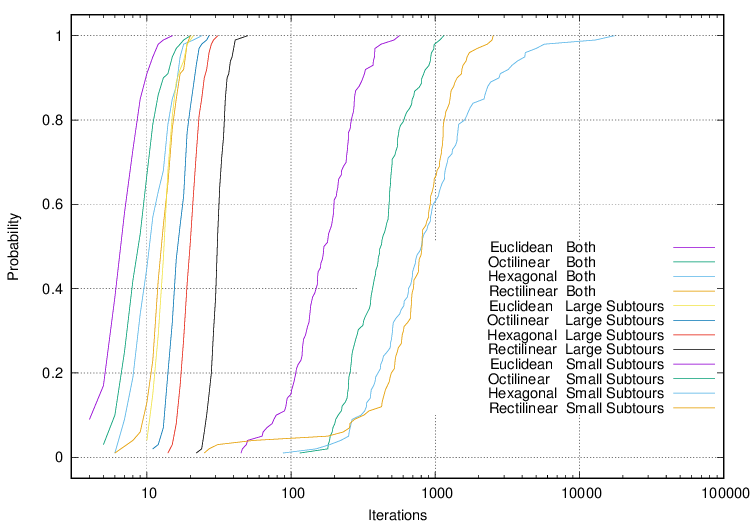}
 \captionof{figure}{CDF of optimize/separate iterations needed to
   converge all subtours, $n=2000$, all metrics, \BothMode{},
   \LargeMode{} and \SmallMode{} modes}
 \label{fig:cdf-2k-both-large-small-subtours}
 \vspace*{\floatsep}
 \includegraphics[width=4.0in,clip=]{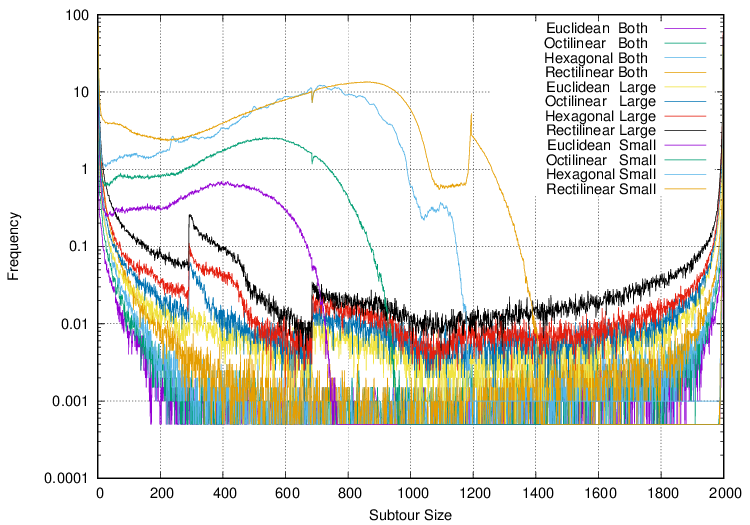}
 \captionof{figure}{Mean number of distinct $k$-vertex subtours ever
   present in LP during solution of subtour relaxation, $n=2000$,
   all metrics, \BothMode{}, \LargeMode{} and \SmallMode{} modes}
 \label{fig:histogram-2k-both-large-small-subtours}
\end{center}
\end{figure}

To assess the trend as instance size increases, we repeated this
procedure at $n=2000$ terminals.
These computations were significantly more costly than for $n=1000$
(especially for \SmallMode{} mode with rectilinear and hexagonal
metrics), so we used only 100 random instances.
The \SmallMode{} rectilinear computation for instance \#20 was highly
degenerate.
We killed this instance after 13,128 optimize/separate iterations, as
there had been no improvement in the lower bound for more than 8000
consecutive iterations.
Being unable to remove slack rows, the LP had grown to more than
635,000 rows.
This data point is therefore missing from the plots, but would make
the CDF (and histogram) for \SmallMode{} rectilinear instances be even
heavier-tailed.
For $n=2000$,
Figure~\ref{fig:cdf-2k-both-large-small-subtours} plots the CDFs of
iterations needed to satisfy all subtour inequalities, and
Figure~\ref{fig:histogram-2k-both-large-small-subtours} plots the
histograms of mean number of $k$-vertex subtours ever present in the
LP during solution of subtour relaxation.
Note that \SmallMode{} mode in particular requires many more
iterations at $n=2000$
(Figure~\ref{fig:cdf-2k-both-large-small-subtours})
than for $n=1000$
(Figure~\ref{fig:cdf-both-large-small-subtours}),
whereas the CDFs remain quite similar for \BothMode{} mode.

This is additional evidence that subtours of large cardinality are
actually stronger (computationally) than subtours of small
cardinality, as predicted by both the EPR and CD indicators.
These computational results also support the predictions
in~\cite{WarmeIndicators1}
(based upon combined effect of EPR/CD indicators and interior angles
between pairs of subtours) that combining subtours of both large and
small cardinality should have very strong computational effect.

\subsection{Steiner Tree in Hypergraph Results}
\label{sec:results-steiner-tree-in-hypergraph}

Although not included in~\cite{WarmeIndicators1}, the EPR and CD for
subtours of the Steiner tree in hypergraph problem (undirected
extended formulation~\cite{GoemansMyung}) yield curves virtually
indistinguishable from the subtours of STHGP.
It should therefore be expected that strengthening via
augmentation would provide a similar speedup on this problem.

To test this, we modified the (very old) version of GeoSteiner
presented in in~\cite{ZachariasenWinter1999} to use augmentation.
Since there are no recognized test instances for the Steiner tree in
hypergraphs problem, we tested this on the Steiner tree in graph
instances from the OR-Library~\cite{OR-library}.

Similar speedups were observed to those reported here for STHGP.
We present no detailed results from these computations, however.
State-of-the-art methods for the Steiner tree in graphs problem
(such
as~\cite{RehfeldtKoch2023}, \cite{Polzin2003}
and~\cite{VahdatiDaneshmand2004})
depend
heavily upon reductions of the graph instance prior to applying
branch-and-cut.
GeoSteiner has no such graph reductions:
our computations were performed on raw, unreduced OR-Library graph
instances.
A highly significant speedup was observed, but the computational
results are not competitive with state-of-the-art.
This provides additional evidence that EPR and CD are highly
correlated with computational strength.

Although the directed and undirected formulations yield the same lower
bound, there has been concensus that the directed formulation (for the
Steiner tree in graphs problem) performs better than the undirected
formulation (used in these experiments).
We conjecture this was largely due to the ``nested cuts'' technique of
Koch and Martin~\cite{KochMartin1998}, for which no corresponding
technique was known for the undirected formulation.
The methods described herein narrow this gap, if not eliminate it.
It has also been observed that the directed formulation seems faster
on instances having more Steiner vertices than terminals, whereas the
undirected formulation seems faster in the opposite case.
More study is certainly warranted.

\section{Discussion}
\label{sec:discussion}

We now discuss various observations and ramifications regarding the
above.

\subsection{The 100,000 Terminal Instances}
\label{sec:100k-instances}

A capstone result in~\cite{JuhlWarmeWinterZachariasen2018} was optimal
solutions for three (out of fifteen) of the 100,000 terminal random
Euclidean instances.
Average FST concatenation time using GeoSteiner 4.0 on that hardware
platform (for the three instances that solved) was 10,370 seconds.
In the present study, GeoSteiner 5.3 solves the same three instances
(FST concatenation phase) in an average of 2880 seconds (due solely to
a much faster CPU).
With the new constraint strengthening methods described herein, all
fifteen of these 100,000 terminal instances now successfully
concatenate in an average of 532.11 seconds.

The difficulty of the twelve instances that fail to solve using
GeoSteiner 5.3 and earlier versions cannot be understated.
After the 14th DIMACS Challenge Workshop
(where~\cite{JuhlWarmeWinterZachariasen2018} was originally presented),
this author attempted to solve 100,000 terminal instance 1 using
GeoSteiner 5.0 on a 5.0 GHz AMD 9590.
After more than 33 CPU {\em months} and 22,846 optimize/separate
iterations, this computation was still slowly finding violated
subtours at the root node.
Although solving three of these instances was quite significant, the
optimize/separate loop clearly exhibits extremely high variance.
The present methods eliminate these difficulties.

\subsection{The 1,000,000 Terminal Instance}
\label{sec:1m-instance}

The optimal solution presented here of a 1,000,000 terminal random
Euclidean Steiner tree instance breaks new ground in the solution of
classical NP-hard optimization problems.
The formulation has 1,402,572 zero-one variables, yet produces only 10
branch-and-bound nodes.
(Two of these are cutoff before being processed.)
Because of the new constraint strengthening procedures (using
reduction together with augmentation), only 163 optimize/separate
iterations are needed to satisfy all of the subtour inequalities.
The gap between the subtour relaxation and the optimal solution is a
mere 1717.037 parts-per-{\em billion}.
Local cuts (393 cuts in 26 rounds) further reduce this gap to just
128.931 parts-per-billion.

Although (at least) NP-hard, uniformly random Euclidean instances are
unusually easy for their size.
The subtour relaxation gaps are very small, and local cuts
significantly reduce these gaps.
They have virtually no dual degeneracy.
When the terminal coordinates are rational or integer, the FST
lengths are algebraic numbers of the form
$\sqrt{a \,+\, b\,\sqrt{3}}$,
where $a$ and $b$ are rational.
Randomized terminal coordinates make integer relationships among
algebraic numbers of this form highly unlikely.
With such a tightly constrained instance and no dual degeneracy, the
perturbation of branching on fractional variables produces substantial
lower bound progress in each child node, quickly closing the narrow
gap.

In contrast, on random rectilinear instances, the largest instance
presently solved with this code is 40,000 terminals.
Rectilinear instances (even random instances) produce heavy dual
degeneracy and large branch-and-bound trees.
Because of the dual degeneracy, the gaps are larger, and the
probability is high that branching on a single fractional variable
yields no progress on the lower bound in at least one of the
children.
% Although having integer FST lengths is an issue, the bigger factor is
% that between two terminals $p$ and $q$ there are so many different
% paths having exactly the same length, even when discretized onto FSTs.

\subsection{Structure of Strengthed Subtours}
\label{sec:strengthened-subtour-structure}

Each violated subtour $S \subset V$ emerging from the flow-separator
is strengthened via reduction (yielding violated subtours
$R_1,\, R_2,\, \ldots,\, R_j$, where
$R_i \subset S$ for $1 \le i \le j$).
The new code also strengthens $S$ via the above augmentation rules
(yielding violated subtours
$A_1,\, A_2,\, \ldots,\, A_k$, where
$S \subset A_i$ for $1 \le i \le k$).

Let $V$ be the vertices and $n=|V|$.
Consider the lattice $L$ of all subsets of $V$.
This lattice has $n+1$ layers -- one layer for each subset cardinality
$i$, $0 \le i \le n$.
The top layer consists of one node representing $V$ of cardinality $n$.
The bottom layer consists of one node representing $\emptyset$ of
cardinality 0.
Lattice $L$ contains a directed arc from node $U$ to $T$ if
$T \subset U$ and $|T|+1 = |U|$.

\begin{figure}[h]
\centering
\includegraphics[width=3.5in]{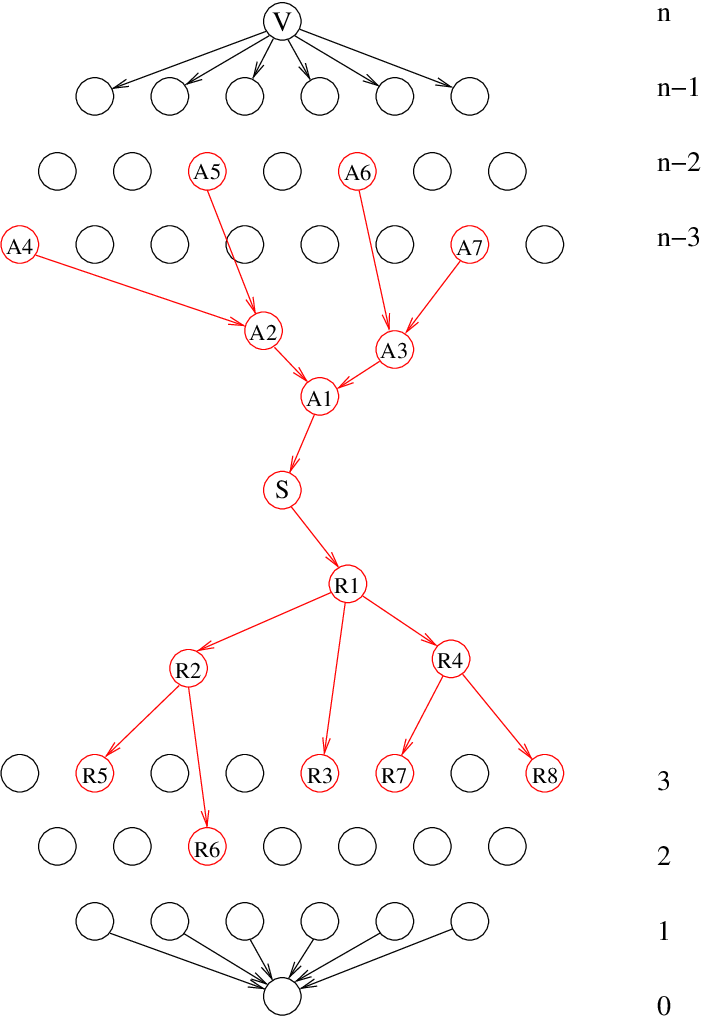}
\caption{Structure of strengthed subtours within lattice $L$ of
  subsets of $V$}
\label{fig:strengthen-lattice}
\end{figure}

Consider the sets
$R_1,\, R_2,\, \ldots,\, R_j,\, S,\, A_1,\, A_2,\, \ldots,\, A_k$
representing the violated subtours obtained by both reduction and
augmentation from $S$.
Each of these subsets corresponds to a node in $L$.
Connecting the violated subtour nodes with directed arcs (indicating
the superset relation, where these arcs are permitted to cross from
one layer to any lower cardinality layer) results in a tree structure
as illustrated in Figure~\ref{fig:strengthen-lattice}.

\FloatBarrier

Let $T$ be a ``top'' node of this tree (no tree arcs enter $T$).
Let $B$ be a ``bottom'' node of this tree (no tree arcs leave $B$).
The path from $T$ to $B$ is unique and always passes through $S$.

We experimentally tested the following algorithmic idea: generate only
the two ``strongest'' violated subtours along the path from $T$ to $B$
(as estimated by interpolating the corresponding subtour size onto the
EPR or CD indicator curve for $n=1000$).
The idea was to minimize the growth of the LP $A$ matrix.
This did not work well in practice because it significantly increased
the number of iterations needed to satisfy all subtours (loss off
convergence rate).

All of the violated subtours in this tree seem to be important for
achieving the fastest convergence of all subtours.
There seems to be some fundamental property of this tree structure (as
obtained by single-vertex reduction/augmentation (to identify violated
subtours $R_1$ and $A_1$) and the additional structure identified by
connected and biconnected components that leads to exquisitely rapid
convergence of all subtour inequalities.

We find this structure to be strongly analogous to the
{\em nested cuts} of Koch and Martin~\cite{KochMartin1998}.
They achieve a similar structure of nested subtours in the directed
formulation for the Steiner problem in graphs by forcing the weight of
all arcs crossing the minimum cut to 1 and then re-solving the
max-flow / min-cut problem, repeating until no further violated cuts
exist.
They then repeat this procedure with the original arc weights, but
with each arc reversed in direction.
Note that when separating cut inequalities for the directed
formulation of the Steiner problem in graphs, one can have multiple
alternate minimum cuts.
It is likely that procedures similar to single-vertex reduction
and single-vertex augmentation might exist for the directed
formulation that find the smallest and largest of these
equally violated subsets, analogous to the violated subtours
$R_1$ and $A_1$ obtained from $S$ by single-vertex reduction and
single-vertex augmentation, respectively.
More generally, similar reduction and augmentation rules might exist
for the directed formulation that emulate and exploit this same tree
structure of cuts to achieve even better convergence rate than
reported in~\cite{KochMartin1998}.

\subsection{Effectiveness of Separation Oracles}
\label{sec:oracle-effectiveness}

We have much computational experience indicating that the particular
choice of which violated inequalities are generated can cause
performance to vary by many orders of magnitude.
Consider the landmark result of
Gr{\"o}tschel, Lov{\'a}sz and Schrijver%
~\cite{GrotschelLovaszSchrijver},
which proved that so long as:
(1) the inequality class is finite; and
(2) the separation oracle is guaranteed to find at least one violated
inequality if any violations exist;
then the ellipsoid algorithm terminates using $O(d^2\,\log(R/r))$
calls to the separation oracle, where $d$ is the dimensionality of the
space, $R$ is the radius of a ball completely containing the polytope
$P$, and $r$ is the radius of a ball completely contained within the
polytope $P$.

By comparison, very little is known regarding bounds on the
number of iterations of the analogous loop wherein the simplex
algorithm is used to optimize over the current set of constraints.
The successive LP solutions $x$ have no relationship
whatsoever to the center points of the successive ellipsoids in the
ellipsoid algorithm, nor is there any corresponding constant reduction
factor as with the successive decrease of ellipsoid volumes.
No matter how similar the two algorithms might appear, nor how
theoretically attractive the iteration bound for the ellipsoid
algorithm, we know of no good upper bounds on the number of iterations
required by the ``more practically used'' loop that separates
successive LP solutions.

Computational experience (including that described herein) clearly
indicates that using a separation oracle that chooses ``strong''
violations over ``weaker'' violations can result in a huge reduction
in the number of iterations needed to satisfy the entire class of
inequalities.
In contrast, a separation oracle that chooses inequalities unwisely
(failing to find ``strong'' violations when they exist) can have
devastating computational consequences, taking so long to converge in
practice that computations are either abandoned, or some
``tailing off'' decision is made to invoke branching.

Note that proving some class of inequalities to be facet-defining
offers no assurance of rapid convergence --- the indicators presented
in~\cite{WarmeIndicators1} demonstrate that strength varies widely within
such classes, even for polyhedra having high levels of symmetry.
The attitude that,
``any violation in such a class is {\itshape good enough},''
should be actively discarded when designing separation algorithms for
real computations.

In the past, ``tuning'' a separation algorithm so that it converges
rapidly has often been at best an art, and even more often an
experimental, trial-and-error process.
We propose that strength indicators (such as those calculated
in~\cite{WarmeIndicators1}) should guide the design of separation
and constraint strengthening algorithms, resulting in a more
principled, empirically based approach to obtaining high-performance
separation algorithms.
% Researchers should favor mathematical insight over intuition or
% guesswork.

\subsection{Common Algorithmic Patterns Yielding Bias}
\label{sec:algorithmic-patterns}

There is a pattern here that is common for many algorithms
that separate a combinatorially large class of inequalities:
\begin{itemize}
 \item The underlying separation algorithm $f(x)$ is computationally
   costly.
 \item Mitigate this by first applying reductions.
 \item Apply $f(x)$ to (smaller) reduced sub-problems.
 \item Resulting violated inequalities are clustered within restricted
   portions of the inequality space.
\end{itemize}
This pattern can result in a form of ``algorithmic bias''
that systematically avoids large portions of the inequality space.
(This technical form of algorithmic bias should not be confused
with the ``social'' type of algorithmic bias being heavily discussed
in the machine learning community.)
%
% It is not unreasonable to assume that this same pattern occurs within
% algorithms for other problems, wherein portions of the inequality
% space containing very strong inequalities (with tremendous potential
% for performance gains) are systematically avoided by algorithmic bias
% within the separation algorithms.
%
The present example of STHGP subtours shows that this pattern can bias
the algorithm {\em away} from very strong --- or even the strongest
inequalities in the class.
Correcting this bias yields enormous performance gains in GeoSteiner.

This provides yet another incentive for calculating such strength
indicators for a particular polyhedron / inequality class ---
if this algorithmic bias systematically avoids regions of the
inequality space containing very strong inequalities, then one can
devise methods that specifically target these strong inequalities that
would otherwise be missed.
We are unaware of existing methods that intentionally counteract such
algorithmic forces.

\section{Conclusion}
\label{sec:conclusion}

For problems having exponential time complexity, increasing the
size of instances that can be practically solved by a factor 2 is
significant progress.
A factor of 10 improvement, however, is truly remarkable, especially
on algorithms and implementations that are already quite mature.
This improvement was not obtained by introducing any new class of
cuts, but by improving the convergence rate over an existing,
exponentially large class of cuts.

The code changes that yielded this performance gain were motivated
solely by the property of $\STHGP{n}$ (as predicted by the indicators
presented in~\cite{WarmeIndicators1}) that subtours of large cardinality
are even stronger than subtours of small cardinality.
This property of $\STHGP{n}$ was neither known nor suspected until
clearly revealed by these indicators.
The accuracy of these indicators as predictors of actual
computational strength is now unambiguously confirmed, as is their
usefulness for obtaining rapidly converging separation algorithms in an
empirical, principled manner.

\begin{credits}
\subsubsection{\ackname}
I am grateful to Martin Zachariasen, Pawel Winter, Jeffrey Salowe,
Karla Hoffman, Maurice Queyranne, Warren D. Smith and Bill Cook:
this work would not exist without their helpful discussions and
encouragement.
I also thank Group W for their support: this work has no direct
connection to their business.
Proverbs 25:2, Soli Deo Gloria.

\subsubsection{\discintname}
The author has financial and intellectual property interests in
GeoSteiner, and has filed for patent protection on ideas in this
paper.
\end{credits}

\bibliographystyle{splncs04}
\bibliography{phd,steiner,local}

\appendix

\section{Appendix --- Results}
\label{sec:appendix}

\renewcommand{\SDEV}[1]{{$\pm #1$}}

% \layout{}
 
%

\begin{figure}[!ht]
\begin{center}
\begin{minipage}[t]{2.25in}
\includegraphics[width=2.25in]{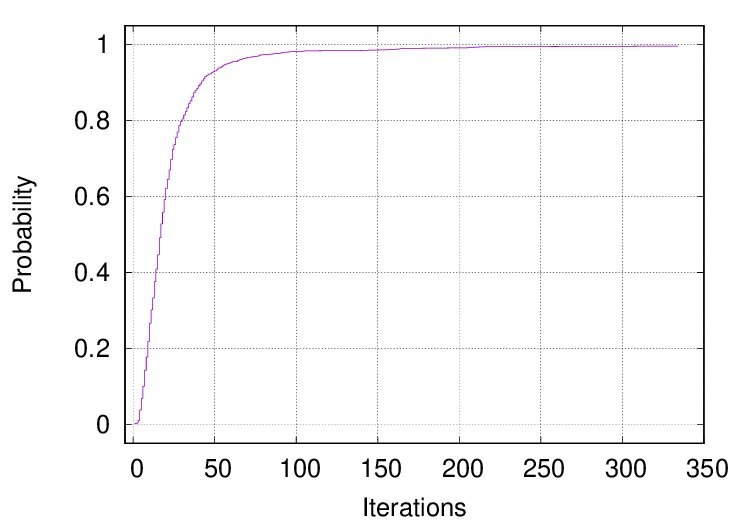}
 \captionof{figure}{CDF of Subtour Relaxation Iterations: New Code}
 \label{fig:cdf-new-iterations}
\end{minipage}
\quad
\begin{minipage}[t]{2.25in}
\includegraphics[width=2.25in]{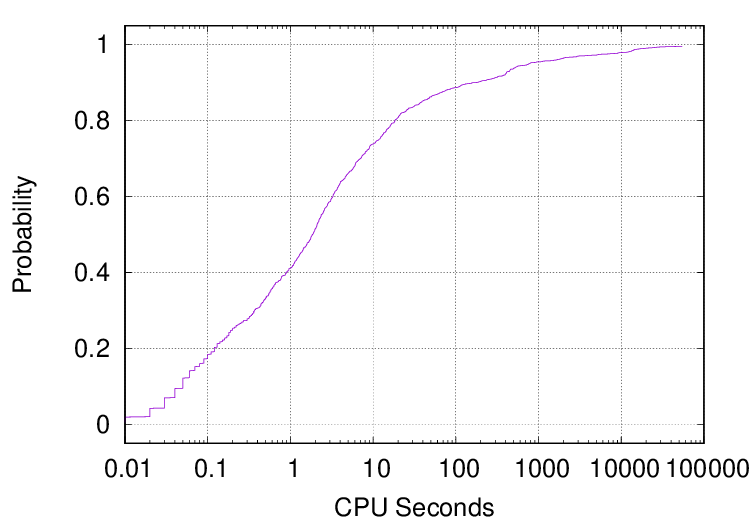}
 \captionof{figure}{CDF of Subtour Relaxation CPU Time: New Code}
 \label{fig:cdf-new-time}
\end{minipage}
\end{center}
\vspace*{\floatsep}
\begin{center}
\begin{minipage}[t]{2.25in}
\includegraphics[width=2.25in]{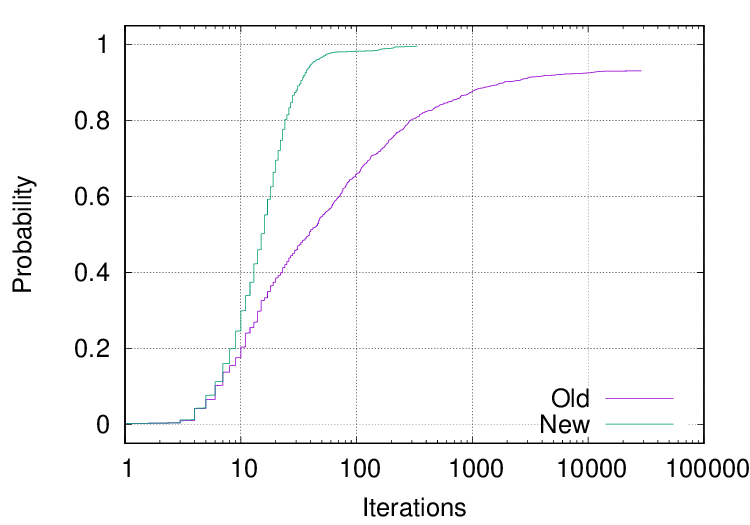}
 \captionof{figure}{CDF of Subtour Relaxation Iterations: Old vs New Code}
 \label{fig:cdf-both-iterations}
\end{minipage}
\quad
\begin{minipage}[t]{2.25in}
\includegraphics[width=2.25in]{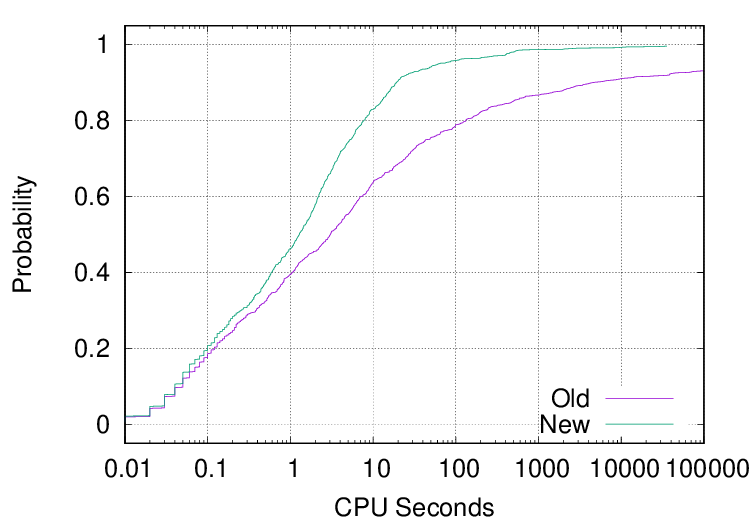}
 \captionof{figure}{CDF of Subtour Relaxation CPU Time: Old vs New Code}
 \label{fig:cdf-both-time}
\end{minipage}
\end{center}
\end{figure}

\begin{figure}[!ht]
\begin{center}
\includegraphics[width=4.75in]{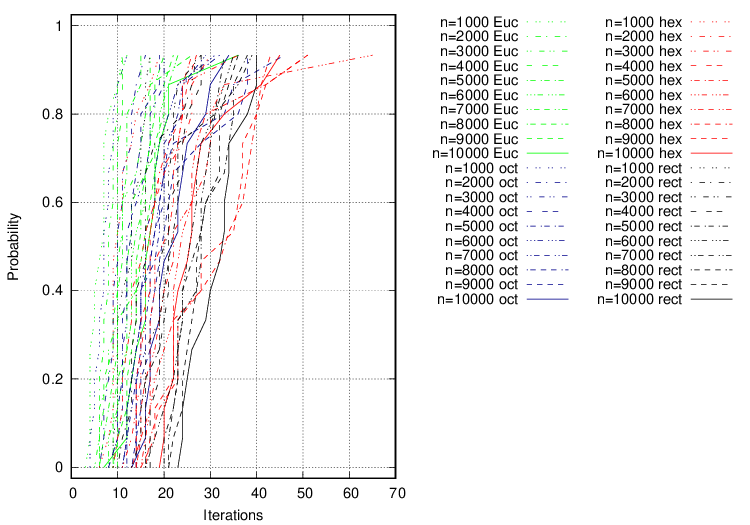}
 \captionof{figure}{CDFs of subtour relaxation iterations, random
   instances with $1000 \le n \le 10000$,
   all four metrics (new code only)}
 \label{fig:cdfs-random-all-metrics}
\end{center}
\end{figure}
\FloatBarrier

% Variables that control cell padding in these tables.
\newlength{\MyTabColSep}
\setlength{\MyTabColSep}{3.0pt}
\newcommand{\MyArrayStretch}{1.0}

{
% Must make this table narrower.
\setlength{\MyTabColSep}{2.0pt}
\MyBeginMeanTable
\caption{Mean FST Concatenation: Random Euclidean Instances}
{\scriptsize
\InitializeTableFootnoteFlags
\setlength{\tabcolsep}{\MyTabColSep}%
\renewcommand{\arraystretch}{\MyArrayStretch}%
\begin{tabular}{r|r|r|r|c|r|r|r|r|c|r|}
\cline{2-11}
 & \multicolumn{5}{c|}{Old} & \multicolumn{5}{c|}{New} \\
\cline{2-11}
 \multicolumn{1}{c|}{Instance}
 & & \multicolumn{1}{c|}{Relax}
 & & & \multicolumn{1}{c|}{Concat}
 & & \multicolumn{1}{c|}{Relax}
 & & & \multicolumn{1}{c|}{Concat} \\
\multicolumn{1}{c|}{Size} &
\multicolumn{1}{c|}{Iter} &
\multicolumn{1}{c|}{Time} &
\multicolumn{1}{c|}{Nodes} &
\multicolumn{1}{c|}{Done} &
\multicolumn{1}{c|}{Time} &
\multicolumn{1}{c|}{Iter} &
\multicolumn{1}{c|}{Time} &
\multicolumn{1}{c|}{Nodes} &
\multicolumn{1}{c|}{Done} &
\multicolumn{1}{c|}{Time} \\
\hline
1000 & 7.000 & .05133 & 1.000 & 15/15 & .05133 & 6.133 & .04400 & 1.000 & 15/15 & .04400 \\
 & \SDEV{2.683} & \SDEV{.01821} & \SDEV{0.0000} & & \SDEV{.01821} & \SDEV{2.247} & \SDEV{.01541} & \SDEV{0.0000} & & \SDEV{.01541} \\
 & & & & & & & & & & \\
2000 & 12.73 & .2320 & 1.000 & 15/15 & .2447 & 8.067 & .1533 & 1.000 & 15/15 & .1647 \\
 & \SDEV{12.76} & \SDEV{.1954} & \SDEV{0.0000} & & \SDEV{.1925} & \SDEV{1.731} & \SDEV{.03553} & \SDEV{0.0000} & & \SDEV{.04097} \\
\end{tabular}
}
\label{tab:mean-random-e}
\MyEndMeanTable
}

\MyBeginMeanTable
\caption{Mean FST Concatenation: Random Rectilinear Instances}
{\scriptsize
\InitializeTableFootnoteFlags
\setlength{\tabcolsep}{\MyTabColSep}%
\renewcommand{\arraystretch}{\MyArrayStretch}%
\begin{tabular}{r|r|r|r|c|r|r|r|r|c|r|}
\cline{2-11}
 & \multicolumn{5}{c|}{Old} & \multicolumn{5}{c|}{New} \\
\cline{2-11}
 \multicolumn{1}{c|}{Instance}
 & & \multicolumn{1}{c|}{Relax}
 & & & \multicolumn{1}{c|}{Concat}
 & & \multicolumn{1}{c|}{Relax}
 & & & \multicolumn{1}{c|}{Concat} \\
\multicolumn{1}{c|}{Size} &
\multicolumn{1}{c|}{Iter} &
\multicolumn{1}{c|}{Time} &
\multicolumn{1}{c|}{Nodes} &
\multicolumn{1}{c|}{Done} &
\multicolumn{1}{c|}{Time} &
\multicolumn{1}{c|}{Iter} &
\multicolumn{1}{c|}{Time} &
\multicolumn{1}{c|}{Nodes} &
\multicolumn{1}{c|}{Done} &
\multicolumn{1}{c|}{Time} \\
\hline
1000 & 21.80 & .3560 & 1.667 & 15/15 & 1467 & 12.47 & .2493 & 2.933 & 15/15 & .4553 \\
 & \SDEV{10.36} & \SDEV{.1676} & \SDEV{2.022} & & \SDEV{5487} & \SDEV{2.705} & \SDEV{.07307} & \SDEV{5.471} & & \SDEV{.2641} \\
 & & & & & & & & & & \\
2000 & 61.40 & 35.88 & 2.067 & 15/15 & 179.9 & 16.47 & .9413 & 3.000 & 15/15 & 2.790 \\
 & \SDEV{81.21} & \SDEV{126.6} & \SDEV{3.991} & & \SDEV{520.3} & \SDEV{5.045} & \SDEV{.3132} & \SDEV{3.327} & & \SDEV{3.775} \\
\end{tabular}
}
\label{tab:mean-random-r}
\MyEndMeanTable

\MyBeginMeanTable
\caption{Mean FST Concatenation: Random Hexagonal Instances}
{\scriptsize
\InitializeTableFootnoteFlags
\setlength{\tabcolsep}{\MyTabColSep}%
\renewcommand{\arraystretch}{\MyArrayStretch}%
\begin{tabular}{r|r|r|r|c|r|r|r|r|c|r|}
\cline{2-11}
 & \multicolumn{5}{c|}{Old} & \multicolumn{5}{c|}{New} \\
\cline{2-11}
 \multicolumn{1}{c|}{Instance}
 & & \multicolumn{1}{c|}{Relax}
 & & & \multicolumn{1}{c|}{Concat}
 & & \multicolumn{1}{c|}{Relax}
 & & & \multicolumn{1}{c|}{Concat} \\
\multicolumn{1}{c|}{Size} &
\multicolumn{1}{c|}{Iter} &
\multicolumn{1}{c|}{Time} &
\multicolumn{1}{c|}{Nodes} &
\multicolumn{1}{c|}{Done} &
\multicolumn{1}{c|}{Time} &
\multicolumn{1}{c|}{Iter} &
\multicolumn{1}{c|}{Time} &
\multicolumn{1}{c|}{Nodes} &
\multicolumn{1}{c|}{Done} &
\multicolumn{1}{c|}{Time} \\
\hline
1000 & 33.33 & .2293 & 1.000 & 15/15 & 10.77 & 10.87 & .1147 & 1.267 & 15/15 & .2333 \\
 & \SDEV{79.22} & \SDEV{.3861} & \SDEV{0.0000} & & \SDEV{30.03} & \SDEV{3.384} & \SDEV{.02553} & \SDEV{.9978} & & \SDEV{.2932} \\
 & & & & & & & & & & \\
2000 & 91.20 & 1.505 & 1.000 & 15/15 & 1.972 & 16.93 & .4373 & 1.133 & 15/15 & .6767 \\
 & \SDEV{133.4} & \SDEV{2.008} & \SDEV{0.0000} & & \SDEV{2.361} & \SDEV{6.516} & \SDEV{.1230} & \SDEV{.3399} & & \SDEV{.5454} \\
\end{tabular}
}
\label{tab:mean-random-u3}
\MyEndMeanTable

\MyBeginMeanTable
\caption{Mean FST Concatenation: Random Octilinear Instances}
{\scriptsize
\InitializeTableFootnoteFlags
\setlength{\tabcolsep}{\MyTabColSep}%
\renewcommand{\arraystretch}{\MyArrayStretch}%
\begin{tabular}{r|r|r|r|c|r|r|r|r|c|r|}
\cline{2-11}
 & \multicolumn{5}{c|}{Old} & \multicolumn{5}{c|}{New} \\
\cline{2-11}
 \multicolumn{1}{c|}{Instance}
 & & \multicolumn{1}{c|}{Relax}
 & & & \multicolumn{1}{c|}{Concat}
 & & \multicolumn{1}{c|}{Relax}
 & & & \multicolumn{1}{c|}{Concat} \\
\multicolumn{1}{c|}{Size} &
\multicolumn{1}{c|}{Iter} &
\multicolumn{1}{c|}{Time} &
\multicolumn{1}{c|}{Nodes} &
\multicolumn{1}{c|}{Done} &
\multicolumn{1}{c|}{Time} &
\multicolumn{1}{c|}{Iter} &
\multicolumn{1}{c|}{Time} &
\multicolumn{1}{c|}{Nodes} &
\multicolumn{1}{c|}{Done} &
\multicolumn{1}{c|}{Time} \\
\hline
1000 & 9.800 & .06400 & 1.000 & 15/15 & .1573 & 7.133 & .05200 & 1.000 & 15/15 & .05600 \\
 & \SDEV{5.924} & \SDEV{.02185} & \SDEV{0.0000} & & \SDEV{.3408} & \SDEV{2.187} & \SDEV{.01275} & \SDEV{0.0000} & & \SDEV{.01405} \\
 & & & & & & & & & & \\
2000 & 17.47 & .2473 & 1.000 & 15/15 & .2820 & 10.20 & .1907 & 1.000 & 15/15 & .2080 \\
 & \SDEV{13.03} & \SDEV{.1076} & \SDEV{0.0000} & & \SDEV{.1449} & \SDEV{2.737} & \SDEV{.06005} & \SDEV{0.0000} & & \SDEV{.07277} \\
\end{tabular}
}
\label{tab:mean-random-u4}
\MyEndMeanTable

\MyBeginSummaryTable
\caption{TSPLIB Subtour Relaxation Summary: Iterations and Time}
{\tiny
\InitializeTableFootnoteFlags
\setlength{\tabcolsep}{\MyTabColSep}%
\renewcommand{\arraystretch}{\MyArrayStretch}%
\begin{tabular}{r|r|r|r|r|r|r|r|r|r|r|r|r|r|r|r|r|}
\cline{2-17}
\multicolumn{1}{c}{} &
\multicolumn{4}{|c}{Euclidean} &
\multicolumn{4}{|c}{Rectilinear} &
\multicolumn{4}{|c}{Hexagonal} &
\multicolumn{4}{|c|}{Octilinear} \\
\cline{2-17}
\multicolumn{1}{c}{} &
\multicolumn{2}{|c}{Iter} &
\multicolumn{2}{|c}{Time} &
\multicolumn{2}{|c}{Iter} &
\multicolumn{2}{|c}{Time} &
\multicolumn{2}{|c}{Iter} &
\multicolumn{2}{|c}{Time} &
\multicolumn{2}{|c}{Iter} &
\multicolumn{2}{|c|}{Time} \\
\cline{2-17}
\multicolumn{1}{c}{Instance} &
\multicolumn{1}{|c}{Old} &
\multicolumn{1}{|c}{New} &
\multicolumn{1}{|c}{Old} &
\multicolumn{1}{|c}{New} &
\multicolumn{1}{|c}{Old} &
\multicolumn{1}{|c}{New} &
\multicolumn{1}{|c}{Old} &
\multicolumn{1}{|c}{New} &
\multicolumn{1}{|c}{Old} &
\multicolumn{1}{|c}{New} &
\multicolumn{1}{|c}{Old} &
\multicolumn{1}{|c}{New} &
\multicolumn{1}{|c}{Old} &
\multicolumn{1}{|c}{New} &
\multicolumn{1}{|c}{Old} &
\multicolumn{1}{|c|}{New} \\
\hline
d198 &
	8 & 5 & 0.01 & 0.01 &
	3 & 3 & 0.00 & 0.00 &
	24 & 7 & 0.03 & 0.02 &
	4 & 4 & 0.00 & 0.00 \\
lin318 &
	7 & 7 & 0.03 & 0.03 &
	7 & 6 & 0.03 & 0.04 &
	6 & 5 & 0.03 & 0.02 &
	83 & 10 & 0.23 & 0.05 \\
fl417 &
	6 & 5 & 0.08 & 0.09 &
	26 & 9 & 0.12 & 0.05 &
	98 & 6 & 0.52 & 0.08 &
	8 & 5 & 0.05 & 0.05 \\
pcb442 &
	15 & 10 & 0.13 & 0.11 &
	3 & 3 & 0.00 & 0.00 &
	6 & 6 & 0.04 & 0.04 &
	14 & 13 & 0.03 & 0.02 \\
att532 &
	19 & 13 & 0.06 & 0.04 &
	213 & 15 & 23.35 & 0.13 &
	19 & 19 & 0.09 & 0.11 &
	11 & 11 & 0.04 & 0.04 \\
\end{tabular}
}
\label{tab:summary-tsplib-relax}
\MyEndSummaryTable

{
% This table mush be made narrower.
\setlength{\MyTabColSep}{0.8pt}
\MyBeginSummaryTable
\caption{TSPLIB FST Concatenation Summary: Nodes and Time}
{\tiny
\InitializeTableFootnoteFlags
\setlength{\tabcolsep}{\MyTabColSep}%
\renewcommand{\arraystretch}{\MyArrayStretch}%
\begin{tabular}{r|r|r|r|r|r|r|r|r|r|r|r|r|r|r|r|r|}
\cline{2-17}
\multicolumn{1}{c}{} &
\multicolumn{4}{|c}{Euclidean} &
\multicolumn{4}{|c}{Rectilinear} &
\multicolumn{4}{|c}{Hexagonal} &
\multicolumn{4}{|c|}{Octilinear} \\
\cline{2-17}
\multicolumn{1}{c}{} &
\multicolumn{2}{|c}{Nodes} &
\multicolumn{2}{|c}{Time} &
\multicolumn{2}{|c}{Nodes} &
\multicolumn{2}{|c}{Time} &
\multicolumn{2}{|c}{Nodes} &
\multicolumn{2}{|c}{Time} &
\multicolumn{2}{|c}{Nodes} &
\multicolumn{2}{|c|}{Time} \\
\cline{2-17}
\multicolumn{1}{c}{Instance} &
\multicolumn{1}{|c}{Old} &
\multicolumn{1}{|c}{New} &
\multicolumn{1}{|c}{Old} &
\multicolumn{1}{|c}{New} &
\multicolumn{1}{|c}{Old} &
\multicolumn{1}{|c}{New} &
\multicolumn{1}{|c}{Old} &
\multicolumn{1}{|c}{New} &
\multicolumn{1}{|c}{Old} &
\multicolumn{1}{|c}{New} &
\multicolumn{1}{|c}{Old} &
\multicolumn{1}{|c}{New} &
\multicolumn{1}{|c}{Old} &
\multicolumn{1}{|c}{New} &
\multicolumn{1}{|c}{Old} &
\multicolumn{1}{|c|}{New} \\
\hline
d198 &
	1 & 1 & 0.01 & 0.01 &
	1 & 1 & 0.00 & 0.00 &
	1 & 1 & 0.69 & 0.32 &
	1 & 1 & 0.00 & 0.00 \\
lin318 &
	1 & 1 & 0.03 & 0.03 &
	1 & 3 & 211.50 & 4.53 &
	1 & 1 & 0.03 & 0.02 &
	1 & 1 & 0.23 & 0.05 \\
fl417 &
	1 & 1 & 0.08 & 0.09 &
	1 & 3 & 40.27 & 1.33 &
	1 & 1 & 0.52 & 0.08 &
	1 & 1 & 0.05 & 0.05 \\
pcb442 &
	19 & 37 & 1.22 & 0.82 &
	1 & 1 & 0.00 & 0.00 &
	12 & 23 & 0.32 & 0.28 &
	7 & 5 & 0.29 & 0.10 \\
att532 &
	1 & 1 & 0.06 & 0.04 &
	1 & 3 & 27.52 & 0.25 &
	1 & 1 & 0.09 & 0.11 &
	1 & 1 & 0.04 & 0.04 \\
\end{tabular}
}
\label{tab:summary-tsplib-concat}
\MyEndSummaryTable
}

{
% This table must be reduced in height.
\renewcommand{\MyArrayStretch}{0.85}
\MyBeginConcatTable
\caption{FST Concatenation: Euclidean OR-Library Instances}
{\tiny
\InitializeTableFootnoteFlags
\setlength{\tabcolsep}{\MyTabColSep}%
\renewcommand{\arraystretch}{\MyArrayStretch}%
\begin{tabular}{r|r|r|r|r|r|r|r|r|r|r|r|r|r|r|r|r|r}
\cline{2-17}
\multicolumn{1}{c}{Size} &
\multicolumn{6}{|c}{Subtour Relaxation} &
\multicolumn{7}{|c}{Root Node} &
\multicolumn{3}{|c}{Branching} &
\multicolumn{1}{|c}{Total} \\
\cline{2-17}
\multicolumn{1}{c}{(Instance)} &
\multicolumn{1}{|c}{Iter} &
\multicolumn{1}{|c}{NZ/Col} &
\multicolumn{1}{|c}{CC} &
\multicolumn{1}{|c}{Gap} &
\multicolumn{1}{|c}{Frac} &
\multicolumn{1}{|c}{Time} &
\multicolumn{1}{|c}{Iter} &
\multicolumn{1}{|c}{Rnds} &
\multicolumn{1}{|c}{Cuts} &
\multicolumn{1}{|c}{NZPC} &
\multicolumn{1}{|c}{Gap} &
\multicolumn{1}{|c}{Frac} &
\multicolumn{1}{|c}{Time} &
\multicolumn{1}{|c}{Nodes} &
\multicolumn{1}{|c}{Proc} &
\multicolumn{1}{|c}{Time} &
\multicolumn{1}{|c}{Time} \\
\hline
500 (1) & 7~ & 3.4745~ & 19~ & 0~ & 0~ & 0.01~ & 0~ & 0~ & 0~ & 0~ & 0~ & 0~ & 0.00~ & 1~ & 1~ & 0.00~ & 0.01~ \\
 & 6 & 3.65341 & 20 & 0 & 0 & 0.02 & 0 & 0 & 0 & 0 & 0 & 0 & 0.00 & 1 & 1 & 0.00 & 0.02 \\
500 (2) & 14~ & 3.58297~ & 46~ & 0~ & 0~ & 0.04~ & 0~ & 0~ & 0~ & 0~ & 0~ & 0~ & 0.00~ & 1~ & 1~ & 0.00~ & 0.04~ \\
 & 13 & 3.95091 & 39 & 0 & 0 & 0.04 & 0 & 0 & 0 & 0 & 0 & 0 & 0.00 & 1 & 1 & 0.00 & 0.04 \\
500 (3) & 4~ & 2.56138~ & 23~ & 0~ & 0~ & 0.01~ & 0~ & 0~ & 0~ & 0~ & 0~ & 0~ & 0.00~ & 1~ & 1~ & 0.00~ & 0.01~ \\
 & 4 & 2.9625 & 23 & 0 & 0 & 0.01 & 0 & 0 & 0 & 0 & 0 & 0 & 0.00 & 1 & 1 & 0.00 & 0.01 \\
\end{tabular}
}
\label{tab:concat-orlib-e}
\MyEndConcatTable
}

{
% This table must be reduced in height.
\renewcommand{\MyArrayStretch}{0.85}
\MyBeginConcatTable
\caption{FST Concatenation: Rectilinear OR-Library Instances}
{\tiny
\InitializeTableFootnoteFlags
\setlength{\tabcolsep}{\MyTabColSep}%
\renewcommand{\arraystretch}{\MyArrayStretch}%
\begin{tabular}{r|r|r|r|r|r|r|r|r|r|r|r|r|r|r|r|r|r}
\cline{2-17}
\multicolumn{1}{c}{Size} &
\multicolumn{6}{|c}{Subtour Relaxation} &
\multicolumn{7}{|c}{Root Node} &
\multicolumn{3}{|c}{Branching} &
\multicolumn{1}{|c}{Total} \\
\cline{2-17}
\multicolumn{1}{c}{(Instance)} &
\multicolumn{1}{|c}{Iter} &
\multicolumn{1}{|c}{NZ/Col} &
\multicolumn{1}{|c}{CC} &
\multicolumn{1}{|c}{Gap} &
\multicolumn{1}{|c}{Frac} &
\multicolumn{1}{|c}{Time} &
\multicolumn{1}{|c}{Iter} &
\multicolumn{1}{|c}{Rnds} &
\multicolumn{1}{|c}{Cuts} &
\multicolumn{1}{|c}{NZPC} &
\multicolumn{1}{|c}{Gap} &
\multicolumn{1}{|c}{Frac} &
\multicolumn{1}{|c}{Time} &
\multicolumn{1}{|c}{Nodes} &
\multicolumn{1}{|c}{Proc} &
\multicolumn{1}{|c}{Time} &
\multicolumn{1}{|c}{Time} \\
\hline
500 (1) & 14~ & 6.76257~ & 39~ & 37.633~ & 13~ & 0.08~ & 3~ & 3~ & 3~ & 11.3333~ & 0~ & 8~ & 0.02~ & 1~ & 1~ & 0.01~ & 0.11~ \\
 & 7 & 5.03916 & 21 & 1954.966 & 56 & 0.05 & 18 & 13 & 16 & 12.3125 & 0 & 36 & 0.75 & 1 & 1 & 0.00 & 0.80 \\
500 (2) & 12~ & 6.36133~ & 36~ & 86.646~ & 19~ & 0.08~ & 2~ & 2~ & 2~ & 6~ & 0~ & 8~ & 0.02~ & 1~ & 1~ & 0.00~ & 0.10~ \\
 & 11 & 7.47525 & 70 & 86.646 & 23 & 0.10 & 1 & 1 & 1 & 5 & 0 & 12 & 0.00 & 1 & 1 & 0.00 & 0.10 \\
500 (3) & 46~ & 10.6206~ & 91~ & 0~ & 32~ & 0.22~ & 0~ & 0~ & 0~ & 0~ & 0~ & 32~ & 0.00~ & 1~ & 1~ & 0.00~ & 0.22~ \\
 & 13 & 7.75492 & 48 & 0 & 0 & 0.07 & 0 & 0 & 0 & 0 & 0 & 0 & 0.00 & 1 & 1 & 0.00 & 0.07 \\
\end{tabular}
}
\label{tab:concat-orlib-r}
\MyEndConcatTable
}

{
% This table must be reduced in height.
\renewcommand{\MyArrayStretch}{0.85}
\MyBeginConcatTable
\caption{FST Concatenation: Hexagonal OR-Library Instances}
{\tiny
\InitializeTableFootnoteFlags
\setlength{\tabcolsep}{\MyTabColSep}%
\renewcommand{\arraystretch}{\MyArrayStretch}%
\begin{tabular}{r|r|r|r|r|r|r|r|r|r|r|r|r|r|r|r|r|r}
\cline{2-17}
\multicolumn{1}{c}{Size} &
\multicolumn{6}{|c}{Subtour Relaxation} &
\multicolumn{7}{|c}{Root Node} &
\multicolumn{3}{|c}{Branching} &
\multicolumn{1}{|c}{Total} \\
\cline{2-17}
\multicolumn{1}{c}{(Instance)} &
\multicolumn{1}{|c}{Iter} &
\multicolumn{1}{|c}{NZ/Col} &
\multicolumn{1}{|c}{CC} &
\multicolumn{1}{|c}{Gap} &
\multicolumn{1}{|c}{Frac} &
\multicolumn{1}{|c}{Time} &
\multicolumn{1}{|c}{Iter} &
\multicolumn{1}{|c}{Rnds} &
\multicolumn{1}{|c}{Cuts} &
\multicolumn{1}{|c}{NZPC} &
\multicolumn{1}{|c}{Gap} &
\multicolumn{1}{|c}{Frac} &
\multicolumn{1}{|c}{Time} &
\multicolumn{1}{|c}{Nodes} &
\multicolumn{1}{|c}{Proc} &
\multicolumn{1}{|c}{Time} &
\multicolumn{1}{|c}{Time} \\
\hline
500 (1) & 4~ & 5.004~ & 19~ & 0~ & 0~ & 0.02~ & 0~ & 0~ & 0~ & 0~ & 0~ & 0~ & 0.00~ & 1~ & 1~ & 0.00~ & 0.02~ \\
 & 4 & 4.94939 & 19 & 0 & 0 & 0.03 & 0 & 0 & 0 & 0 & 0 & 0 & 0.00 & 1 & 1 & 0.00 & 0.03 \\
500 (2) & 6~ & 7.55532~ & 29~ & 0~ & 0~ & 0.04~ & 0~ & 0~ & 0~ & 0~ & 0~ & 0~ & 0.00~ & 1~ & 1~ & 0.00~ & 0.04~ \\
 & 5 & 5.87797 & 29 & 0 & 0 & 0.04 & 0 & 0 & 0 & 0 & 0 & 0 & 0.00 & 1 & 1 & 0.00 & 0.04 \\
500 (3) & 5~ & 5.75853~ & 17~ & 0~ & 0~ & 0.03~ & 0~ & 0~ & 0~ & 0~ & 0~ & 0~ & 0.00~ & 1~ & 1~ & 0.00~ & 0.03~ \\
 & 5 & 5.5714 & 27 & 0 & 0 & 0.03 & 0 & 0 & 0 & 0 & 0 & 0 & 0.00 & 1 & 1 & 0.00 & 0.03 \\
\end{tabular}
}
\label{tab:concat-orlib-u3}
\MyEndConcatTable
}

{
% This table must be reduced in height.
\renewcommand{\MyArrayStretch}{0.85}
\MyBeginConcatTable
\caption{FST Concatenation: Octilinear OR-Library Instances}
{\tiny
\InitializeTableFootnoteFlags
\setlength{\tabcolsep}{\MyTabColSep}%
\renewcommand{\arraystretch}{\MyArrayStretch}%
\begin{tabular}{r|r|r|r|r|r|r|r|r|r|r|r|r|r|r|r|r|r}
\cline{2-17}
\multicolumn{1}{c}{Size} &
\multicolumn{6}{|c}{Subtour Relaxation} &
\multicolumn{7}{|c}{Root Node} &
\multicolumn{3}{|c}{Branching} &
\multicolumn{1}{|c}{Total} \\
\cline{2-17}
\multicolumn{1}{c}{(Instance)} &
\multicolumn{1}{|c}{Iter} &
\multicolumn{1}{|c}{NZ/Col} &
\multicolumn{1}{|c}{CC} &
\multicolumn{1}{|c}{Gap} &
\multicolumn{1}{|c}{Frac} &
\multicolumn{1}{|c}{Time} &
\multicolumn{1}{|c}{Iter} &
\multicolumn{1}{|c}{Rnds} &
\multicolumn{1}{|c}{Cuts} &
\multicolumn{1}{|c}{NZPC} &
\multicolumn{1}{|c}{Gap} &
\multicolumn{1}{|c}{Frac} &
\multicolumn{1}{|c}{Time} &
\multicolumn{1}{|c}{Nodes} &
\multicolumn{1}{|c}{Proc} &
\multicolumn{1}{|c}{Time} &
\multicolumn{1}{|c}{Time} \\
\hline
500 (1) & 30~ & 5.774~ & 43~ & 0~ & 0~ & 0.09~ & 0~ & 0~ & 0~ & 0~ & 0~ & 0~ & 0.00~ & 1~ & 1~ & 0.00~ & 0.09~ \\
 & 4 & 4.35728 & 21 & 0 & 0 & 0.02 & 0 & 0 & 0 & 0 & 0 & 0 & 0.00 & 1 & 1 & 0.00 & 0.02 \\
500 (2) & 8~ & 3.78013~ & 30~ & 23.705~ & 4~ & 0.02~ & 5~ & 4~ & 4~ & 6~ & 0~ & 0~ & 0.01~ & 1~ & 1~ & 0.00~ & 0.03~ \\
 & 8 & 4.53921 & 30 & 23.705 & 4 & 0.02 & 5 & 4 & 4 & 6 & 0 & 0 & 0.01 & 1 & 1 & 0.00 & 0.03 \\
500 (3) & 6~ & 4.22269~ & 36~ & 0~ & 0~ & 0.02~ & 0~ & 0~ & 0~ & 0~ & 0~ & 0~ & 0.00~ & 1~ & 1~ & 0.00~ & 0.02~ \\
 & 6 & 4.58519 & 36 & 0 & 0 & 0.02 & 0 & 0 & 0 & 0 & 0 & 0 & 0.00 & 1 & 1 & 0.00 & 0.02 \\
\end{tabular}
}
\label{tab:concat-orlib-u4}
\MyEndConcatTable
}

\FloatBarrier

{
% This table must be reduced in height.
\renewcommand{\MyArrayStretch}{0.85}
\MyBeginConcatTable
\caption{FST Concatenation: Euclidean Random Instances}
{\tiny
\InitializeTableFootnoteFlags
\setlength{\tabcolsep}{\MyTabColSep}%
\renewcommand{\arraystretch}{\MyArrayStretch}%
\begin{tabular}{r|r|r|r|r|r|r|r|r|r|r|r|r|r|r|r|r|r}
\cline{2-17}
\multicolumn{1}{c}{Size} &
\multicolumn{6}{|c}{Subtour Relaxation} &
\multicolumn{7}{|c}{Root Node} &
\multicolumn{3}{|c}{Branching} &
\multicolumn{1}{|c}{Total} \\
\cline{2-17}
\multicolumn{1}{c}{(Instance)} &
\multicolumn{1}{|c}{Iter} &
\multicolumn{1}{|c}{NZ/Col} &
\multicolumn{1}{|c}{CC} &
\multicolumn{1}{|c}{Gap} &
\multicolumn{1}{|c}{Frac} &
\multicolumn{1}{|c}{Time} &
\multicolumn{1}{|c}{Iter} &
\multicolumn{1}{|c}{Rnds} &
\multicolumn{1}{|c}{Cuts} &
\multicolumn{1}{|c}{NZPC} &
\multicolumn{1}{|c}{Gap} &
\multicolumn{1}{|c}{Frac} &
\multicolumn{1}{|c}{Time} &
\multicolumn{1}{|c}{Nodes} &
\multicolumn{1}{|c}{Proc} &
\multicolumn{1}{|c}{Time} &
\multicolumn{1}{|c}{Time} \\
\hline
1000 (1) & 8~ & 2.79704~ & 15~ & 0~ & 0~ & 0.06~ & 0~ & 0~ & 0~ & 0~ & 0~ & 0~ & 0.00~ & 1~ & 1~ & 0.00~ & 0.06~ \\
 & 8 & 4.0401 & 15 & 0 & 0 & 0.06 & 0 & 0 & 0 & 0 & 0 & 0 & 0.00 & 1 & 1 & 0.00 & 0.06 \\
1000 (2) & 6~ & 2.30946~ & 35~ & 0~ & 0~ & 0.03~ & 0~ & 0~ & 0~ & 0~ & 0~ & 0~ & 0.00~ & 1~ & 1~ & 0.00~ & 0.03~ \\
 & 4 & 2.88729 & 35 & 0 & 0 & 0.03 & 0 & 0 & 0 & 0 & 0 & 0 & 0.00 & 1 & 1 & 0.00 & 0.03 \\
1000 (3) & 6~ & 2.60257~ & 29~ & 0~ & 0~ & 0.05~ & 0~ & 0~ & 0~ & 0~ & 0~ & 0~ & 0.00~ & 1~ & 1~ & 0.00~ & 0.05~ \\
 & 7 & 4.11681 & 42 & 0 & 0 & 0.05 & 0 & 0 & 0 & 0 & 0 & 0 & 0.00 & 1 & 1 & 0.00 & 0.05 \\
\end{tabular}
}
\label{tab:concat-rp-1-e}
\MyEndConcatTable
}

{
% This table must be reduced in height.
\renewcommand{\MyArrayStretch}{0.85}
\MyBeginConcatTable
\caption{FST Concatenation: Euclidean Random Instances (continued)}
{\tiny
\InitializeTableFootnoteFlags
\setlength{\tabcolsep}{\MyTabColSep}%
\renewcommand{\arraystretch}{\MyArrayStretch}%
\begin{tabular}{r|r|r|r|r|r|r|r|r|r|r|r|r|r|r|r|r|r}
\cline{2-17}
\multicolumn{1}{c}{Size} &
\multicolumn{6}{|c}{Subtour Relaxation} &
\multicolumn{7}{|c}{Root Node} &
\multicolumn{3}{|c}{Branching} &
\multicolumn{1}{|c}{Total} \\
\cline{2-17}
\multicolumn{1}{c}{(Instance)} &
\multicolumn{1}{|c}{Iter} &
\multicolumn{1}{|c}{NZ/Col} &
\multicolumn{1}{|c}{CC} &
\multicolumn{1}{|c}{Gap} &
\multicolumn{1}{|c}{Frac} &
\multicolumn{1}{|c}{Time} &
\multicolumn{1}{|c}{Iter} &
\multicolumn{1}{|c}{Rnds} &
\multicolumn{1}{|c}{Cuts} &
\multicolumn{1}{|c}{NZPC} &
\multicolumn{1}{|c}{Gap} &
\multicolumn{1}{|c}{Frac} &
\multicolumn{1}{|c}{Time} &
\multicolumn{1}{|c}{Nodes} &
\multicolumn{1}{|c}{Proc} &
\multicolumn{1}{|c}{Time} &
\multicolumn{1}{|c}{Time} \\
\hline
3000 (1) & 7~ & 2.96838~ & 82~ & .632~ & 6~ & 0.21~ & 2~ & 2~ & 2~ & 5.5~ & 0~ & 0~ & 0.09~ & 1~ & 1~ & 0.00~ & 0.30~ \\
 & 8 & 3.60458 & 82 & .632 & 6 & 0.25 & 2 & 2 & 2 & 5.5 & 0 & 0 & 0.07 & 1 & 1 & 0.00 & 0.32 \\
3000 (2) & 12~ & 3.43526~ & 41~ & 0~ & 0~ & 0.39~ & 0~ & 0~ & 0~ & 0~ & 0~ & 0~ & 0.00~ & 1~ & 1~ & 0.00~ & 0.39~ \\
 & 11 & 4.23872 & 31 & 0 & 0 & 0.37 & 0 & 0 & 0 & 0 & 0 & 0 & 0.00 & 1 & 1 & 0.00 & 0.37 \\
3000 (3) & 9~ & 3.52081~ & 43~ & 11.032~ & 7~ & 0.30~ & 2~ & 2~ & 2~ & 5.5~ & 0~ & 0~ & 0.09~ & 1~ & 1~ & 0.00~ & 0.39~ \\
 & 10 & 4.73702 & 43 & 11.032 & 7 & 0.32 & 2 & 2 & 2 & 5.5 & 0 & 0 & 0.08 & 1 & 1 & 0.00 & 0.40 \\
\end{tabular}
}
\label{tab:concat-rp-2-e}
\MyEndConcatTable
}

{
% This table must be reduced in height.
\renewcommand{\MyArrayStretch}{0.85}
\MyBeginConcatTable
\caption{FST Concatenation: Euclidean Random Instances (continued)}
{\tiny
\InitializeTableFootnoteFlags
\setlength{\tabcolsep}{\MyTabColSep}%
\renewcommand{\arraystretch}{\MyArrayStretch}%
\begin{tabular}{r|r|r|r|r|r|r|r|r|r|r|r|r|r|r|r|r|r}
\cline{2-17}
\multicolumn{1}{c}{Size} &
\multicolumn{6}{|c}{Subtour Relaxation} &
\multicolumn{7}{|c}{Root Node} &
\multicolumn{3}{|c}{Branching} &
\multicolumn{1}{|c}{Total} \\
\cline{2-17}
\multicolumn{1}{c}{(Instance)} &
\multicolumn{1}{|c}{Iter} &
\multicolumn{1}{|c}{NZ/Col} &
\multicolumn{1}{|c}{CC} &
\multicolumn{1}{|c}{Gap} &
\multicolumn{1}{|c}{Frac} &
\multicolumn{1}{|c}{Time} &
\multicolumn{1}{|c}{Iter} &
\multicolumn{1}{|c}{Rnds} &
\multicolumn{1}{|c}{Cuts} &
\multicolumn{1}{|c}{NZPC} &
\multicolumn{1}{|c}{Gap} &
\multicolumn{1}{|c}{Frac} &
\multicolumn{1}{|c}{Time} &
\multicolumn{1}{|c}{Nodes} &
\multicolumn{1}{|c}{Proc} &
\multicolumn{1}{|c}{Time} &
\multicolumn{1}{|c}{Time} \\
\hline
5000 (1) & 15~ & 3.7235~ & 48~ & 0~ & 0~ & 0.92~ & 0~ & 0~ & 0~ & 0~ & 0~ & 0~ & 0.00~ & 1~ & 1~ & 0.00~ & 0.92~ \\
 & 12 & 5.09879 & 48 & 0 & 22 & 0.79 & 0 & 0 & 0 & 0 & 0 & 22 & 0.00 & 1 & 1 & 0.00 & 0.79 \\
5000 (2) & 15~ & 3.31213~ & 71~ & 0~ & 0~ & 0.83~ & 0~ & 0~ & 0~ & 0~ & 0~ & 0~ & 0.00~ & 1~ & 1~ & 0.00~ & 0.83~ \\
 & 10 & 4.33203 & 104 & 0 & 0 & 0.55 & 0 & 0 & 0 & 0 & 0 & 0 & 0.00 & 1 & 1 & 0.00 & 0.55 \\
5000 (3) & 7~ & 3.09815~ & 25~ & 0~ & 0~ & 0.54~ & 0~ & 0~ & 0~ & 0~ & 0~ & 0~ & 0.00~ & 1~ & 1~ & 0.00~ & 0.54~ \\
 & 7 & 4.50515 & 33 & 0 & 0 & 0.54 & 0 & 0 & 0 & 0 & 0 & 0 & 0.00 & 1 & 1 & 0.00 & 0.54 \\
\end{tabular}
}
\label{tab:concat-rp-3-e}
\MyEndConcatTable
}

{
% This table must be reduced in height.
\renewcommand{\MyArrayStretch}{0.85}
\MyBeginConcatTable
\caption{FST Concatenation: Euclidean Random Instances (continued)}
{\tiny
\InitializeTableFootnoteFlags
\setlength{\tabcolsep}{\MyTabColSep}%
\renewcommand{\arraystretch}{\MyArrayStretch}%
\begin{tabular}{r|r|r|r|r|r|r|r|r|r|r|r|r|r|r|r|r|r}
\cline{2-17}
\multicolumn{1}{c}{Size} &
\multicolumn{6}{|c}{Subtour Relaxation} &
\multicolumn{7}{|c}{Root Node} &
\multicolumn{3}{|c}{Branching} &
\multicolumn{1}{|c}{Total} \\
\cline{2-17}
\multicolumn{1}{c}{(Instance)} &
\multicolumn{1}{|c}{Iter} &
\multicolumn{1}{|c}{NZ/Col} &
\multicolumn{1}{|c}{CC} &
\multicolumn{1}{|c}{Gap} &
\multicolumn{1}{|c}{Frac} &
\multicolumn{1}{|c}{Time} &
\multicolumn{1}{|c}{Iter} &
\multicolumn{1}{|c}{Rnds} &
\multicolumn{1}{|c}{Cuts} &
\multicolumn{1}{|c}{NZPC} &
\multicolumn{1}{|c}{Gap} &
\multicolumn{1}{|c}{Frac} &
\multicolumn{1}{|c}{Time} &
\multicolumn{1}{|c}{Nodes} &
\multicolumn{1}{|c}{Proc} &
\multicolumn{1}{|c}{Time} &
\multicolumn{1}{|c}{Time} \\
\hline
7000 (1) & 28~ & 5.28259~ & 150~ & .088~ & 5~ & 3.12~ & 1~ & 1~ & 1~ & 5~ & 0~ & 0~ & 0.14~ & 1~ & 1~ & 0.00~ & 3.26~ \\
 & 12 & 6.30334 & 150 & .088 & 5 & 1.34 & 1 & 1 & 1 & 5 & 0 & 0 & 0.15 & 1 & 1 & 0.00 & 1.49 \\
7000 (2) & 40~ & 5.67071~ & 218~ & 2.913~ & 17~ & 3.86~ & 4~ & 3~ & 4~ & 5.75~ & 0~ & 3~ & 0.48~ & 1~ & 1~ & 0.12~ & 4.46~ \\
 & 13 & 5.76922 & 108 & 2.913 & 17 & 1.45 & 3 & 3 & 4 & 4.5 & 0 & 3 & 0.39 & 1 & 1 & 0.12 & 1.96 \\
7000 (3) & 15~ & 4.06466~ & 148~ & .474~ & 7~ & 1.51~ & 1~ & 1~ & 1~ & 3~ & 0~ & 0~ & 0.14~ & 1~ & 1~ & 0.00~ & 1.65~ \\
 & 18 & 6.90952 & 101 & .474 & 7 & 2.24 & 1 & 1 & 1 & 3 & 0 & 0 & 0.13 & 1 & 1 & 0.00 & 2.37 \\
\end{tabular}
}
\label{tab:concat-rp-4-e}
\MyEndConcatTable
}

{
% This table must be reduced in height.
\renewcommand{\MyArrayStretch}{0.85}
\MyBeginConcatTable
\caption{FST Concatenation: Euclidean Random Instances (continued)}
{\tiny
\InitializeTableFootnoteFlags
\setlength{\tabcolsep}{\MyTabColSep}%
\renewcommand{\arraystretch}{\MyArrayStretch}%
\begin{tabular}{r|r|r|r|r|r|r|r|r|r|r|r|r|r|r|r|r|r}
\cline{2-17}
\multicolumn{1}{c}{Size} &
\multicolumn{6}{|c}{Subtour Relaxation} &
\multicolumn{7}{|c}{Root Node} &
\multicolumn{3}{|c}{Branching} &
\multicolumn{1}{|c}{Total} \\
\cline{2-17}
\multicolumn{1}{c}{(Instance)} &
\multicolumn{1}{|c}{Iter} &
\multicolumn{1}{|c}{NZ/Col} &
\multicolumn{1}{|c}{CC} &
\multicolumn{1}{|c}{Gap} &
\multicolumn{1}{|c}{Frac} &
\multicolumn{1}{|c}{Time} &
\multicolumn{1}{|c}{Iter} &
\multicolumn{1}{|c}{Rnds} &
\multicolumn{1}{|c}{Cuts} &
\multicolumn{1}{|c}{NZPC} &
\multicolumn{1}{|c}{Gap} &
\multicolumn{1}{|c}{Frac} &
\multicolumn{1}{|c}{Time} &
\multicolumn{1}{|c}{Nodes} &
\multicolumn{1}{|c}{Proc} &
\multicolumn{1}{|c}{Time} &
\multicolumn{1}{|c}{Time} \\
\hline
9000 (1) & 57~ & 7.57055~ & 227~ & 0~ & 49~ & 7.25~ & 0~ & 0~ & 0~ & 0~ & 0~ & 49~ & 0.00~ & 1~ & 1~ & 0.00~ & 7.25~ \\
 & 18 & 7.62928 & 288 & 0 & 0 & 3.05 & 0 & 0 & 0 & 0 & 0 & 0 & 0.00 & 1 & 1 & 0.00 & 3.05 \\
9000 (2) & 14~ & 3.78585~ & 74~ & 4.426~ & 9~ & 2.18~ & 2~ & 1~ & 2~ & 6.5~ & 0~ & 0~ & 0.27~ & 1~ & 1~ & 0.00~ & 2.45~ \\
 & 10 & 4.50626 & 83 & 4.426 & 9 & 1.59 & 1 & 1 & 2 & 4.5 & 0 & 0 & 0.19 & 1 & 1 & 0.00 & 1.78 \\
9000 (3) & 20~ & 4.53387~ & 102~ & .691~ & 7~ & 2.84~ & 2~ & 2~ & 2~ & 5.5~ & 0~ & 0~ & 0.34~ & 1~ & 1~ & 0.00~ & 3.18~ \\
 & 10 & 4.93361 & 125 & .691 & 7 & 1.71 & 2 & 2 & 2 & 5.5 & 0 & 0 & 0.32 & 1 & 1 & 0.00 & 2.03 \\
\end{tabular}
}
\label{tab:concat-rp-5-e}
\MyEndConcatTable
}

\FloatBarrier

{
% This table must be reduced in height.
\renewcommand{\MyArrayStretch}{0.85}
\MyBeginConcatTable
\caption{FST Concatenation: Rectilinear Random Instances}
{\tiny
\InitializeTableFootnoteFlags
\setlength{\tabcolsep}{\MyTabColSep}%
\renewcommand{\arraystretch}{\MyArrayStretch}%
\begin{tabular}{r|r|r|r|r|r|r|r|r|r|r|r|r|r|r|r|r|r}
\cline{2-17}
\multicolumn{1}{c}{Size} &
\multicolumn{6}{|c}{Subtour Relaxation} &
\multicolumn{7}{|c}{Root Node} &
\multicolumn{3}{|c}{Branching} &
\multicolumn{1}{|c}{Total} \\
\cline{2-17}
\multicolumn{1}{c}{(Instance)} &
\multicolumn{1}{|c}{Iter} &
\multicolumn{1}{|c}{NZ/Col} &
\multicolumn{1}{|c}{CC} &
\multicolumn{1}{|c}{Gap} &
\multicolumn{1}{|c}{Frac} &
\multicolumn{1}{|c}{Time} &
\multicolumn{1}{|c}{Iter} &
\multicolumn{1}{|c}{Rnds} &
\multicolumn{1}{|c}{Cuts} &
\multicolumn{1}{|c}{NZPC} &
\multicolumn{1}{|c}{Gap} &
\multicolumn{1}{|c}{Frac} &
\multicolumn{1}{|c}{Time} &
\multicolumn{1}{|c}{Nodes} &
\multicolumn{1}{|c}{Proc} &
\multicolumn{1}{|c}{Time} &
\multicolumn{1}{|c}{Time} \\
\hline
1000 (1) & 17~ & 6.64013~ & 51~ & 60.981~ & 7~ & 0.24~ & 2~ & 2~ & 2~ & 12.5~ & 0~ & 11~ & 0.03~ & 1~ & 1~ & 0.00~ & 0.27~ \\
 & 16 & 8.13145 & 135 & 60.981 & 7 & 0.30 & 3 & 2 & 2 & 8.5 & 0 & 11 & 0.03 & 1 & 1 & 0.00 & 0.33 \\
1000 (2) & 29~ & 8.18195~ & 165~ & 14.967~ & 17~ & 0.54~ & 2~ & 2~ & 2~ & 4.5~ & 0~ & 12~ & 0.04~ & 1~ & 1~ & 0.00~ & 0.58~ \\
 & 13 & 7.19992 & 118 & 14.967 & 21 & 0.20 & 2 & 2 & 2 & 4.5 & 0 & 16 & 0.03 & 1 & 1 & 0.01 & 0.24 \\
1000 (3) & 11~ & 5.96086~ & 78~ & 115.818~ & 35~ & 0.20~ & 33~ & 5~ & 7~ & 13.7143~ & 0~ & 9~ & 0.40~ & 1~ & 1~ & 0.01~ & 0.61~ \\
 & 11 & 6.87728 & 96 & 115.818 & 35 & 0.20 & 15 & 5 & 7 & 10.5714 & 0 & 9 & 0.21 & 1 & 1 & 0.00 & 0.41 \\
\end{tabular}
}
\label{tab:concat-rp-1-r}
\MyEndConcatTable
}

{
% This table must be reduced in height.
\renewcommand{\MyArrayStretch}{0.85}
\MyBeginConcatTable
\caption{FST Concatenation: Rectilinear Random Instances (continued)}
{\tiny
\InitializeTableFootnoteFlags
\setlength{\tabcolsep}{\MyTabColSep}%
\renewcommand{\arraystretch}{\MyArrayStretch}%
\begin{tabular}{r|r|r|r|r|r|r|r|r|r|r|r|r|r|r|r|r|r}
\cline{2-17}
\multicolumn{1}{c}{Size} &
\multicolumn{6}{|c}{Subtour Relaxation} &
\multicolumn{7}{|c}{Root Node} &
\multicolumn{3}{|c}{Branching} &
\multicolumn{1}{|c}{Total} \\
\cline{2-17}
\multicolumn{1}{c}{(Instance)} &
\multicolumn{1}{|c}{Iter} &
\multicolumn{1}{|c}{NZ/Col} &
\multicolumn{1}{|c}{CC} &
\multicolumn{1}{|c}{Gap} &
\multicolumn{1}{|c}{Frac} &
\multicolumn{1}{|c}{Time} &
\multicolumn{1}{|c}{Iter} &
\multicolumn{1}{|c}{Rnds} &
\multicolumn{1}{|c}{Cuts} &
\multicolumn{1}{|c}{NZPC} &
\multicolumn{1}{|c}{Gap} &
\multicolumn{1}{|c}{Frac} &
\multicolumn{1}{|c}{Time} &
\multicolumn{1}{|c}{Nodes} &
\multicolumn{1}{|c}{Proc} &
\multicolumn{1}{|c}{Time} &
\multicolumn{1}{|c}{Time} \\
\hline
3000 (1) & 22~ & 8.29299~ & 205~ & 43.671~ & 73~ & 1.64~ & 6~ & 2~ & 7~ & 5.85714~ & 0~ & 42~ & 0.21~ & 1~ & 1~ & 0.10~ & 1.95~ \\
 & 16 & 9.10845 & 276 & 43.671 & 56 & 1.52 & 4 & 2 & 7 & 5.14286 & 0 & 25 & 0.25 & 1 & 1 & 0.00 & 1.77 \\
3000 (2) & 317~ & 61.3332~ & 941~ & 22.436~ & 57~ & 483.06~ & 32~ & 7~ & 17~ & 6.88235~ & 6.999~ & 33~ & 61.13~ & 1~ & 1~ & 7.37~ & 551.56~ \\
 & 22 & 11.4814 & 197 & 22.436 & 40 & 1.97 & 15 & 13 & 30 & 17.3333 & 0 & 16 & 4.47 & 1 & 1 & 0.00 & 6.44 \\
3000 (3) & 82~ & 14.3099~ & 749~ & 43.631~ & 86~ & 40.68~ & 12~ & 7~ & 13~ & 6.15385~ & 0~ & 47~ & 5.12~ & 1~ & 1~ & 0.67~ & 46.47~ \\
 & 23 & 12.3976 & 324 & 43.631 & 72 & 2.79 & 10 & 7 & 13 & 5.92308 & 0 & 33 & 0.91 & 1 & 1 & 0.00 & 3.70 \\
\end{tabular}
}
\label{tab:concat-rp-2-r}
\MyEndConcatTable
}

{
% This table must be reduced in height.
\renewcommand{\MyArrayStretch}{0.85}
\MyBeginConcatTable
\caption{FST Concatenation: Rectilinear Random Instances (continued)}
{\tiny
\InitializeTableFootnoteFlags
\setlength{\tabcolsep}{\MyTabColSep}%
\renewcommand{\arraystretch}{\MyArrayStretch}%
\begin{tabular}{r|r|r|r|r|r|r|r|r|r|r|r|r|r|r|r|r|r}
\cline{2-17}
\multicolumn{1}{c}{Size} &
\multicolumn{6}{|c}{Subtour Relaxation} &
\multicolumn{7}{|c}{Root Node} &
\multicolumn{3}{|c}{Branching} &
\multicolumn{1}{|c}{Total} \\
\cline{2-17}
\multicolumn{1}{c}{(Instance)} &
\multicolumn{1}{|c}{Iter} &
\multicolumn{1}{|c}{NZ/Col} &
\multicolumn{1}{|c}{CC} &
\multicolumn{1}{|c}{Gap} &
\multicolumn{1}{|c}{Frac} &
\multicolumn{1}{|c}{Time} &
\multicolumn{1}{|c}{Iter} &
\multicolumn{1}{|c}{Rnds} &
\multicolumn{1}{|c}{Cuts} &
\multicolumn{1}{|c}{NZPC} &
\multicolumn{1}{|c}{Gap} &
\multicolumn{1}{|c}{Frac} &
\multicolumn{1}{|c}{Time} &
\multicolumn{1}{|c}{Nodes} &
\multicolumn{1}{|c}{Proc} &
\multicolumn{1}{|c}{Time} &
\multicolumn{1}{|c}{Time} \\
\hline
5000 (1) & 70~ & 14.9799~ & 637~ & 62.242~ & 161~ & 21.25~ & 39~ & 7~ & 36~ & 8.61111~ & 0~ & 96~ & 23.93~ & 1~ & 1~ & 5.31~ & 50.49~ \\
 & 23 & 11.075 & 229 & 62.242 & 148 & 5.52 & 12 & 7 & 36 & 7.36111 & 0 & 83 & 1.50 & 1 & 1 & 0.00 & 7.02 \\
5000 (2) & 155~ & 19.988~ & 1090~ & 57.148~ & 118~ & 157.17~ & 90~ & 14~ & 46~ & 18.087~ & 1.984~ & 51~ & 284.62~ & 5~ & 3~ & 14.97~ & 456.76~ \\
 & 20 & 14.0233 & 528 & 57.148 & 101 & 7.77 & 28 & 11 & 35 & 11.1143 & 2.810 & 64 & 3.69 & 6 & 4 & 6.52 & 17.98 \\
5000 (3) & 277~ & 27.4761~ & 1229~ & 50.820~ & 89~ & 540.26~ & 41~ & 34~ & 34~ & 14.7647~ & 30.834~ & 85~ & 14.34~ & 7~ & 6~ & 97.56~ & 652.16~ \\
 & 16 & 10.5932 & 183 & 50.820 & 90 & 4.16 & 36 & 32 & 32 & 13.6875 & 32.993 & 127 & 2.57 & 27 & 14 & 12.67 & 19.40 \\
\end{tabular}
}
\label{tab:concat-rp-3-r}
\MyEndConcatTable
}

{
% This table must be reduced in height.
\renewcommand{\MyArrayStretch}{0.85}
\MyBeginConcatTable
\caption{FST Concatenation: Rectilinear Random Instances (continued)}
{\tiny
\InitializeTableFootnoteFlags
\setlength{\tabcolsep}{\MyTabColSep}%
\renewcommand{\arraystretch}{\MyArrayStretch}%
\begin{tabular}{r|r|r|r|r|r|r|r|r|r|r|r|r|r|r|r|r|r}
\cline{2-17}
\multicolumn{1}{c}{Size} &
\multicolumn{6}{|c}{Subtour Relaxation} &
\multicolumn{7}{|c}{Root Node} &
\multicolumn{3}{|c}{Branching} &
\multicolumn{1}{|c}{Total} \\
\cline{2-17}
\multicolumn{1}{c}{(Instance)} &
\multicolumn{1}{|c}{Iter} &
\multicolumn{1}{|c}{NZ/Col} &
\multicolumn{1}{|c}{CC} &
\multicolumn{1}{|c}{Gap} &
\multicolumn{1}{|c}{Frac} &
\multicolumn{1}{|c}{Time} &
\multicolumn{1}{|c}{Iter} &
\multicolumn{1}{|c}{Rnds} &
\multicolumn{1}{|c}{Cuts} &
\multicolumn{1}{|c}{NZPC} &
\multicolumn{1}{|c}{Gap} &
\multicolumn{1}{|c}{Frac} &
\multicolumn{1}{|c}{Time} &
\multicolumn{1}{|c}{Nodes} &
\multicolumn{1}{|c}{Proc} &
\multicolumn{1}{|c}{Time} &
\multicolumn{1}{|c}{Time} \\
\hline
7000 (1) & 1343~ & 55.9813~ & 2523~ & 31.760~ & 236~ & 38266.91~ & 108~ & 41~ & 61~ & 16.1148~ & 5.975~ & 194~ & 1690.58~ & 52~ & 28~ & 310.62~ & 40268.11~ \\
 & 39 & 19.5118 & 828 & 31.760 & 208 & 30.41 & 62 & 51 & 67 & 15.4478 & 6.187 & 169 & 17.93 & 11 & 6 & 41.67 & 90.01 \\
7000 (2) & 65~ & 13.0198~ & 844~ & 79.512~ & 217~ & 25.26~ & 57~ & 14~ & 42~ & 11.4524~ & 17.579~ & 164~ & 8.26~ & 43~ & 24~ & 110.45~ & 143.97~ \\
 & 28 & 18.1777 & 966 & 79.512 & 192 & 19.47 & 25 & 14 & 50 & 13.28 & 15.808 & 132 & 9.49 & 9 & 5 & 33.91 & 62.87 \\
7000 (3) & $ > 2525 $~ & 72.9893~ & $ > 2735 $~ & 17.647~ & 2166~ & $ > 86496.30 $~ & $ > 0 $~ & 0~ & 0~ & 0~ & 17.647~ & 2166~ & $ > 0.00 $~ & $ > 1 $~ & $ > 1 $~ & $ > 0.00 $~ & $ > 86496.30 $~ \\
 & 29 & 12.0283 & 631 & 15.459 & 128 & 12.99 & 0 & 0 & 0 & 0 & 15.459 & 128 & 0.00 & 9 & 6 & 9.85 & 22.84 \\
\end{tabular}
}
\label{tab:concat-rp-4-r}
\MyEndConcatTable
}

\MyBeginConcatTable
\caption{FST Concatenation: Rectilinear Random Instances (continued,
  new code only)}
{\tiny
\InitializeTableFootnoteFlags
\setlength{\tabcolsep}{\MyTabColSep}%
\renewcommand{\arraystretch}{\MyArrayStretch}%
\begin{tabular}{r|r|r|r|r|r|r|r|r|r|r|r|r|r|r|r|r|r}
\cline{2-17}
\multicolumn{1}{c}{Size} &
\multicolumn{6}{|c}{Subtour Relaxation} &
\multicolumn{7}{|c}{Root Node} &
\multicolumn{3}{|c}{Branching} &
\multicolumn{1}{|c}{Total} \\
\cline{2-17}
\multicolumn{1}{c}{(Instance)} &
\multicolumn{1}{|c}{Iter} &
\multicolumn{1}{|c}{NZ/Col} &
\multicolumn{1}{|c}{CC} &
\multicolumn{1}{|c}{Gap} &
\multicolumn{1}{|c}{Frac} &
\multicolumn{1}{|c}{Time} &
\multicolumn{1}{|c}{Iter} &
\multicolumn{1}{|c}{Rnds} &
\multicolumn{1}{|c}{Cuts} &
\multicolumn{1}{|c}{NZPC} &
\multicolumn{1}{|c}{Gap} &
\multicolumn{1}{|c}{Frac} &
\multicolumn{1}{|c}{Time} &
\multicolumn{1}{|c}{Nodes} &
\multicolumn{1}{|c}{Proc} &
\multicolumn{1}{|c}{Time} &
\multicolumn{1}{|c}{Time} \\
\hline
9000 (1) & 24 & 12.3961 & 868 & 89.471 & 290 & 18.02 & 26 & 14 & 84 & 8.57143 & 11.300 & 183 & 9.41 & 4 & 3 & 30.97 & 58.40 \\
9000 (2) & 32 & 19.7643 & 1224 & 56.188 & 215 & 40.75 & 46 & 34 & 79 & 10.4937 & 1.251 & 122 & 23.68 & 1 & 1 & 36.89 & 101.32 \\
9000 (3) & 22 & 11.7654 & 339 & 78.275 & 315 & 12.37 & 25 & 17 & 55 & 8.74545 & 10.917 & 172 & 6.60 & 5 & 3 & 42.86 & 61.83 \\
9000 (4) & 25 & 14.3108 & 811 & 58.719 & 267 & 20.62 & 7 & 4 & 22 & 9 & 43.542 & 221 & 4.97 & 67 & 45 & 165.70 & 191.29 \\
9000 (5) & 27 & 12.5738 & 779 & 106.637 & 287 & 25.17 & 17 & 11 & 59 & 8.40678 & 28.033 & 200 & 5.92 & 17 & 11 & 29.76 & 60.85 \\
9000 (6) & 28 & 17.4819 & 769 & 134.459 & 298 & 30.21 & 27 & 17 & 63 & 6.19048 & 75.177 & 222 & 9.89 & 212 & 175 & 524.20 & 564.30 \\
\end{tabular}
}
\label{tab:concat-rp-5-r}
\MyEndConcatTable

\FloatBarrier

{
% This table must be reduced in height.
\renewcommand{\MyArrayStretch}{0.85}
\MyBeginConcatTable
\caption{FST Concatenation: Hexagonal Random Instances}
{\tiny
\InitializeTableFootnoteFlags
\setlength{\tabcolsep}{\MyTabColSep}%
\renewcommand{\arraystretch}{\MyArrayStretch}%
\begin{tabular}{r|r|r|r|r|r|r|r|r|r|r|r|r|r|r|r|r|r}
\cline{2-17}
\multicolumn{1}{c}{Size} &
\multicolumn{6}{|c}{Subtour Relaxation} &
\multicolumn{7}{|c}{Root Node} &
\multicolumn{3}{|c}{Branching} &
\multicolumn{1}{|c}{Total} \\
\cline{2-17}
\multicolumn{1}{c}{(Instance)} &
\multicolumn{1}{|c}{Iter} &
\multicolumn{1}{|c}{NZ/Col} &
\multicolumn{1}{|c}{CC} &
\multicolumn{1}{|c}{Gap} &
\multicolumn{1}{|c}{Frac} &
\multicolumn{1}{|c}{Time} &
\multicolumn{1}{|c}{Iter} &
\multicolumn{1}{|c}{Rnds} &
\multicolumn{1}{|c}{Cuts} &
\multicolumn{1}{|c}{NZPC} &
\multicolumn{1}{|c}{Gap} &
\multicolumn{1}{|c}{Frac} &
\multicolumn{1}{|c}{Time} &
\multicolumn{1}{|c}{Nodes} &
\multicolumn{1}{|c}{Proc} &
\multicolumn{1}{|c}{Time} &
\multicolumn{1}{|c}{Time} \\
\hline
1000 (1) & 11~ & 6.55263~ & 63~ & 0~ & 0~ & 0.13~ & 0~ & 0~ & 0~ & 0~ & 0~ & 0~ & 0.00~ & 1~ & 1~ & 0.00~ & 0.13~ \\
 & 9 & 6.84017 & 63 & 0 & 0 & 0.11 & 0 & 0 & 0 & 0 & 0 & 0 & 0.00 & 1 & 1 & 0.00 & 0.11 \\
1000 (2) & 329~ & 14.9025~ & 90~ & 0~ & 0~ & 1.67~ & 0~ & 0~ & 0~ & 0~ & 0~ & 0~ & 0.00~ & 1~ & 1~ & 0.00~ & 1.67~ \\
 & 17 & 7.15256 & 16 & 0 & 0 & 0.13 & 0 & 0 & 0 & 0 & 0 & 0 & 0.00 & 1 & 1 & 0.00 & 0.13 \\
1000 (3) & 17~ & 6.53276~ & 28~ & 0~ & 0~ & 0.16~ & 0~ & 0~ & 0~ & 0~ & 0~ & 0~ & 0.00~ & 1~ & 1~ & 0.00~ & 0.16~ \\
 & 9 & 6.05374 & 38 & 0 & 0 & 0.09 & 0 & 0 & 0 & 0 & 0 & 0 & 0.00 & 1 & 1 & 0.00 & 0.09 \\
\end{tabular}
}
\label{tab:concat-rp-1-u3}
\MyEndConcatTable
}

{
% This table must be reduced in height.
\renewcommand{\MyArrayStretch}{0.85}
\MyBeginConcatTable
\caption{FST Concatenation: Hexagonal Random Instances (continued)}
{\tiny
\InitializeTableFootnoteFlags
\setlength{\tabcolsep}{\MyTabColSep}%
\renewcommand{\arraystretch}{\MyArrayStretch}%
\begin{tabular}{r|r|r|r|r|r|r|r|r|r|r|r|r|r|r|r|r|r}
\cline{2-17}
\multicolumn{1}{c}{Size} &
\multicolumn{6}{|c}{Subtour Relaxation} &
\multicolumn{7}{|c}{Root Node} &
\multicolumn{3}{|c}{Branching} &
\multicolumn{1}{|c}{Total} \\
\cline{2-17}
\multicolumn{1}{c}{(Instance)} &
\multicolumn{1}{|c}{Iter} &
\multicolumn{1}{|c}{NZ/Col} &
\multicolumn{1}{|c}{CC} &
\multicolumn{1}{|c}{Gap} &
\multicolumn{1}{|c}{Frac} &
\multicolumn{1}{|c}{Time} &
\multicolumn{1}{|c}{Iter} &
\multicolumn{1}{|c}{Rnds} &
\multicolumn{1}{|c}{Cuts} &
\multicolumn{1}{|c}{NZPC} &
\multicolumn{1}{|c}{Gap} &
\multicolumn{1}{|c}{Frac} &
\multicolumn{1}{|c}{Time} &
\multicolumn{1}{|c}{Nodes} &
\multicolumn{1}{|c}{Proc} &
\multicolumn{1}{|c}{Time} &
\multicolumn{1}{|c}{Time} \\
\hline
3000 (1) & 40~ & 10.7235~ & 125~ & 2.689~ & 13~ & 1.50~ & 70~ & 70~ & 70~ & 10.4~ & 0~ & 18~ & 8.82~ & 1~ & 1~ & 0.00~ & 10.32~ \\
 & 17 & 8.71793 & 144 & 2.689 & 12 & 0.79 & 12 & 12 & 12 & 7.58333 & 2.689 & 15 & 0.05 & 2 & 2 & 0.56 & 1.40 \\
3000 (2) & 332~ & 16.681~ & 409~ & 9.351~ & 9~ & 13.89~ & 292~ & 269~ & 269~ & 20.9108~ & 0~ & 0~ & 15.94~ & 1~ & 1~ & 0.00~ & 29.83~ \\
 & 18 & 7.98093 & 117 & 9.351 & 9 & 0.93 & 10 & 10 & 10 & 11.5 & 9.351 & 13 & 0.00 & 2 & 2 & 0.53 & 1.46 \\
3000 (3) & 16~ & 8.00201~ & 73~ & 0~ & 0~ & 0.80~ & 0~ & 0~ & 0~ & 0~ & 0~ & 0~ & 0.00~ & 1~ & 1~ & 0.00~ & 0.80~ \\
 & 16 & 7.40664 & 84 & 0 & 7 & 0.78 & 0 & 0 & 0 & 0 & 0 & 7 & 0.00 & 1 & 1 & 0.00 & 0.78 \\
\end{tabular}
}
\label{tab:concat-rp-2-u3}
\MyEndConcatTable
}

{
% This table must be reduced in height.
\renewcommand{\MyArrayStretch}{0.85}
\MyBeginConcatTable
\caption{FST Concatenation: Hexagonal Random Instances (continued)}
{\tiny
\InitializeTableFootnoteFlags
\setlength{\tabcolsep}{\MyTabColSep}%
\renewcommand{\arraystretch}{\MyArrayStretch}%
\begin{tabular}{r|r|r|r|r|r|r|r|r|r|r|r|r|r|r|r|r|r}
\cline{2-17}
\multicolumn{1}{c}{Size} &
\multicolumn{6}{|c}{Subtour Relaxation} &
\multicolumn{7}{|c}{Root Node} &
\multicolumn{3}{|c}{Branching} &
\multicolumn{1}{|c}{Total} \\
\cline{2-17}
\multicolumn{1}{c}{(Instance)} &
\multicolumn{1}{|c}{Iter} &
\multicolumn{1}{|c}{NZ/Col} &
\multicolumn{1}{|c}{CC} &
\multicolumn{1}{|c}{Gap} &
\multicolumn{1}{|c}{Frac} &
\multicolumn{1}{|c}{Time} &
\multicolumn{1}{|c}{Iter} &
\multicolumn{1}{|c}{Rnds} &
\multicolumn{1}{|c}{Cuts} &
\multicolumn{1}{|c}{NZPC} &
\multicolumn{1}{|c}{Gap} &
\multicolumn{1}{|c}{Frac} &
\multicolumn{1}{|c}{Time} &
\multicolumn{1}{|c}{Nodes} &
\multicolumn{1}{|c}{Proc} &
\multicolumn{1}{|c}{Time} &
\multicolumn{1}{|c}{Time} \\
\hline
5000 (1) & 351~ & 14.8427~ & 234~ & 0~ & 0~ & 17.57~ & 0~ & 0~ & 0~ & 0~ & 0~ & 0~ & 0.00~ & 1~ & 1~ & 0.00~ & 17.57~ \\
 & 22 & 10.3174 & 126 & 0 & 0 & 2.26 & 0 & 0 & 0 & 0 & 0 & 0 & 0.00 & 1 & 1 & 0.00 & 2.26 \\
5000 (2) & $ > 23701 $~ & 4607.98~ & $ > 932 $~ & 7.485~ & 456~ & $ > 86415.15 $~ & $ > 0 $~ & 0~ & 0~ & 0~ & 7.485~ & 456~ & $ > 0.00 $~ & $ > 1 $~ & $ > 1 $~ & $ > 0.00 $~ & $ > 86415.15 $~ \\
 & 22 & 9.54686 & 204 & 6.886 & 18 & 2.23 & 11 & 11 & 11 & 6.54545 & 6.886 & 23 & 0.06 & 4 & 3 & 1.95 & 4.24 \\
5000 (3) & 75~ & 12.2326~ & 154~ & 21.612~ & 15~ & 4.14~ & 18~ & 17~ & 19~ & 6.63158~ & 0~ & 0~ & 0.99~ & 1~ & 1~ & 0.00~ & 5.13~ \\
 & 13 & 8.43488 & 208 & 21.612 & 15 & 1.59 & 11 & 11 & 12 & 6.66667 & 0 & 0 & 0.64 & 1 & 1 & 0.00 & 2.23 \\
\end{tabular}
}
\label{tab:concat-rp-3-u3}
\MyEndConcatTable
}

{
% This table must be reduced in height.
\renewcommand{\MyArrayStretch}{0.85}
\MyBeginConcatTable
\caption{FST Concatenation: Hexagonal Random Instances (continued)}
{\tiny
\InitializeTableFootnoteFlags
\setlength{\tabcolsep}{\MyTabColSep}%
\renewcommand{\arraystretch}{\MyArrayStretch}%
\begin{tabular}{r|r|r|r|r|r|r|r|r|r|r|r|r|r|r|r|r|r}
\cline{2-17}
\multicolumn{1}{c}{Size} &
\multicolumn{6}{|c}{Subtour Relaxation} &
\multicolumn{7}{|c}{Root Node} &
\multicolumn{3}{|c}{Branching} &
\multicolumn{1}{|c}{Total} \\
\cline{2-17}
\multicolumn{1}{c}{(Instance)} &
\multicolumn{1}{|c}{Iter} &
\multicolumn{1}{|c}{NZ/Col} &
\multicolumn{1}{|c}{CC} &
\multicolumn{1}{|c}{Gap} &
\multicolumn{1}{|c}{Frac} &
\multicolumn{1}{|c}{Time} &
\multicolumn{1}{|c}{Iter} &
\multicolumn{1}{|c}{Rnds} &
\multicolumn{1}{|c}{Cuts} &
\multicolumn{1}{|c}{NZPC} &
\multicolumn{1}{|c}{Gap} &
\multicolumn{1}{|c}{Frac} &
\multicolumn{1}{|c}{Time} &
\multicolumn{1}{|c}{Nodes} &
\multicolumn{1}{|c}{Proc} &
\multicolumn{1}{|c}{Time} &
\multicolumn{1}{|c}{Time} \\
\hline
7000 (1) & 92~ & 20.4616~ & 301~ & 9.135~ & 15~ & 17.17~ & 17~ & 2~ & 2~ & 3~ & 0~ & 0~ & 2.12~ & 1~ & 1~ & 0.00~ & 19.29\UNPRUNED{}~ \\
 & 19 & 15.4626 & 361 & 9.135 & 19 & 7.56 & 6 & 2 & 2 & 3 & 0 & 0 & 1.22 & 1 & 1 & 0.00 & 8.78\UNPRUNED{} \\
7000 (2) & 60~ & 14.187~ & 277~ & 0~ & 10~ & 8.99~ & 3~ & 3~ & 3~ & 3.33333~ & 0~ & 0~ & 0.33~ & 1~ & 1~ & 0.00~ & 9.32~ \\
 & 22 & 10.9319 & 319 & 0 & 4 & 4.68 & 1 & 1 & 1 & 3 & 0 & 0 & 0.10 & 1 & 1 & 0.00 & 4.78 \\
7000 (3) & $ > 52555 $~ & 2187.99~ & $ > 160 $~ & 12.747~ & 52~ & $ > 86412.86 $~ & $ > 0 $~ & 0~ & 0~ & 0~ & 12.747~ & 52~ & $ > 0.00 $~ & $ > 1 $~ & $ > 1 $~ & $ > 0.00 $~ & $ > 86412.86 $~ \\
 & 14 & 8.37911 & 83 & 12.336 & 21 & 2.56 & 18 & 18 & 42 & 21.2857 & 4.793 & 28 & 4.39 & 5 & 4 & 3.56 & 10.51 \\
\end{tabular}
}
\label{tab:concat-rp-4-u3}
\MyEndConcatTable
}

{
% This table must be reduced in height and width
\setlength{\MyTabColSep}{2.75pt}
\renewcommand{\MyArrayStretch}{0.4}
\MyBeginConcatTable
\caption{FST Concatenation: Hexagonal Random Instances (continued)}
{\tiny
\InitializeTableFootnoteFlags
\setlength{\tabcolsep}{\MyTabColSep}%
\renewcommand{\arraystretch}{\MyArrayStretch}%
\begin{tabular}{r|r|r|r|r|r|r|r|r|r|r|r|r|r|r|r|r|r}
\cline{2-17}
\multicolumn{1}{c}{Size} &
\multicolumn{6}{|c}{Subtour Relaxation} &
\multicolumn{7}{|c}{Root Node} &
\multicolumn{3}{|c}{Branching} &
\multicolumn{1}{|c}{Total} \\
\cline{2-17}
\multicolumn{1}{c}{(Instance)} &
\multicolumn{1}{|c}{Iter} &
\multicolumn{1}{|c}{NZ/Col} &
\multicolumn{1}{|c}{CC} &
\multicolumn{1}{|c}{Gap} &
\multicolumn{1}{|c}{Frac} &
\multicolumn{1}{|c}{Time} &
\multicolumn{1}{|c}{Iter} &
\multicolumn{1}{|c}{Rnds} &
\multicolumn{1}{|c}{Cuts} &
\multicolumn{1}{|c}{NZPC} &
\multicolumn{1}{|c}{Gap} &
\multicolumn{1}{|c}{Frac} &
\multicolumn{1}{|c}{Time} &
\multicolumn{1}{|c}{Nodes} &
\multicolumn{1}{|c}{Proc} &
\multicolumn{1}{|c}{Time} &
\multicolumn{1}{|c}{Time} \\
\hline
9000 (1) & 117~ & 15.2591~ & 342~ & 0~ & 0~ & 20.01~ & 0~ & 0~ & 0~ & 0~ & 0~ & 0~ & 0.00~ & 1~ & 1~ & 0.00~ & 20.01~ \\
 & 35 & 13.4054 & 300 & 0 & 0 & 8.94 & 0 & 0 & 0 & 0 & 0 & 0 & 0.00 & 1 & 1 & 0.00 & 8.94 \\
9000 (2) & 446~ & 21.4202~ & 370~ & 15.376~ & 18~ & 58.24~ & $ > 44542 $~ & 44379~ & 44439~ & 63.0989~ & 4.308~ & 30~ & $ > 86343.67 $~ & $ > 1 $~ & $ > 1 $~ & $ > 0.00 $~ & $ > 86401.91 $~ \\
 & 37 & 14.1714 & 314 & 15.376 & 21 & 9.09 & 25 & 25 & 58 & 19.6897 & 12.004 & 30 & 3.37 & 5 & 5 & 8.34 & 20.80 \\
9000 (3) & 26~ & 8.70583~ & 101~ & 7.967~ & 25~ & 6.45~ & 2806~ & 2792~ & 3248~ & 42.3436~ & 1.590~ & 24~ & 37513.17~ & 1~ & 1~ & 313.16~ & 37832.78~ \\
 & 18 & 9.01288 & 120 & 7.967 & 24 & 4.73 & 74 & 73 & 209 & 20.8565 & 2.318 & 33 & 14.59 & 3 & 3 & 7.80 & 27.12 \\
\end{tabular}
}
\label{tab:concat-rp-5-u3}
\MyEndConcatTable
}

\FloatBarrier

{
% This table must be reduced in height.
\renewcommand{\MyArrayStretch}{0.85}
\MyBeginConcatTable
\caption{FST Concatenation: Octilinear Random Instances}
{\tiny
\InitializeTableFootnoteFlags
\setlength{\tabcolsep}{\MyTabColSep}%
\renewcommand{\arraystretch}{\MyArrayStretch}%
\begin{tabular}{r|r|r|r|r|r|r|r|r|r|r|r|r|r|r|r|r|r}
\cline{2-17}
\multicolumn{1}{c}{Size} &
\multicolumn{6}{|c}{Subtour Relaxation} &
\multicolumn{7}{|c}{Root Node} &
\multicolumn{3}{|c}{Branching} &
\multicolumn{1}{|c}{Total} \\
\cline{2-17}
\multicolumn{1}{c}{(Instance)} &
\multicolumn{1}{|c}{Iter} &
\multicolumn{1}{|c}{NZ/Col} &
\multicolumn{1}{|c}{CC} &
\multicolumn{1}{|c}{Gap} &
\multicolumn{1}{|c}{Frac} &
\multicolumn{1}{|c}{Time} &
\multicolumn{1}{|c}{Iter} &
\multicolumn{1}{|c}{Rnds} &
\multicolumn{1}{|c}{Cuts} &
\multicolumn{1}{|c}{NZPC} &
\multicolumn{1}{|c}{Gap} &
\multicolumn{1}{|c}{Frac} &
\multicolumn{1}{|c}{Time} &
\multicolumn{1}{|c}{Nodes} &
\multicolumn{1}{|c}{Proc} &
\multicolumn{1}{|c}{Time} &
\multicolumn{1}{|c}{Time} \\
\hline
1000 (1) & 17~ & 3.82482~ & 33~ & 0~ & 0~ & 0.08~ & 0~ & 0~ & 0~ & 0~ & 0~ & 0~ & 0.00~ & 1~ & 1~ & 0.00~ & 0.08~ \\
 & 9 & 4.97554 & 27 & 0 & 0 & 0.07 & 0 & 0 & 0 & 0 & 0 & 0 & 0.00 & 1 & 1 & 0.00 & 0.07 \\
1000 (2) & 6~ & 3.94706~ & 33~ & 0~ & 0~ & 0.06~ & 0~ & 0~ & 0~ & 0~ & 0~ & 0~ & 0.00~ & 1~ & 1~ & 0.00~ & 0.06~ \\
 & 6 & 4.54984 & 55 & 0 & 0 & 0.05 & 0 & 0 & 0 & 0 & 0 & 0 & 0.00 & 1 & 1 & 0.00 & 0.05 \\
1000 (3) & 6~ & 3.36257~ & 37~ & 0~ & 0~ & 0.05~ & 0~ & 0~ & 0~ & 0~ & 0~ & 0~ & 0.00~ & 1~ & 1~ & 0.00~ & 0.05~ \\
 & 6 & 3.85321 & 37 & 0 & 0 & 0.04 & 0 & 0 & 0 & 0 & 0 & 0 & 0.00 & 1 & 1 & 0.00 & 0.04 \\
\end{tabular}
}
\label{tab:concat-rp-1-u4}
\MyEndConcatTable
}

\FloatBarrier

{
% This table must be reduced in height.
\renewcommand{\MyArrayStretch}{0.85}
\MyBeginConcatTable
\caption{FST Concatenation: Octilinear Random Instances (continued)}
{\tiny
\InitializeTableFootnoteFlags
\setlength{\tabcolsep}{\MyTabColSep}%
\renewcommand{\arraystretch}{\MyArrayStretch}%
\begin{tabular}{r|r|r|r|r|r|r|r|r|r|r|r|r|r|r|r|r|r}
\cline{2-17}
\multicolumn{1}{c}{Size} &
\multicolumn{6}{|c}{Subtour Relaxation} &
\multicolumn{7}{|c}{Root Node} &
\multicolumn{3}{|c}{Branching} &
\multicolumn{1}{|c}{Total} \\
\cline{2-17}
\multicolumn{1}{c}{(Instance)} &
\multicolumn{1}{|c}{Iter} &
\multicolumn{1}{|c}{NZ/Col} &
\multicolumn{1}{|c}{CC} &
\multicolumn{1}{|c}{Gap} &
\multicolumn{1}{|c}{Frac} &
\multicolumn{1}{|c}{Time} &
\multicolumn{1}{|c}{Iter} &
\multicolumn{1}{|c}{Rnds} &
\multicolumn{1}{|c}{Cuts} &
\multicolumn{1}{|c}{NZPC} &
\multicolumn{1}{|c}{Gap} &
\multicolumn{1}{|c}{Frac} &
\multicolumn{1}{|c}{Time} &
\multicolumn{1}{|c}{Nodes} &
\multicolumn{1}{|c}{Proc} &
\multicolumn{1}{|c}{Time} &
\multicolumn{1}{|c}{Time} \\
\hline
3000 (1) & 110~ & 5.86097~ & 75~ & 3.526~ & 15~ & 1.85~ & 4~ & 4~ & 7~ & 21.7143~ & 0~ & 0~ & 0.27~ & 1~ & 1~ & 0.00~ & 2.12~ \\
 & 11 & 5.48 & 72 & 3.526 & 15 & 0.32 & 4 & 4 & 7 & 8.85714 & 0 & 0 & 0.09 & 1 & 1 & 0.00 & 0.41 \\
3000 (2) & 23~ & 5.78807~ & 86~ & 14.734~ & 4~ & 0.69~ & 8~ & 6~ & 6~ & 11.5~ & 0~ & 0~ & 0.23~ & 1~ & 1~ & 0.00~ & 0.92~ \\
 & 14 & 6.46043 & 70 & 14.734 & 4 & 0.44 & 7 & 6 & 6 & 10.8333 & 0 & 0 & 0.18 & 1 & 1 & 0.00 & 0.62 \\
3000 (3) & 15~ & 5.54731~ & 65~ & 15.792~ & 4~ & 0.39~ & 25~ & 12~ & 12~ & 21.25~ & 0~ & 0~ & 1.46~ & 1~ & 1~ & 0.00~ & 1.85~ \\
 & 10 & 5.82215 & 73 & 15.792 & 4 & 0.38 & 14 & 12 & 12 & 16.3333 & 0 & 0 & 0.45 & 1 & 1 & 0.00 & 0.83 \\
\end{tabular}
}
\label{tab:concat-rp-2-u4}
\MyEndConcatTable
}

\FloatBarrier

{
% This table must be reduced in height.
\renewcommand{\MyArrayStretch}{0.85}
\MyBeginConcatTable
\caption{FST Concatenation: Octilinear Random Instances (continued)}
{\tiny
\InitializeTableFootnoteFlags
\setlength{\tabcolsep}{\MyTabColSep}%
\renewcommand{\arraystretch}{\MyArrayStretch}%
\begin{tabular}{r|r|r|r|r|r|r|r|r|r|r|r|r|r|r|r|r|r}
\cline{2-17}
\multicolumn{1}{c}{Size} &
\multicolumn{6}{|c}{Subtour Relaxation} &
\multicolumn{7}{|c}{Root Node} &
\multicolumn{3}{|c}{Branching} &
\multicolumn{1}{|c}{Total} \\
\cline{2-17}
\multicolumn{1}{c}{(Instance)} &
\multicolumn{1}{|c}{Iter} &
\multicolumn{1}{|c}{NZ/Col} &
\multicolumn{1}{|c}{CC} &
\multicolumn{1}{|c}{Gap} &
\multicolumn{1}{|c}{Frac} &
\multicolumn{1}{|c}{Time} &
\multicolumn{1}{|c}{Iter} &
\multicolumn{1}{|c}{Rnds} &
\multicolumn{1}{|c}{Cuts} &
\multicolumn{1}{|c}{NZPC} &
\multicolumn{1}{|c}{Gap} &
\multicolumn{1}{|c}{Frac} &
\multicolumn{1}{|c}{Time} &
\multicolumn{1}{|c}{Nodes} &
\multicolumn{1}{|c}{Proc} &
\multicolumn{1}{|c}{Time} &
\multicolumn{1}{|c}{Time} \\
\hline
5000 (1) & 100~ & 7.40848~ & 308~ & 4.681~ & 10~ & 4.25~ & 1~ & 1~ & 2~ & 4.5~ & 0~ & 0~ & 0.06~ & 1~ & 1~ & 0.00~ & 4.31~ \\
 & 19 & 7.73459 & 213 & 4.681 & 10 & 1.40 & 1 & 1 & 2 & 4 & 0 & 0 & 0.06 & 1 & 1 & 0.00 & 1.46 \\
5000 (2) & 51~ & 6.33049~ & 150~ & 7.617~ & 12~ & 2.25~ & 5~ & 5~ & 6~ & 6.5~ & 0~ & 0~ & 0.24~ & 1~ & 1~ & 0.00~ & 2.49~ \\
 & 15 & 7.00644 & 266 & 7.617 & 12 & 1.04 & 5 & 5 & 6 & 6.5 & 0 & 0 & 0.24 & 1 & 1 & 0.00 & 1.28 \\
5000 (3) & 11~ & 4.60892~ & 75~ & 4.190~ & 14~ & 0.83~ & 3~ & 3~ & 4~ & 8~ & 0~ & 0~ & 0.15~ & 1~ & 1~ & 0.00~ & 0.98~ \\
 & 12 & 5.72446 & 70 & 4.190 & 15 & 0.87 & 3 & 3 & 4 & 6.5 & 0 & 0 & 0.15 & 1 & 1 & 0.00 & 1.02 \\
\end{tabular}
}
\label{tab:concat-rp-3-u4}
\MyEndConcatTable
}

\FloatBarrier

{
% This table must be reduced in height.
\renewcommand{\MyArrayStretch}{0.85}
\MyBeginConcatTable
\caption{FST Concatenation: Octilinear Random Instances (continued)}
{\tiny
\InitializeTableFootnoteFlags
\setlength{\tabcolsep}{\MyTabColSep}%
\renewcommand{\arraystretch}{\MyArrayStretch}%
\begin{tabular}{r|r|r|r|r|r|r|r|r|r|r|r|r|r|r|r|r|r}
\cline{2-17}
\multicolumn{1}{c}{Size} &
\multicolumn{6}{|c}{Subtour Relaxation} &
\multicolumn{7}{|c}{Root Node} &
\multicolumn{3}{|c}{Branching} &
\multicolumn{1}{|c}{Total} \\
\cline{2-17}
\multicolumn{1}{c}{(Instance)} &
\multicolumn{1}{|c}{Iter} &
\multicolumn{1}{|c}{NZ/Col} &
\multicolumn{1}{|c}{CC} &
\multicolumn{1}{|c}{Gap} &
\multicolumn{1}{|c}{Frac} &
\multicolumn{1}{|c}{Time} &
\multicolumn{1}{|c}{Iter} &
\multicolumn{1}{|c}{Rnds} &
\multicolumn{1}{|c}{Cuts} &
\multicolumn{1}{|c}{NZPC} &
\multicolumn{1}{|c}{Gap} &
\multicolumn{1}{|c}{Frac} &
\multicolumn{1}{|c}{Time} &
\multicolumn{1}{|c}{Nodes} &
\multicolumn{1}{|c}{Proc} &
\multicolumn{1}{|c}{Time} &
\multicolumn{1}{|c}{Time} \\
\hline
7000 (1) & 39~ & 6.50383~ & 136~ & .392~ & 10~ & 3.45~ & 3~ & 3~ & 3~ & 19~ & 0~ & 0~ & 0.28~ & 1~ & 1~ & 0.00~ & 3.73~ \\
 & 16 & 6.7628 & 157 & .392 & 10 & 1.90 & 3 & 3 & 3 & 8.66667 & 0 & 0 & 0.27 & 1 & 1 & 0.00 & 2.17 \\
7000 (2) & 19~ & 6.2104~ & 260~ & 7.409~ & 27~ & 1.90~ & 7~ & 6~ & 10~ & 6.3~ & 0~ & 0~ & 0.57~ & 1~ & 1~ & 0.00~ & 2.47~ \\
 & 24 & 9.06898 & 157 & 7.409 & 27 & 2.55 & 9 & 8 & 12 & 6.08333 & 0 & 0 & 0.78 & 1 & 1 & 0.00 & 3.33 \\
7000 (3) & 25~ & 5.46993~ & 98~ & 0~ & 0~ & 2.52~ & 0~ & 0~ & 0~ & 0~ & 0~ & 0~ & 0.00~ & 1~ & 1~ & 0.00~ & 2.52~ \\
 & 12 & 5.55574 & 112 & 0 & 0 & 1.56 & 0 & 0 & 0 & 0 & 0 & 0 & 0.00 & 1 & 1 & 0.00 & 1.56 \\
\end{tabular}
}
\label{tab:concat-rp-4-u4}
\MyEndConcatTable
}

{
% This table must be reduced in height.
\renewcommand{\MyArrayStretch}{0.85}
\MyBeginConcatTable
\caption{FST Concatenation: Octilinear Random Instances (continued)}
{\tiny
\InitializeTableFootnoteFlags
\setlength{\tabcolsep}{\MyTabColSep}%
\renewcommand{\arraystretch}{\MyArrayStretch}%
\begin{tabular}{r|r|r|r|r|r|r|r|r|r|r|r|r|r|r|r|r|r}
\cline{2-17}
\multicolumn{1}{c}{Size} &
\multicolumn{6}{|c}{Subtour Relaxation} &
\multicolumn{7}{|c}{Root Node} &
\multicolumn{3}{|c}{Branching} &
\multicolumn{1}{|c}{Total} \\
\cline{2-17}
\multicolumn{1}{c}{(Instance)} &
\multicolumn{1}{|c}{Iter} &
\multicolumn{1}{|c}{NZ/Col} &
\multicolumn{1}{|c}{CC} &
\multicolumn{1}{|c}{Gap} &
\multicolumn{1}{|c}{Frac} &
\multicolumn{1}{|c}{Time} &
\multicolumn{1}{|c}{Iter} &
\multicolumn{1}{|c}{Rnds} &
\multicolumn{1}{|c}{Cuts} &
\multicolumn{1}{|c}{NZPC} &
\multicolumn{1}{|c}{Gap} &
\multicolumn{1}{|c}{Frac} &
\multicolumn{1}{|c}{Time} &
\multicolumn{1}{|c}{Nodes} &
\multicolumn{1}{|c}{Proc} &
\multicolumn{1}{|c}{Time} &
\multicolumn{1}{|c}{Time} \\
\hline
9000 (1) & 235~ & 11.7208~ & 534~ & 14.344~ & 27~ & 33.92~ & 206~ & 4~ & 5~ & 6.6~ & 3.497~ & 19~ & 21.86~ & 1~ & 1~ & 0.60~ & 56.38~ \\
 & 38 & 13.3942 & 253 & 14.344 & 27 & 6.18 & 5 & 4 & 5 & 5 & 3.497 & 19 & 0.57 & 2 & 2 & 1.16 & 7.91 \\
9000 (2) & 222~ & 12.4784~ & 725~ & 1.834~ & 5~ & 64.71~ & 2~ & 2~ & 2~ & 5.5~ & 0~ & 0~ & 0.25~ & 1~ & 1~ & 0.00~ & 64.96~ \\
 & 17 & 7.8632 & 194 & 1.834 & 5 & 3.13 & 2 & 2 & 2 & 5.5 & 0 & 0 & 0.22 & 1 & 1 & 0.00 & 3.35 \\
9000 (3) & 523~ & 12.1164~ & 619~ & 4.103~ & 12~ & 69.11~ & 59~ & 2~ & 3~ & 8.66667~ & 0~ & 277~ & 16.66~ & 1~ & 1~ & 0.00~ & 85.77~ \\
 & 15 & 5.9222 & 99 & 4.103 & 12 & 2.20 & 2 & 2 & 3 & 7.66667 & 0 & 0 & 0.23 & 1 & 1 & 0.00 & 2.43 \\
\end{tabular}
}
\label{tab:concat-rp-5-u4}
\MyEndConcatTable
}

\FloatBarrier

{
% This table must be reduced in height.
\renewcommand{\MyArrayStretch}{0.85}
\MyBeginConcatTable
\caption{FST Concatenation: Large Euclidean Random Instances}
{\tiny
\InitializeTableFootnoteFlags
\setlength{\tabcolsep}{\MyTabColSep}%
\renewcommand{\arraystretch}{\MyArrayStretch}%
\begin{tabular}{r|r|r|r|r|r|r|r|r|r|r|r|r|r|r|r|r|r}
\cline{2-17}
\multicolumn{1}{c}{Size} &
\multicolumn{6}{|c}{Subtour Relaxation} &
\multicolumn{7}{|c}{Root Node} &
\multicolumn{3}{|c}{Branching} &
\multicolumn{1}{|c}{Total} \\
\cline{2-17}
\multicolumn{1}{c}{(Instance)} &
\multicolumn{1}{|c}{Iter} &
\multicolumn{1}{|c}{NZ/Col} &
\multicolumn{1}{|c}{CC} &
\multicolumn{1}{|c}{Gap} &
\multicolumn{1}{|c}{Frac} &
\multicolumn{1}{|c}{Time} &
\multicolumn{1}{|c}{Iter} &
\multicolumn{1}{|c}{Rnds} &
\multicolumn{1}{|c}{Cuts} &
\multicolumn{1}{|c}{NZPC} &
\multicolumn{1}{|c}{Gap} &
\multicolumn{1}{|c}{Frac} &
\multicolumn{1}{|c}{Time} &
\multicolumn{1}{|c}{Nodes} &
\multicolumn{1}{|c}{Proc} &
\multicolumn{1}{|c}{Time} &
\multicolumn{1}{|c}{Time} \\
\hline
15,000 (1) & 58~ & 6.00472~ & 258~ & 2.608~ & 17~ & 15.60~ & 21~ & 10~ & 20~ & 34.65~ & 0~ & 0~ & 39.46~ & 1~ & 1~ & 0.00~ & 55.06~ \\
 & 26 & 8.22534 & 178 & 2.608 & 17 & 9.34 & 13 & 12 & 23 & 15.9565 & 0 & 0 & 7.28 & 1 & 1 & 0.00 & 16.62 \\
15,000 (2) & 147~ & 9.14707~ & 630~ & .114~ & 4~ & 63.14~ & 1~ & 1~ & 1~ & 13~ & 0~ & 4~ & 0.34~ & 1~ & 1~ & 0.00~ & 63.48~ \\
 & 14 & 8.76321 & 315 & .114 & 7 & 5.15 & 1 & 1 & 1 & 7 & 0 & 7 & 0.38 & 1 & 1 & 0.00 & 5.53 \\
15,000 (3) & 90~ & 7.6671~ & 393~ & 1.666~ & 24~ & 25.89~ & 2~ & 2~ & 7~ & 5.85714~ & 0~ & 0~ & 0.66~ & 1~ & 1~ & 0.00~ & 26.55~ \\
 & 19 & 6.5117 & 260 & 1.666 & 24 & 6.62 & 2 & 2 & 7 & 4.71429 & 0 & 0 & 0.70 & 1 & 1 & 0.00 & 7.32 \\
\end{tabular}
}
\label{tab:concat-huge-1-e}
\MyEndConcatTable
}

\FloatBarrier

{
% This table must be reduced in height.
\renewcommand{\MyArrayStretch}{0.85}
\MyBeginConcatTable
\caption{FST Concatenation: Large Euclidean Random Instances (continued)}
{\tiny
\InitializeTableFootnoteFlags
\setlength{\tabcolsep}{\MyTabColSep}%
\renewcommand{\arraystretch}{\MyArrayStretch}%
\begin{tabular}{r|r|r|r|r|r|r|r|r|r|r|r|r|r|r|r|r|r}
\cline{2-17}
\multicolumn{1}{c}{Size} &
\multicolumn{6}{|c}{Subtour Relaxation} &
\multicolumn{7}{|c}{Root Node} &
\multicolumn{3}{|c}{Branching} &
\multicolumn{1}{|c}{Total} \\
\cline{2-17}
\multicolumn{1}{c}{(Instance)} &
\multicolumn{1}{|c}{Iter} &
\multicolumn{1}{|c}{NZ/Col} &
\multicolumn{1}{|c}{CC} &
\multicolumn{1}{|c}{Gap} &
\multicolumn{1}{|c}{Frac} &
\multicolumn{1}{|c}{Time} &
\multicolumn{1}{|c}{Iter} &
\multicolumn{1}{|c}{Rnds} &
\multicolumn{1}{|c}{Cuts} &
\multicolumn{1}{|c}{NZPC} &
\multicolumn{1}{|c}{Gap} &
\multicolumn{1}{|c}{Frac} &
\multicolumn{1}{|c}{Time} &
\multicolumn{1}{|c}{Nodes} &
\multicolumn{1}{|c}{Proc} &
\multicolumn{1}{|c}{Time} &
\multicolumn{1}{|c}{Time} \\
\hline
25,000 (1) & 192~ & 6.95488~ & 613~ & 1.026~ & 14~ & 132.09~ & 5~ & 2~ & 4~ & 8.25~ & .010~ & 3~ & 5.07~ & 1~ & 1~ & 2.75~ & 139.91~ \\
 & 31 & 9.04825 & 343 & 1.026 & 14 & 25.39 & 2 & 2 & 5 & 15.6 & 0 & 0 & 2.06 & 1 & 1 & 0.00 & 27.45 \\
25,000 (2) & 63~ & 7.08368~ & 220~ & .312~ & 4~ & 43.63~ & 9~ & 5~ & 5~ & 11.2~ & 0~ & 19~ & 7.30~ & 1~ & 1~ & 0.00~ & 50.93~ \\
 & 19 & 8.54578 & 261 & .312 & 4 & 14.42 & 4 & 4 & 4 & 10.25 & 0 & 0 & 3.39 & 1 & 1 & 0.00 & 17.81 \\
25,000 (3) & 94~ & 7.60495~ & 464~ & .615~ & 34~ & 64.86~ & 3~ & 2~ & 7~ & 7.28571~ & 0~ & 0~ & 2.17~ & 1~ & 1~ & 0.00~ & 67.03~ \\
 & 32 & 10.7617 & 326 & .615 & 34 & 24.70 & 2 & 2 & 7 & 4.71429 & 0 & 7 & 1.59 & 1 & 1 & 0.00 & 26.29 \\
\end{tabular}
}
\label{tab:concat-huge-2-e}
\MyEndConcatTable
}

\FloatBarrier

{
% This table must be reduced in height.
\renewcommand{\MyArrayStretch}{0.85}
\MyBeginConcatTable
\caption{FST Concatenation: Large Euclidean Random Instances (continued)}
\input{Table_A_033.tex}
\label{tab:concat-huge-3-e}
\MyEndConcatTable
}

\FloatBarrier

\MyBeginConcatTable
\caption{FST Concatenation: Large Euclidean Random Instances (continued)}
\input{Table_A_034.tex}
\label{tab:concat-huge-4-e}
\MyEndConcatTable

\FloatBarrier

\MyBeginConcatTable
\caption{FST Concatenation: Large Euclidean Random Instances
  (continued, new code only)}
\input{Table_A_035.tex}
\label{tab:concat-huge-5-e}
\MyEndConcatTable

\FloatBarrier

{
% This table must be reduced in height.
\renewcommand{\MyArrayStretch}{0.85}
\MyBeginConcatTable
\caption{FST Concatenation: Large Rectilinear Random Instances (new
  code only)}
\input{Table_A_036.tex}
\label{tab:concat-huge-1-r}
\MyEndConcatTable
}

\FloatBarrier

\MyBeginConcatTable
\caption{FST Concatenation: Euclidean TSPLIB Instances}
\input{Table_A_037.tex}
\label{tab:concat-tsplib-1-e}
\MyEndConcatTable

{
% This table must be made narrower
\setlength{\MyTabColSep}{2.3pt}
\MyBeginConcatTable
\caption{FST Concatenation: Euclidean TSPLIB Instances (continued)}
\input{Table_A_038.tex}
\label{tab:concat-tsplib-2-e}
\MyEndConcatTable
}

\FloatBarrier

\MyBeginConcatTable
\caption{FST Concatenation: Rectilinear TSPLIB Instances}
\input{Table_A_039.tex}
\label{tab:concat-tsplib-1-r}
\MyEndConcatTable

\MyBeginConcatTable
\caption{FST Concatenation: Rectilinear TSPLIB Instances (continued)}
\input{Table_A_040.tex}
\label{tab:concat-tsplib-2-r}
\MyEndConcatTable

\FloatBarrier

\MyBeginConcatTable
\caption{FST Concatenation: Hexagonal TSPLIB Instances}
\input{Table_A_041.tex}
\label{tab:concat-tsplib-1-u3}
\MyEndConcatTable

{
% This table must be made narrower
\setlength{\MyTabColSep}{1.75pt}
\MyBeginConcatTable
\caption{FST Concatenation: Hexagonal TSPLIB Instances (continued)}
\input{Table_A_042.tex}
\label{tab:concat-tsplib-2-u3}
\MyEndConcatTable
}

\FloatBarrier

\MyBeginConcatTable
\caption{FST Concatenation: Octilinear TSPLIB Instances}
\input{Table_A_043.tex}
\label{tab:concat-tsplib-1-u4}
\MyEndConcatTable

{
% This table must be made narrower
\setlength{\MyTabColSep}{2.5pt}
\MyBeginConcatTable
\caption{FST Concatenation: Octilinear TSPLIB Instances (continued)}
\input{Table_A_044.tex}
\label{tab:concat-tsplib-2-u4}
\MyEndConcatTable
}

\FloatBarrier

\MyBeginFSTTable
\caption{FST Generation and Pruning: Euclidean OR-Library Instances}
\input{Table_A_045.tex}
\label{tab:fsts-orlib-e}
\MyEndFSTTable

\MyBeginFSTTable
\caption{FST Generation and Pruning: Rectilinear OR-Library Instances}
\input{Table_A_046.tex}
\label{tab:fsts-orlib-r}
\MyEndFSTTable

\MyBeginFSTTable
\caption{FST Generation and Pruning: Hexagonal OR-Library Instances}
\input{Table_A_047.tex}
\label{tab:fsts-orlib-u3}
\MyEndFSTTable

\MyBeginFSTTable
\caption{FST Generation and Pruning: Octilinear OR-Library Instances}
\input{Table_A_048.tex}
\label{tab:fsts-orlib-u4}
\MyEndFSTTable

\FloatBarrier

\MyBeginFSTTable
\caption{FST Generation and Pruning: Random Euclidean Instances}
\input{Table_A_049.tex}
\label{tab:fsts-rp-1-e}
\MyEndFSTTable

\MyBeginFSTTable
\caption{FST Generation and Pruning: Random Euclidean Instances (continued)}
\input{Table_A_050.tex}
\label{tab:fsts-rp-2-e}
\MyEndFSTTable

\MyBeginFSTTable
\caption{FST Generation and Pruning: Random Euclidean Instances (continued)}
\input{Table_A_051.tex}
\label{tab:fsts-rp-3-e}
\MyEndFSTTable

\MyBeginFSTTable
\caption{FST Generation and Pruning: Random Euclidean Instances (continued)}
\input{Table_A_052.tex}
\label{tab:fsts-rp-4-e}
\MyEndFSTTable

\FloatBarrier

\MyBeginFSTTable
\caption{FST Generation and Pruning: Random Rectilinear Instances}
\input{Table_A_053.tex}
\label{tab:fsts-rp-1-r}
\MyEndFSTTable

\MyBeginFSTTable
\caption{FST Generation and Pruning: Random Rectilinear Instances (continued)}
\input{Table_A_054.tex}
\label{tab:fsts-rp-2-r}
\MyEndFSTTable

\MyBeginFSTTable
\caption{FST Generation and Pruning: Random Rectilinear Instances (continued)}
\input{Table_A_055.tex}
\label{tab:fsts-rp-3-r}
\MyEndFSTTable

\MyBeginFSTTable
\caption{FST Generation and Pruning: Random Rectilinear Instances (continued)}
\input{Table_A_056.tex}
\label{tab:fsts-rp-4-r}
\MyEndFSTTable

\FloatBarrier

\MyBeginFSTTable
\caption{FST Generation and Pruning: Random Hexagonal Instances}
\input{Table_A_057.tex}
\label{tab:fsts-rp-1-u3}
\MyEndFSTTable

\MyBeginFSTTable
\caption{FST Generation and Pruning: Random Hexagonal Instances (continued)}
\input{Table_A_058.tex}
\label{tab:fsts-rp-2-u3}
\MyEndFSTTable

\MyBeginFSTTable
\caption{FST Generation and Pruning: Random Hexagonal Instances (continued)}
\input{Table_A_059.tex}
\label{tab:fsts-rp-3-u3}
\MyEndFSTTable

\MyBeginFSTTable
\caption{FST Generation and Pruning: Random Hexagonal Instances (continued)}
\input{Table_A_060.tex}
\label{tab:fsts-rp-4-u3}
\MyEndFSTTable

\FloatBarrier

\MyBeginFSTTable
\caption{FST Generation and Pruning: Random Octilinear Instances}
\input{Table_A_061.tex}
\label{tab:fsts-rp-1-u4}
\MyEndFSTTable

\MyBeginFSTTable
\caption{FST Generation and Pruning: Random Octilinear Instances (continued)}
\input{Table_A_062.tex}
\label{tab:fsts-rp-2-u4}
\MyEndFSTTable

\MyBeginFSTTable
\caption{FST Generation and Pruning: Random Octilinear Instances (continued)}
\input{Table_A_063.tex}
\label{tab:fsts-rp-3-u4}
\MyEndFSTTable

\MyBeginFSTTable
\caption{FST Generation and Pruning: Random Octilinear Instances (continued)}
\input{Table_A_064.tex}
\label{tab:fsts-rp-4-u4}
\MyEndFSTTable

\FloatBarrier

\MyBeginFSTTable
\caption{FST Generation and Pruning: Large Random Euclidean Instances}
\input{Table_A_065.tex}
\label{tab:fsts-huge-1-e}
\MyEndFSTTable

\MyBeginFSTTable
\caption{FST Generation and Pruning: Large Random Euclidean Instances (continued)}
\input{Table_A_066.tex}
\label{tab:fsts-huge-2-e}
\MyEndFSTTable

\MyBeginFSTTable
\caption{FST Generation and Pruning: Large Random Euclidean Instances (continued)}
\input{Table_A_067.tex}
\label{tab:fsts-huge-3-e}
\MyEndFSTTable

\FloatBarrier

\MyBeginFSTTable
\caption{FST Generation and Pruning: Large Random Rectilinear Instances}
\input{Table_A_068.tex}
\label{tab:fsts-huge-1-r}
\MyEndFSTTable

\MyBeginFSTTable
\caption{FST Generation and Pruning: Large Random Rectilinear Instances (continued)}
\input{Table_A_069.tex}
\label{tab:fsts-huge-2-r}
\MyEndFSTTable

\FloatBarrier

\MyBeginFSTTable
\caption{FST Generation and Pruning: Euclidean TSPLIB Instances}
\input{Table_A_070.tex}
\label{tab:fsts-tsplib-e}
\MyEndFSTTable

\FloatBarrier

\MyBeginFSTTable
\caption{FST Generation and Pruning: Rectilinear TSPLIB Instances}
\input{Table_A_071.tex}
\label{tab:fsts-tsplib-r}
\MyEndFSTTable

\FloatBarrier

\MyBeginFSTTable
\caption{FST Generation and Pruning: Hexagonal TSPLIB Instances}
\input{Table_A_072.tex}
\label{tab:fsts-tsplib-u3}
\MyEndFSTTable

\FloatBarrier

\MyBeginFSTTable
\caption{FST Generation and Pruning: Octilinear TSPLIB Instances}
\input{Table_A_073.tex}
\label{tab:fsts-tsplib-u4}
\MyEndFSTTable

\eject

\strut
\vskip 0.25in
\noindent
{\huge \bf Revision History}

\vskip 0.25in

\noindent
03/21/2024: Original version.

\vskip 0.25in

\noindent
07/08/2025: Version 2:
\begin{itemize}
  \item Improvements from review of condensed journal version.
  \item Added computational results for 250k and 500k terminal random
    Euclidean instances \# 2--15.
  \item Added \BothMode{} mode results.
  \item Added Figures~\ref{fig:cdf-2k-both-large-small-subtours},
    \ref{fig:histogram-2k-both-large-small-subtours} and accompanying
    text for $n=2000$ terminal instances to show trend as $n$
    increases.
\end{itemize}

\end{document}